\documentclass[final,3p,times]{Classes/elsarticle}
\usepackage{graphicx}
\usepackage{amssymb}
\usepackage{amsthm}
\usepackage{amsmath}
\usepackage{tabularx}
\usepackage{multirow}
\usepackage{booktabs}
\usepackage{microtype} 
\usepackage[parfill]{parskip}
\usepackage{ulem}
\usepackage{cancel}
\usepackage{caption}
\usepackage{subcaption}
\usepackage{verbatim}
\usepackage[table,xcdraw]{xcolor}
\usepackage{enumitem}
\usepackage{nomencl}
\usepackage[capitalise]{cleveref}

\journal{Journal of Theoretical, Computational and Applied Mechanics}

\makenomenclature
\begin{document}

%%%%BEGIN FRONT MATTER
\begin{frontmatter}
% Title
\title{Mathematics of Stable Tensegrity Structures}

% List of authors
\author[1,2]{Ajay B. Harish \corref{one}}
\ead{ajaybangalore.harish@manchester.ac.uk}
\author[1]{Shubham Deshpande$^{\dagger}$}
%\ead{shubhamd201097@gmail.com}
\author[3]{Stephanie R. Andress$^{\dagger}$}
%\ead{stevie.andress@gmail.com}
\cortext[one]{Corresponding author at: Department of Mechanical, Aerospace and Civil Engineering, University of Manchester, UK; $^{\dagger}$ Equal contributors}

% Addresses
\address[1]{Department of Civil and Environmental Engineering, University of California, Berkeley, USA}
\address[2]{School of Mechanical, Aerospace and Civil Engineering (MACE), University of Manchester, UK}
\address[3]{Department of Mechanical Engineering, Purdue University, USA}

% Abstract to be restricted to 100 words
\begin{abstract}
Tensegrity structures have been extensively studied over the last years due to their potential applications in modern engineering like metamaterials, deployable structures, planetary lander modules, etc. Many of the form-finding methods proposed continue to produce structures with one or more soft/swinging modes. These modes have been vividly highlighted and outlined as the grounds for these structures to be unsuitable as engineering structures. This work proposes a relationship between the number of rods and strings to satisfy the full-rank convexity criterion as a part of the form-finding process. Using the proposed form-finding process for the famous three-rod tensegrity, the work proposes an alternative three-rod ten-string that is stable. The work demonstrates that the stable tensegrities suitable for engineering are feasible and can be designed.
\end{abstract}

% Keywords of the paper
\begin{keyword}
Tensegrity \sep Linearisation stability \sep Form-finding \sep Modal analysis \sep Space structure; Deployable structure
\end{keyword}

\end{frontmatter}

% Nomenclature
%\textbf{Nomenclature}
%\mbox{}
\nomenclature{$\mathbf{A}$}{Small motion prescribed on the tensegrity, particularly representing the end points of the rods, points of beads and pivot points related to the tensegrity}
\nomenclature{$\left[A \right]_{(9\times 3)}$}{Matrix representing the rigid-body rotations}
\nomenclature{$\left[B \right]_{(9\times 3)}$}{Matrix representing the rigid-body translations}
\nomenclature{$\left[C \right]_{(9\times 6)}$}{Matrix representing the standard basis function representing the motion of points $p_j$ and $p_k$}
\nomenclature{$b$}{Number of beads in the tensegrity}
\nomenclature{$c_{ij}$}{Length of the cable}
\nomenclature{$C$}{List of unordered pairs ${i,j}$ describing the list of cables}
\nomenclature{$d_{ij}$}{Distance between any two points ${i,j}$}
\nomenclature{$f_k$}{Force acting on a member of the tensegrity}
\nomenclature{$f^{k} (\mathbf{v})$}{Functional representation for the $c-e$-inequalities, i.e. cable constraints}
\nomenclature{$g^{\ell} (\mathbf{v})$}{Functional representation for the $e$-equalities, i.e. rod constraints}
\nomenclature{$\ell_k$}{Length of the member of the tensegrity}
\nomenclature{$m$}{Number of cables in the tensegrity}
\nomenclature{$n$}{Number of rods in the tensegrity}
\nomenclature{$N$}{Number of points making the tensegrity}
\nomenclature{$p_i$}{Position vector of a single point, represented by $(x_i,y_i,z_i)$, of the tensegrity}
\nomenclature{$P$}{Position of the tensegrity, represented by $3N$-tuples of $(x_i,y_i,z_i)$}
\nomenclature{$q_k$}{Force density in a member of the tensegrity}
\nomenclature{$Q$}{Position of the tensegrity after movement, represented by $3N$-tuples of $(x_i,y_i,z_i)$}
\nomenclature{$r_{ij}$}{Length of the rod}
\nomenclature{$\mathbf{r}_i$}{Rigid-body motion vectors with $i = 0, \cdots, 5$}
\nomenclature{$R$}{List of unordered pairs ${i,j}$ describing the list of rods}
\nomenclature{$\mathbb{R}^{3N}$}{Space spanning all the $3N$ points of the tensegrity}
\nomenclature{$\mathcal{R}$}{Sub-space of all the motions of the tensegrity}
\nomenclature{$T$}{Translation vector applied on the tensegrity}
\nomenclature{$\mathbf{v}_j$}{The vector representing each connection of the tensegrity and also referred to as constraint vectors. Here, $j =6, \cdots, 3N-1$}
\nomenclature{$\mathcal{V}$}{Sub-space of all constraint vectors}
\nomenclature{$(x_i,y_i,z_i)$}{The coordinate of each point in the tensegrity}
\nomenclature{$\varepsilon$}{Small motion prescribed on the tensegrity}
\nomenclature{$\Pi$}{The set of all points in the tensegrity}
\nomenclature{$S$}{A plane that constraints the motion of a bead}
\nomenclature{$\mathcal{S}$}{The space of all the string index pairs $\{i,j\}$}
\nomenclature{$c_r$}{Magnitude of compressive force acting on the rod in the tensegrity}
\nomenclature{$t_{ij}$}{Magnitude of tensile force acting on the string in the tensegrity}
\nomenclature{$\hat{\mathbf{v}}_{i,j}$}{Unit vector depicting the direction along which the force acts, i.e from point $p_j$ to $p_i$}
\nomenclature{$\mathbf{t}_{ij}$}{Tensile force vector acting on the wire}
\nomenclature{$\mathbf{r}_{(\cdot)}$}{Compressive force vector acting on the rod}
\printnomenclature

% Section 1 (Introduction)
\section{Introduction}

Tensegrities go back to several decades earlier when \citet{Sn1973} developed some of his first structures. In 1961, \citet{Fu1961} presented a class of cable-bar structures where bars were arranged in compression, structural integrity maintained by strings, known as Tensegrities. In addition, tensegrity structures are defined to have not two rods joined at the same point. Today, people use both broader and narrower definitions of `tensegrity.' However, in this work, we will consider the definition of \citet{Fu1961} and \citet{RoWh1996}. The structures developed by Snelson and co-workers led to a significant renewal in interest related to tensegrities. However, like the famous ``Snelson Tower'', many of these are unstable, making them harder to use and integrate into actual engineering structures. This work demonstrates that this is due to these designs developed based on soft elastic strings rather than rigid strings. Over the years, researchers have focused attention to mathematically understand the source of these instabilities and develop rigid and stable structures for various applications.

In the last years, there has been a renewed interest in the area of tensegrities, including ideas and discussions in the areas of tensegrity metamaterials \cite{fraternali2012, Lietal2019, wang2020b}, bio-tensegrities \cite{chai2020a}, planetary lander modules \cite{garanger2020}, tensegrity-inspired structures \cite{luo2017a, Paetal2019, y2020a, ma2020a, mirzaaghazadeh2020a}, foldable structures \cite{YaSu2019}. Of particular interest to this work are those that discuss the mechanics of tensegrities and form-finding methods and shall be highlighted. While the two topics appear diverse, it is necessary to develop a thorough comprehension of these structures' mechanics to develop robust form-finding methodologies. 

\subsection{Literature review}
General discussion and detailed literature reviews on many of the topics of interest to the tensegrity community have been discussed in \citet{TiPe2003, oliveira2009a, Muetal2019, wang2021a}. Some common conclusions common to many of the studies into the mechanics of the tensegrity structures include:
\begin{enumerate}
\item The structures are statically indeterminate and dynamically quasi-stable
\item Most of the structures have a soft mode, that is also often referred to as infinitesimal flex or mechanism mode or swinging mode
\item Such soft modes are reduced by the addition of pre-stress that can also be equivalent to requiring an additional string
\item There is a need to understand the source of non-linearity / snapping behavior/instability observed in the force-displacement curve.
\item There is a need for a solution that can do away with these instabilities to help design stable, rigid tensegrities that can be used in actual engineering applications.
\end{enumerate}

\subsection{Mechanics of tensegrities}
The 1980s and 90s saw tremendous growth in the mathematical literature surrounding tensegrity structures' mechanics and stability. This included several works from Connelly \cite{Co1982, CoWh1992, Co1995}, Whiteley \cite{Wh1981, RoWh1981, Wh1982, Wh1984, CoWh1996, RoWh1996}, Calladine and Pelligrino \cite{Ca1978, PeCa1986, Pe1992, KwPe1994} and several other mathematicians and structural engineers. In \citet{CoWh1992} present four mathematical concepts to uniquely define the stability of a tensegrity, i.e. infinitesimally rigid $\subset$ pre-stress stability $\subset$ second-order rigidity $\subset$ rigid. This rigidity hierarchy has formed the basis for several works to this day, including the present work where a second-order rigidity is presented.

Alongside the seminal preliminary works of Connelly and co-workers, several recent works from the early 2000s re-visit the mechanics of tensegrities. \citet{OpWi2000} state that tensegrities are under-constrained and undergo an infinitesimal flex and display nonlinear geometric stiffening upon loading. Such a definition automatically paves the way for unstable and soft modes in tensegrities, as this work, later, demonstrates. They discuss instability in the force-displacement relations by considering the famous 3-rod and 9-string tensegrity. This is further supported by a vibration analysis \cite{OpWi2001} where they demonstrate insufficient damping due to the \texttt{elastic} cables used to build the tensegrity structure. \citet{OpWi2001} conclude that this inefficient damping that is apparent in the models is a serious drawback for practical usage of tensegrity structures. The same three-rod structure is further again considered for dynamic analysis by \citet{Mu2001}. They unequivocally demonstrate the soft mode in this structure, which is mentioned as an \texttt{infinitesimal mechanism mode of swinging bars}. They further demonstrate that the pre-stressing leads to an increase in the eigen-frequency of the lowest mode of oscillation. As will be discussed later in this work, this pre-stressing acts as the additional string that would have been required to build a stable three-rod and ten-string tensegrity. This non-linearity in the force-displacement behavior of such a tensegrity is again explored and discussed as a snapping instability by \citet{zhang2016a}. Here, they explore the snapping behavior as a function of the applied torque and the structure's symmetry. Such a snap behavior is an obvious effect of the soft modes that are inherent in the structure. Other prominent works that discussed statics and dynamics of tensegrities include \citet{muni2001a, oliveto2011a, ashwear2014a, li2016a, hsu2020a}.

\subsubsection{Tensegrity towers}
Another area that has received attention includes deployable structures and tensegrity towers. \citet{sultan2003a} further use the tower, made of three-rod nine-string units and proposed by \citet{Mu2001}, for the design of a deployable structure. The work starts with the pre-conceived notion that the tensegrities are composed of soft and hard members. While the work explores possible deployable paths, such would work in the absence of soft modes. Juan and Mirats Tur \cite{juan2008a, mirats2009a}, motivated by robotics and controls, discuss the static and dynamic analysis of tensegrity structures. They use the ideas of Connelly and co-workers and provide a comprehensive review of the static analysis of such structures to date. \citet{bel2011a} propose an algorithm to consider deployable structures. These deployable structures are made of unit cells of tensegrities and continuous cables that can be adjusted to actively deploy and activate structures. However, while they discuss the active deformation behavior, local soft modes that arise in such structures are not considered, limiting the serious applicability for tensegrity engineering. Some of the other works that consider the behavior of tensegrity towers include \citet{luo2017a, oh2019a, kahla2020a, li2020a, y2020a}. The concept of deployable tensegrity is re-visited by \citet{liu2017a} where they propose the usage of soft elastomers for the programmable deployment of these tensegrities. The work again considers the famous three-rod and nine-string tensegrity. The proposed work's potential applications are immense but again does not consider the swinging mode instability in these tensegrities repeatedly discussed in the mechanics' literature. However, exploring such deployment for more rigid tensegrities of the kind proposed in this work can help make deployable tensegrity-based stable space structures. 

\subsubsection{Design philosophies}
Specific designs based on the hypothesis of stability arising due to symmetry have been another direction of common interest explored in many works on tensegrity structures. \citet{MuNi2001,MuNi2001b} considered the regular truncated icosahedral and dodecahedral shaped tensegrities. The authors again confirm the presence of 15 distinct infinitesimal mechanism modes in these structures. The same idea is used by \citet{rimoli2017a} to construct tensegrity structures based on truncated octahedron. However, the work is motivated towards building lattice-based composite structures that show superior impact behavior. They use a homogenization technique to obtain the dynamic response of a structure and pitch its potential applications in metamaterials. While the averaged temporal behavior during impact shows promise, the oscillations in the results indicate several soft modes that warrant significant further investigations.

\subsubsection{Form-finding techniques}
While mechanics of tensegrities has been one of the focus areas, the other area remains form-finding. Various form-finding techniques have been proposed over the years. In general, the form-finding techniques can be classified based on (a) force method (b) energy method. \citet{barnes1999a} discussed the application of dynamic relaxation technique for tensegrity nets and membranes. The idea here is the application of an explicit solution technique for the static behavior of structures. In the last two decades, this has been used extensively. A thorough review of the form-finding techniques till the early 2000s can be found in the work of \cite{TiPe2003}. \cite{lu2017a} consider the idea of dynamic relaxation for pre-stressed strut structures and determine equilibrium configurations for three and four-rod tensegrities. The minimization starts with a pre-set initial configuration to obtain a relatively near-final configuration. 

The force density methods are based on the formulation of a constrained optimization problem to minimize the overall force density. Here, force density, $q_k$, is defined as 
\begin{equation}
q_k = \frac{f_k}{\ell_k}
\end{equation}
where for any member $k$ has a member force of $f_k$ and a length of $\ell_k$. \citet{lee2016a} consider the constrained optimization problem, aiming to minimize the overall force density in the cable structures. Similar to \citet{lee2016a}, \citet{cai2018a} also consider a minimization through grouping of elements. Alternatively, \citet{XuWaLu2018} propose nonlinear programming as a possible methodology for minimization of the force density in the cables. \citet{liu2019a} consider ideas from computer graphics to create zones of non-intersection and propose a force maximization method of form-finding. The work demonstrates tensegrities with a varying number of rods and strings. While they denote the structures as super-stable, the work does not show any static/dynamic/vibrational analysis to prove the same. Some of the other recent works using form-finding based on force density idea include those of \citet{zhang2020a, wang2020a, koohestani2020a,wang2021a}. \citet{uzun2016a} alternative uses the idea of genetic algorithms to perform potential energy minimization for various configurations. There are several more that have explored the concept of form-finding including \citet{rieffel2009a, tran2010b, tran2010a, Ko2013, gan2015a, yuan2017a, feng2017a, dong2019a, aloui2019}.

One of the common aspects of all the above literature in form-finding is a general lack of discussions on the presence or absence of a soft swinging mode discussed thoroughly in mechanics literature, stated earlier. A general overview of the form-finding papers shows that most of the structures obtained from these form-finding methods lack stability and have one or more soft modes. A simple check is the total number of strings in the structure, which will be discussed in the following sections. \Cref{comparison-strings} compares some of the structures produced from the form-finding techniques that show evidence of string-deficiency and definite inherent soft modes. While exceptions could exist, most often, the proposed structures are deficient in one more strings.

\begin{table}[!htb]
\centering
\caption{Tensegrities proposed in literature using form-finding methods. The number of strings proposed are compared with actual number of strings required for linearization stability. The concept of linearization stability will be introduced later in this work.}
\label{comparison-strings}
\begin{tabular}{ccccc}
\hline
Work & Figure & Number of rods & Number of strings & \begin{tabular}[c]{@{}c@{}}Required number of\\ strings for stability \\ $5\left(n-1\right)$ \end{tabular} \\ \hline
 & 11 & 3 & 9 & 10 \\ \cline{2-5} 
\multirow{-2}{*}{\citet{lu2017a}} & 12 & 4 & 12 & 15 \\ \hline
\citet{cai2018a} & 7 & 6 & 24 & 25 \\ \hline
 & 2c & 3 & 9 & 10 \\ \cline{2-5} 
 & 2d & 4 & 12 & 15 \\ \cline{2-5} 
 & 2e & 4 & 11 & 15 \\ \cline{2-5} 
\multirow{-4}{*}{\citet{XuWaLu2018}} & 4e & 6 & 9 & 25 \\ \hline
 &  & 3 & 9 & 10 \\ \cline{3-5} 
 &  & 4 & 12 & 15 \\ \cline{3-5} 
\multirow{-3}{*}{\citet{liu2019a}} & \multirow{-3}{*}{19} & 6 & 15 & 25 \\ \hline
\citet{Paetal2019} & 1 & 12 & 32 & 55 \\ \hline
\citet{zhang2020a} & 5 & 4 & 12 & 15 \\ \hline
\citet{koohestani2020a} & 3 & 20 & 60 & 95 \\ \hline
 & 4 & 6 & 24 & 25 \\ \hline
 & 7 & 30 & 90 & 145 \\ \hline
\citet{wang2021a} & 7 & 3 & 12 & 10 \\ \hline
 & 9 & 6 & 24 & 25 \\ \hline
\end{tabular}
\end{table}

\subsection{Overview and need of this work}
As evident, most of the designed tensegrities employ methods that lead to ``string-deficient'' structures that have inherent soft swinging modes. Such structures lack the rigidity and are not suitable for actual engineering structures. As will be discussed in this work, a possible reasoning for form-finding leading to such string-deficient structures is that many of the designs are based on usage of elastic/rubbery bands as many researchers have pointed out that the tensegrities are harder to build with rigid strings. In contrast, the large elastic behavior of elastic strings can help through small adjustments in length. However, small adjustments in length automatically introduce such swinging modes. While the obvious question is, ``why can we not add more cables to stabilize them''. Unfortunately, such an addition is not feasible, without using more cables than required. The only possible way is to design them from the ground up, considering constraints to eliminate such soft modes. Firstly, this work demonstrates the mathematical framework that integrates the mechanics into the form-finding process to ensure robust tensegrity structural designs.

Nevertheless, in certain applications, the presence of swinging modes are a necessity. A simple example is an application to mimic water plants that sway with the water current; a planetary lander module that can have a soft mode that allows it to deform and diffuse impact energy; meta-materials that block or allow certain frequency modes. However, if the form-finding technique does not account for the mechanics of the swinging modes, they cannot eliminate or control them in preferred directions. Thus, this work stems from the above need for understanding string-deficiency, leading to soft swinging modes, and further incorporating this into the form-finding process to design robust tensegrity structures suitable for any engineering application. 

It is also important to note here that this work tries to eliminate the assumption that tensegrity structures are necessarily pre-stressed. It demonstrates that the notion of pre-stressing arose owing the form-finding and design philosophies that use rubber bands (rather than rigid ropes/strings).

This work demonstrates the mathematics of mechanics and stability involved in designing these structures and derives a relation between the number of rods/beads/cables required to ensure a stable structure. The work also further outlines possible form-finding methods, that integrate the mechanics, to build larger stable tensegrity structures. To demonstrate the developed concepts, the work considers the famous three-rod and nine-string tensegrity to compare with the three-rod and ten-string tensegrity derived from the methods proposed here.
% Section 2 (Mathematics of tensegrities)
\section{Mathematics of tensegrities}

It is important to define some terminology required to identify a tensegrity and the space it spans uniquely. It can be assumed that given some arbitrary $N$ points in space, these points need to be connected by $n$ rods, $b$ beads, and $m$ cables. Considering generality, it can be assumed that it is unknown at this stage as to how the rods and cables need to be placed nor if the selected points are suitable to form a tensegrity. To illustrate this point with an example, a kite is the simplest form of a tensegrity. In this case, there are four points and connected by two rods and four strings. A kite is feasible only if the four points form a quadrilateral. However, if 3/4 points were co-linear, then it is not feasible to form a kite tensegrity. Thus, it is essential to check if the selected $N>0$ points are suitable for forming a tensegrity even before form-finding. In this endeavor, some relevant mathematical concepts are defined. 

The $n > 0$ rods, $b \ge 0$ beads, and $N = 2n+b$ cable attachment points are indexed by $0,1,\cdots, N-1$. Here, the relation between the number of cables $(m)$ / rods $(n)$ / beads $(b)$ is unknown and yet to be determined. Each rod has two attachments points while a bead has only one. The topology of rods are uniquely defined by a list $R$ of $n$ unordered pairs $\{i,j\}$. This typical pair $\{2k,2k + 1\}$ for $k = 0, \cdots,n$ are to be joined by rods. Similarly, a list $C$ defines the $m$ unordered pairs to be joined by cables. At this point, the dependency between $n$ and $m$ remains unknown and shall be discussed in the upcoming section. 

Topology enforces obvious conditions like length $r_{ij}$ for each rod pair and a length $c_{ij}$ for each cable pair. Here the position includes specific location given as $\left(x_i,y_i,z_i\right)$ for each node $i$. A single position $P$ shall be defined to include where the tensegrity is and what shape it is, and this needs to be disentangled to study the shape properly. 

Note here that a single tensegrity can occupy multiple positions due to soft modes like the swinging mode. The set of all these positions are denoted by $\Pi$. To be mathematically precise, $\Pi$ defines the space of all possible positions where each single position $P$ (of a tensegrity) corresponds to $3N$ tuples $\left(x_0, y_0, z_0, \cdots , x_{N-1}, y_{N-1}, z_{N-1}\right)$. In other words, this looks like the space $\mathbb{R}^{3N}$ of $3N$ tuples of real numbers.

It is natural that the position $P$ is allowed if it satisfies, the ``rod constraints''
\begin{equation}
\left(x_i-x_j\right)^2 +\left(x_i-x_j\right)^2 + \left(x_i-x_j\right)^2  = d_{ij} = r_{ij}^{2} \ \forall \left\{i,j\right\} \in R
\label{rodconstraint}
\end{equation}
and the ``cable constraints''
\begin{equation}
\left(x_i-x_j\right)^2 +\left(x_i-x_j\right)^2 + \left(x_i-x_j\right)^2 = d_{ij} \le c_{ij}^{2} \ \forall \left\{i,j\right\} \in C
\label{cableconstraint}
\end{equation}
where $d_{ij}$ is referred to as distance between any given ends $i$ and $j$, $r_{ij}$ is the length of the rod and $c_{ij}$ is the original undeformed length of the cable. In addition, ``anchorage constraints'', or boundary conditions also need to be considered to define the boundary value problem uniquely. Thus,
\begin{equation}
\left( x_i, y_i, z_i\right) = \left( \overline{x}_i, \overline{y}_i, \overline{z}_i\right) 
\end{equation}
are also considered into the mathematical description of the tensegrity.

\subsection{Tensegrity shape coordinates}
Further on, here, to avoid the soft modes and help uniquely characterize a tensegrity, it is essential to define additional mathematical constraints. If the original tensegrity defined by position $P$ is moved to $P'$ due to some combination of translations and rotations, and the only possible positions $P'$ are such that they are not near, i.e., $P \approx P'$ is not true, then $P$ represents a uniquely defined stable tensegrity. In other words, it is necessary to ensure that there are no possible positions $P'$ in the proximity of $P$ that can be reached without appreciable change in energy (i.e., like a swinging mode or soft mode). For example, the famous three-rod and nine-string tensegrity has a swinging mode instability since there are possible configurations near the equilibrium position where it can move into.
 
Given a Position $P = \left\{\left(x_i, y_i, z_i \right) \ \forall i = 0, \cdots , N-1\right\}$ then a translation $T = \left(t_x,t_y,t_z\right)$ moves each point $p_i = \left(x_i, y_i, z_i \right)$, such that the rod and cable constraints are not violated. This is equivalent to a virtual displacement offered to test the stability of the structure at a particular position $P$.

Now assuming that a small motion $\mathbf{A} = \mathbf{I} + \varepsilon \mathbf{B}$ is applied to position $P$. We can write the squared distance, originally $d_{ij}^{2} = \left( p_i - p_j \right) \cdot  \left( p_i - p_j \right)$, as
\begin{equation}
\left( \mathbf{A}p_i -  \mathbf{A}p_j \right) \cdot \left( \mathbf{A}p_i -  \mathbf{A}p_j \right) = d_{ij}^{2} + 0 + \varepsilon^{2} \left( \mathbf{B}v_{ij} \right) \cdot \left( \mathbf{B}v_{ij} \right)
\end{equation}
where $v_{ij} = \left( p_i - p_j \right)$. It is evident that this is constant to the first order in $\varepsilon$ but has higher order terms. Thus, turning the tensegrity moves $P$ in a vector direction that does not change the $d_{ij}$ or their match to the constraints. Like all other separations, rods and cables are changing length at zero rate. At $P = \left(x_0, y_0, z_0, x_1, y_1, z_1, \cdots, x_i, y_i, z_i, \cdots, x_{N-1}, y_{N-1}, z_{N-1} \right)$, all vectors giving zero (rate of) change in shape and are exactly the linear combinations of the rigid body motion vectors given as
\begin{equation}
\begin{split}
\mathbf{r}_{0} &= \left[ 1 \ 0 \ 0 \ 1 \ 0 \ 0 \cdots 1 \ 0 \ 0 \cdots 1 \ 0 \ 0\right]^{T} \\
\mathbf{r}_{1} &= \left[ 0 \ 1 \ 0 \ 0 \ 1 \ 0 \cdots 0 \ 1 \ 0 \cdots 0 \ 1 \ 0\right]^{T} \\
\mathbf{r}_{2} &= \left[ 0 \ 0 \ 1 \ 0 \ 0 \ 1 \cdots 0 \ 0 \ 1 \cdots 0 \ 0 \ 1\right]^{T}
\end{split}
\end{equation}
for translation, and 
\begin{equation}
\begin{split}
\mathbf{r}_{3} &= \left[ 0 \ z_1 \ -y_1 \ 0 \ z_2 \ -y_2 \cdots 0 \ z_i \ -y_i \cdots 0 \ z_{N-1} \ -y_{N-1} \right]^{T} \\
\mathbf{r}_{4} &= \left[ z_1 \ 0 \ -x_1 \ z_1 \ 0 \ -x_2 \cdots z_i \ 0 \ -x_i \cdots z_{N-1} \ 0 \ -x_{N-1} \right]^{T} \\
\mathbf{r}_{5} &= \left[ y_1 \ -x_1 \ 0 \ y_2 \ -x_2 \ 0 \cdots y_i \ -x_i \ 0 \cdots y_{N-1} \ -x_{N-1} \ 0 \right]^{T}
\end{split}
\end{equation}
for rotations. The $3N$-vector that describes the rigid body motion of the $N$ points will necessarily be a unique combination of the above six vectors. 

Further on, one can assume that there exists additional $(3N-6)$ vectors $\mathbf{v}_j$ with $j = 6, ..., 3N-1$, such that the set comprising of $r_0,r_1, \cdots, r_5, v_6, v_7, \cdots,v_{3N-1}$ span the subspace $\mathcal{R}$, where $\mathcal{R}$ is the subspace of all possible motions for $P$. While the vectors $\mathbf{r}_i$ represent the rigid body motions, the vectors $\mathbf{v}_i$ represent the swinging modes or soft modes. A more rigorous definition is available in the upcoming sections. For now, it can be said that $\mathbf{v}_j$ are linearly independent and that any vector in the space $\mathcal{R}$ can be written uniquely as a combination of the $\mathbf{r}_i$ and $\mathbf{v}_j$. This also means that any shape-and-position $Q$ can be reached from $P$ as
\begin{equation}
Q - P = r_0 \mathbf{r}_0 + r_1 \mathbf{r}_1 + \cdots + r_5 \mathbf{r}_5 + v_6 \mathbf{v}_6 + v_7 \mathbf{v}_7 + \cdots + v_{3N-1} \mathbf{v}_{3N-1}
\end{equation}
for some $3N$ tuple $\left(r_0, \cdots, r_5, v_6, \cdots, v_{3N-1} \right)$. The $3N$ tuples can be decomposed into rigid body rotation and translation 
\begin{equation}
\begin{split}
P\left(r_0, \cdots, r_5, v_6, \cdots, v_{3N-1}\right) &= R_{\left(r_0, \cdots, r_5\right)} \left( P + \left( v_6 \mathbf{v}_6 + v_7 \mathbf{v}_7 + \cdots + v_{3N-1} \mathbf{v}_{3N-1}  \right) \right) \\
&= R_{\left(r_0, \cdots, r_5\right)} \left( P + \sum_{i=6}^{N-1} {v_{i} \mathbf{v}_{i}}  \right)
\end{split}
\end{equation}

The operator $R$, for a single point, can be illustrated to be
\begin{equation}
\left(x,y,z\right) \rightarrow \begin{bmatrix}
1 & 0 & 0 \\
0 & \cos r_3 & \sin r_3 \\
0 & -\sin r_3 & \cos r_3 
\end{bmatrix} \left( \begin{bmatrix}
\cos r_4 & 0 & \sin r_4 \\
0 & 1 & 0 \\
-\sin r_4 & 0 & \cos r_4 
\end{bmatrix} \left( \begin{bmatrix}
\cos r_5 & \sin r_5 & 0 \\
-\sin r_5 & \cos r_5 & 0 \\
0 & 0 & 1 
\end{bmatrix} \begin{bmatrix}
x \\
y \\
z 
\end{bmatrix} \right) \right) + \begin{bmatrix}
r_0 \\
r_1 \\
r_2 
\end{bmatrix}
\end{equation}
For the $N$ tuples of points, $3N$ tuples of coordinates, the operator performs the same operation, though only three at a time. If one were to add and write out the entire $3N$ vector, it would lead to an inelegant and large matrix formula with the result dependent on the order of rotation $r_3$ around the $x$-axis, $r_4$ around the $y$-axis, and $r_5$ around the $z$-axis.

Additionally, since it is known that the position $P$ with all $N$ points can itself be written as a vector of size $3N$
\begin{equation}
P = \left(x_0,y_0,z_0,x_1,y_1,z_1, \cdots , x_{N-1},y_{N-1},z_{N-1} \right)
\end{equation}
it can also be expressed as a linear combination of the $3N$ tuple vectors as
\begin{equation}
P = r_0 \mathbf{r}_0 + \cdots + r_5 \mathbf{r}_5 + v_6 \mathbf{v}_6 + \cdots + v_{N-1} \mathbf{v}_{3N-1}.
\end{equation}
This implies that the coefficients $\left(r_0, \cdots, r_5, v_6, \cdots, v_{3N-1} \right)$ can be uniquely determined if the Jacobian matrix resulting from the column vectors, i.e. $\left[ \mathbf{r}_0, \cdots , \mathbf{r}_5, \mathbf{v}_6, \cdots, \mathbf{v}_{3N-1} \right]$ has a rank of $3N$. In such a case, it can also be said that the selected set of points can lead to a stable tensegrity. It is pertinent to note here that the vectors $\mathbf{v}_i$ are yet to be determined and calculated.

\subsection{Constraint vectors}
Elaborating from earlier, this sub-section outlines the determination of the constraint vectors $\left[ \mathbf{v}_6, \mathbf{v}_7,\cdots, \mathbf{v}_{N-1}, \right]$. One reasonable starting point to determine them is to consider the usual basis vectors that span the space $\mathbb{R}^{3N}$ consisting of the $N$ points of the tensegrity. These are
\begin{equation}
\begin{bmatrix} 1 \\ 0 \\ 0 \\ \vdots \\ 0 \\ 0 \\ 0 \end{bmatrix}, \begin{bmatrix} 0 \\ 1 \\ 0 \\ \vdots \\ 0 \\ 0 \\ 0 \end{bmatrix}, \begin{bmatrix} 0 \\ 0 \\ 1 \\ \vdots \\ 0 \\ 0 \\ 0 \end{bmatrix}, \cdots, \begin{bmatrix} 0 \\ 0 \\ 0 \\ \vdots \\ 1 \\ 0 \\ 0 \end{bmatrix}, \begin{bmatrix} 1 \\ 0 \\ 0 \\ \vdots \\ 0 \\ 1 \\ 0 \end{bmatrix}, \begin{bmatrix} 1 \\ 0 \\ 0 \\ \vdots \\ 0 \\ 0 \\ 1 \end{bmatrix}
\end{equation}

Thus, this results in $3N$ vectors and in addition to the six vectors $\left(\mathbf{r}_i\right)$, there are a total of $3N+6$ vectors. This implies that six of them need to be discarded to obtain a final set for the basis vectors $\left(\mathbf{v}_i\right)$. Additionally, the final set of vectors can be set to span a space $\mathcal{V}$ such that $\left(\mathbf{v}_{6}, \mathbf{v}_{7}, \cdots, \mathbf{v}_{N-1} \right) \in \mathcal{V}$. The question remains on the conditionality to choose the six vectors that need to be discarded and yet maintain a rank of $3N$. For example: if one were to discard the first six vectors, this would allow $\mathcal{V}$ to contain any pattern of motion direction for the point $p_i$ that keeps points $p_0$ and $p_1$ fixed. The points are kept fixed by the line $\overline{p_0 p_1}$. Thus, the uniqueness of the configuration $P$ is lost. Similar random omissions can lead to choices that include a non-unique configuration for $P$. 

A more systematic manner to choose the set of six vectors to be discarded could be to find the three points that give the biggest triangle, i.e. to determine the three indices
\begin{equation}
0 \le i < j < k < N \text{ such that } \| \left( p_j - p_i \right) \| \text{ is greatest.}
\end{equation}
If one were to construct a matrix by just restricting to these three points, then
\begin{equation}
\begin{bmatrix}
\left[A \right]_{(9\times 3)} & \left[B \right]_{(9\times 3)} & \left[C \right]_{(9\times 6)}
\end{bmatrix}
\end{equation}
where $\left[A \right]_{(9\times 3)}$ relates to the rigid body rotations along $x-$, $y-$ and $z-$ axis; $\left[B \right]_{(9\times 3)}$ relates to the rigid body translations about $x-$, $y-$ and $z-$ axis; $\left[C \right]_{(9\times 3)}$ relates to the standard basis functions to move points $p_j$ and $p_k$ (three each). The matrices can be elaborated as
\begin{equation}
\left[A\right] = \begin{bmatrix}
1 & 0 & 0 \\
0 & 1 & 0 \\
0 & 0 & 1 \\
1 & 0 & 0 \\
0 & 1 & 0 \\
0 & 0 & 1 \\
1 & 0 & 0 \\
0 & 1 & 0 \\
0 & 0 & 1 
\end{bmatrix};
\left[B\right] = \begin{bmatrix}
0  & z_i   & y_i \\
z_i   & 0  & -x_i\\
-y_i  & x_i   & 0\\
0  & z_j   & y_j \\
z_j   & 0  & -x_j\\
-y_j  & x_j   & 0\\
0  & z_k   & y_k \\
z_k   & 0  & -x_k\\
-y_k  & x_k   & 0
\end{bmatrix}; 
\left[C\right] = \begin{bmatrix}
0 & 0 & 0 & 0 & 0 & 0 \\
0 & 0 & 0 & 0 & 0 & 0 \\
0 & 0 & 0 & 0 & 0 & 0 \\
1 & 0 & 0 & 0 & 0 & 0 \\
0 & 1 & 0 & 0 & 0 & 0 \\
0 & 0 & 1 & 0 & 0 & 0 \\
0 & 0 & 0 & 1 & 0 & 0 \\
0 & 0 & 0 & 0 & 1 & 0 \\
0 & 0 & 0 & 0 & 0 & 1 \\
\end{bmatrix}
\end{equation}

Now, if one were to drop three vectors from the combined matrix, there are a total of 20 possible combinations (choosing 3 out of 6, without regard to order). This requires calculating the determinant of these combinations and choosing one that allows for a full-rank matrix. With these, now one can also determine the the coefficients $\left( v_6, v_7, \cdots, v_{3N-1} \right)$ that allow $P$ to be represented in a unique manner. If there was any motion that leads to $P'$, the above ensures $P \ne P'$. This implies that no soft-modes are present, and the only way to change would be a completely new shape. The coming sections will discuss how the developed ideas could help identify actual rod and cable positions, or otherwise commonly known as form-finding, in the tensegrity literature.

\subsection{Linear constraints and convexity}
At this point, the physical meaning of the vectors $\mathbf{v}$ need to be defined. Considering the constraint equations as a functional representation of the vectors $\mathbf{v}$, they can be given as
\begin{equation}
\begin{split}
\mathbf{g}^{\ell} \left( \mathbf{v} \right) &= 0 \ \forall \ \ell = 0, \cdots, e-1 \ \text{(e equalities, representing rod constraints)} \\
\mathbf{f}^{k} \left( \mathbf{v} \right) &\le 0 \ \forall \ k = e, \cdots, c-1 \ \text{(c-e inequalities, representing cable constraints)}
\end{split}
\end{equation}
Here, the functionals $\mathbf{f}^{k} \left(\mathbf{v}\right) $ are given to be
\begin{equation}
\mathbf{f}^{k} \left(\mathbf{v}\right) = \left[ f_{0}^{k} \ f_{1}^{k} \ \cdots f_{m-1}^{k} \right] \begin{bmatrix}
v_0 \\
v_1 \\
\vdots \\
v_{m-1}
\end{bmatrix} = f_{0}^{k} v_0 + f_{1}^{k} v_1 + \cdots + f_{m-1}^{k} v_{m-1}
\end{equation}
where $f^k$ approximates the change in the length from its value at $P$. The above automatically satisfy the zero vector or otherwise known as the trivial solution. However, the interest is often in the non-trivial solutions where $\mathbf{v} \neq 0$. The vectors $\mathbf{v} \in \mathcal{V}$ represent the small changes from $P$ under the cable constraints $\mathbf{f}^{k} \left(\mathbf{v}\right)$. As discussed earlier, these cable constraints represent the change in length of the string/cables from its value at $P$. These are zero when $\mathbf{v}=0$. This is feasible only if all coefficients are zero or if the vectors are linearly independent.

In the case that ${\mathbf{g}^0, \cdots, \mathbf{g}^{e-1}, \mathbf{f}^{e}, \cdots, \mathbf{g}^{c-1}}$ satisfies the full rank condition and 
\begin{equation}
b_0 \mathbf{g}^0 + \cdots + b_{e-1} \mathbf{g}^{e-1}+ a_{e} \mathbf{f}^{e} + \cdots + a_{c-1}\mathbf{g}^{c-1} = 0
\end{equation}
for some $a_{e} > 0, \cdots, a_{c-1} >0$, then it is said to satisfy the spanning convexity test.
% Section 3 (Beam-based models)
\section{Bead-based models}
The previous section discussed some mathematical constructs to define a tensegrity uniquely. The introduction also discussed that many tensegrities designed to date are based on concepts of using an elastic cable (like an elastic band) rather than actual rigid cables (like the behavior of an actual rope). Before proceeding further, it is imperative to comprehend the difference between the behavior ensuing the usage of elastic vs. rigid cables. Here, simple structures with $b=1$ (not rods but only beads) are considered to illustrate the difference. In a curved space, such as a sphere's surface, the bead can be fixed firmly on the surface. However, in 3D space, it collapses unless anchored at fixed points. While these do not satisfy the conventional definition of ``Tensegrities,'' they provide low-dimensional examples to demonstrate the calculations that can also apply to free-standing stable tensegrities. 

\begin{figure}[!htb]
\centering
\includegraphics[scale=0.35]{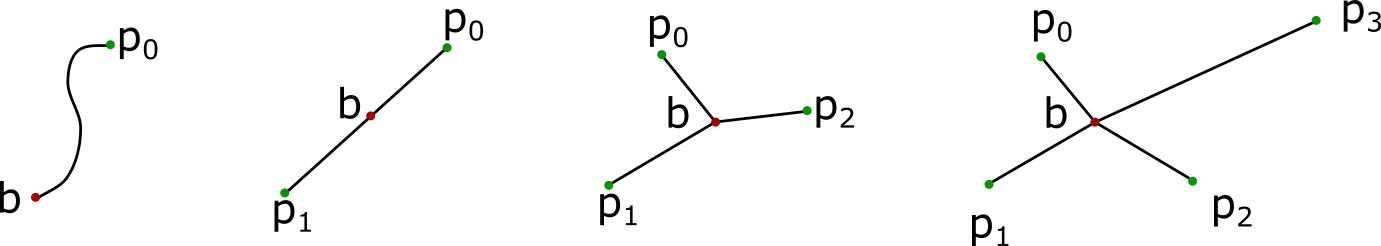}
\caption{Holding bead $b$ with one, two, three or four cables.\label{beadheld}}
\end{figure}

Consider a bead $b$ held on a table top by $N$ cables of length $c_j$ as shown in \Cref{beadheld}.
The cable is anchored to the table at anchor points $p_0 = \left(x_0, y_0 \right)$ up to $p_{N-1} = \left(x_{N-1}, y_{N-1} \right)$. To start with, if the bead was constrained only by a single cable, then it is sure to lie around loose; with two, it can be held tight between $p_0$ and $p_1$ but does not feel very firm in the direction perpendicular to the line $\overline{p_0p_1}$. With three, it can be held ``rather'' securely. However, four gets tricky. In the case of four cables, if the cable from $\overline{bp_0}$ was tightened, then the cable from $\overline{p_0p_3}$ goes slack and vice-versa. However, this is different if it was an elastic band where such a configuration with four cables is possible but can remain unstable.

\begin{figure}[!htb]
\centering
\includegraphics[scale=0.35]{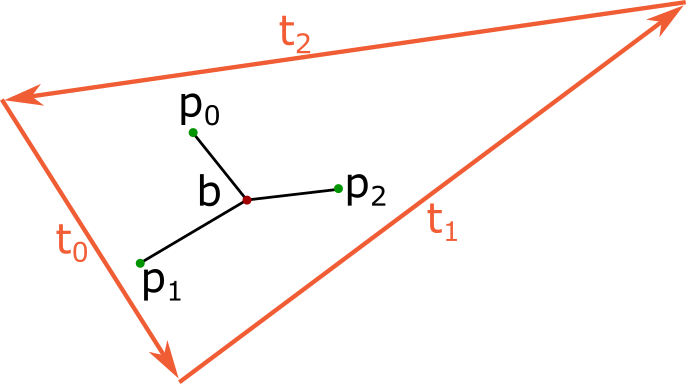}
\caption{The tensions (up to a scale factor) in the three cables balanced at $b$. The vectors are parallel to the cables, but unrelated to the cable lengths: we could move $\mathbf{t}_0$ twice as far from $P$ and yet pull in the same direction with the same force $\mathbf{t}_0$.
\label{ThreeStrings}}
\end{figure}

The case of three strings is a simple case of force balance from engineering mechanics as shown in \Cref{ThreeStrings}.
For this scenario, considering the force vectors $t_1$, $t_2$ and $t_3$ and the position vector of bead $b$, we can say that
\begin{equation}
\left( t_2 - b \right) = a_0 \left( t_0 - b \right) + a_1 \left( t_1 - b \right)
\end{equation}
for unique values of $a_0$ and $a_1$. It is also important to note here that their values are strictly negative since the cables can only consider tension loads. However, the above discussion on uniqueness fails once four cables are considered.

\begin{figure}[!htb]
\centering
\includegraphics[scale=0.5]{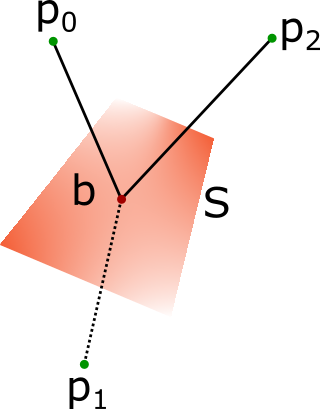}
\caption{With $b$ confined to a plane \textcolor{red}{$S$}, three cables (from above, below or on \textcolor{red}{$S$}) trap it.\label{beadin3D}}
\end{figure}

Now considering the same in 3D as shown in \Cref{beadin3D}. Here, four cables in 3D can hold a bead securely in the same analogy to three in a plane.

If $b$ is stuck to a plane $S$ but cables can reach it from points above, below and points on $S$, then three cables are sufficient. Each $\overline{bp_i}$ sets up a constraint restricted to $S$, giving a circle centered at the nearest $T$-point to $p_i$. Here, the four inequality constraints (with $p_0, p_1, p_2, p_3$ not in the same plane) is exchanged for three inequalities and one $S$-defining equality. For example, if $S$ is the plane $z = 0$, and suppose $t_0 = \left(-2,-1,1\right)$, $t_1 = \left(1,0,2\right)$ and $t_0 = \left(-1,1,-1\right)$, then the constraints on $P$ are
\begin{equation}
\begin{split}
Z &= 0 \\
\left(X+2 \right)^2 + \left(Y+1 \right)^2 + \left(Z-1 \right)^2 &\le 6\\
\left(X-1 \right)^2 + Y^2 + \left(Z-2 \right)^2 &\le 5\\
\left(X+1 \right)^2 + \left(Y-1 \right)^2 + \left(Z+1 \right)^2 & \le 3
\end{split}
\end{equation}
and solving first for $Z$ leads to a two-variable constraint.

Similarly, if $b$ was on the end of a free-turning stick with its other end attached at $Q$, so an equality $\|b-Q\| = radius$ holds it to a spherical surface. Thus, three cables can fix its position, holding it to the tangent plane near $b$ approximates the spherical surface.

\subsection{Rubber vs. cables}
There is an easy way to make the situation with four cables right, and that is to use elastic/rubbery cables, such as elastic bands.

It is fairly easy to show that joining $b$ to any number of fixed points $p_i$ gives an elastic energy function $E(P)$ of the position which is ``strictly convex'' as long as all the cables stay stretched. The configuration also has a unique equilibrium $P_0$ at the minimum-$E$ point, which is stable in two senses: move the bead off $P_0$, and it will move back (elastically stable); change any of the $p_i$ or the elastic constants a little, and $P_0$ will move just a little (structurally stable). Thus, one can easily get close to the designed system. Hence, it is no coincidence that most of the tutorials available on building a tensegrity suggest the usage of rubber bands. 

This also remains one of the origins for the long-standing notion that tensegrities are necessarily pre-stressed structures. As iterated in this discussion and mathematically illustrated in the coming sections, such notion has arisen from the frequent usage of elastic bands for the design of such tensegrity structures. These elastic bands need to be stretched, atleast slightly, if they are to remain in tension which is a pre-requisite for the balance of forces and static equilibrium. However, as will be demonstrated in the later sections, the proposed three-rod and ten-string structure does not necessarily need the strings to be in tension.
% Section 4 (Stability in tensegrities)
\section{Mechanics of stable tensegrities}
As defined earlier, an $n$-rod tensegrity is said to comprise of a set of points $p_0 = (x_0,y_0,z_0), \cdots, p_{2n-1} = (x_{2n-1},y_{2n-1},z_{2n-1})$ with a rod between each pair $\{p_0,p_1\}, \{p_2, p_3\},$ $\cdots, \{p_{2r}, p_{2r+1}\},$ $\cdots, \{p_{2n-2}, p_{2n-1}\}$. The strings in the tensegrity join end pairs $p_{i}$ and $p_{j}$ on different rods, for index pairs $\{i,j\} \in \mathcal{S}$. For each point $p_{i}$ has a set $\mathcal{S}_i$ of indices $j$ for which $p_{i}$ and $p_{j}$ are joined by a string, of which there are $\sigma$ altogether. At present, the relation between $n$, $r$ and $\sigma$ are unknown and the concepts related to these are introduced here.

Irrespective of whether the pair $\{i,j\}$ belongs to rod or string, it is possible to define a vector such as
\begin{equation}
\mathbf{v}_{i,j} = p_{j} - p_{i} \text{ and } \mathbf{v}_{j,i} = p_{i} - p_{j}
\end{equation}
and its modulus as $m_{i,j} = m_{j,i} = \| \mathbf{v}_{i,j} \|$ and the corresponding unit vectors as
\begin{equation}
\hat{\mathbf{v}}_{i,j} = \frac{\mathbf{v}_{i,j}}{m_{i,j}} \text{ and } \hat{\mathbf{v}}_{j,i} = -\hat{\mathbf{v}}_{i,j}.
\end{equation}

Assuming no inelastic effects like buckling, each rod has a compressive force $c_r$ and each string has a tensile force $t_{ij}$. The compressive force acts on the points $p_{2r}$ and $p_{2r+1}$ along the rod and given to be
\begin{equation}
\begin{split}
\mathbf{r}_{2r} &= c_{r} \hat{\mathbf{v}}_{2r,2r+1} = \frac{c_r}{\| \mathbf{v}_{2r,2r+1} \|}  \mathbf{v}_{2r,2r+1} = \Tilde{c}_r  \mathbf{v}_{2r,2r+1} \\
\mathbf{r}_{2r+1} &= c_{r} \hat{\mathbf{v}}_{2r+1,2r} = \frac{c_r}{\| \mathbf{v}_{2r+1,2r} \|}  \mathbf{v}_{2r+1,2r} = \Tilde{c}_r  \mathbf{v}_{2r+1,2r}
\end{split}
\end{equation}
Similarly, the tensile forces can be given to be
\begin{equation}
\mathbf{t}_{ij} = t_{ij} \hat{\mathbf{v}}_{i,j} = \frac{t_{ij}}{\| \mathbf{v}_{i,j} \|} \mathbf{v}_{i,j} = \Tilde{t}_{ij} \mathbf{v_{i,j}}
\end{equation}
It is also evident that $\mathbf{r}_{2r+1} = -\mathbf{r}_{2r}$ and $\mathbf{t}_{ij} = -\mathbf{t}_{ji}$ and the equilibrium at any point $\mathbf{p}_{i}$ requires
\begin{equation}
\mathbf{r}_{i} = \sum_{j \in \mathcal{S}_i} {\mathbf{t}_{ij}}
\end{equation}

\subsection{Constraint matrices}
Considering the vectors in terms of constraints, the configuration change vector can be given as
\begin{equation}
\delta \mathbf{p} = \left( \delta \mathbf{p}_1, \cdots, \delta \mathbf{p}_{2n} \right) = \left(\delta x_1, \delta y_1, \delta z_1, \cdots, \delta x_{2n}, \delta y_{2n}, \delta z_{2n} \right).
\end{equation}
The change in distance between the points $\mathbf{p}_{i}$ and $\mathbf{p}_{j}$, to the first order, can be given to be
\begin{equation}
\begin{split}
\hat{\mathbf{v}}_{i,j} \cdot \left(\delta \mathbf{p}_j - \delta \mathbf{p}_i \right) &= \hat{\mathbf{v}}_{i,j} \cdot \delta \mathbf{p}_j + \hat{\mathbf{v}}_{j,i} \cdot \delta \mathbf{p}_i \\
&= \left[ 0 \ 0 \ 0 \ \cdots \ 0 \ \hat{\mathbf{v}}_{i,j}^{x} \ \hat{\mathbf{v}}_{i,j}^{y} \ \hat{\mathbf{v}}_{i,j}^{z} \ 0 \ \cdots \ 0 \ \hat{\mathbf{v}}_{j,i}^{x} \ \hat{\mathbf{v}}_{j,i}^{y} \ \hat{\mathbf{v}}_{j,i}^{z} \ 0 \ \cdots \ 0 \ 0 \ 0 \right] \begin{bmatrix}
\delta x_0 \\
\delta y_0 \\
\delta z_0 \\
\vdots \\
\delta x_{2n-1} \\
\delta y_{2n-1} \\
\delta z_{2n-1} \\
\end{bmatrix} \\
&= \left[ \hat{C}_{i,j} \right] \left[ \delta \mathbf{p} \right] 
\end{split}
\end{equation}
Here, the non-zero entries are at the locations $3i,3i+1,3i+2$ and their negatives at $3j,3j+1,3j+2$. This can also be re-written as 
\begin{equation}
{\mathbf{v}}_{i,j} \cdot \left(\delta \mathbf{p}_j - \delta \mathbf{p}_i \right) = {\mathbf{v}}_{i,j} \cdot \delta \mathbf{p}_j + {\mathbf{v}}_{j,i} \cdot \delta \mathbf{p}_i = \left[ {C}_{i,j} \right] \left[ \delta \mathbf{p} \right] 
\end{equation}
This leads to the equality constraint $\left[ {C}_{2r,2r+1} \right]  \delta \mathbf{p} = 0$ for each rod $r$ and inequality constraint $\left[ {C}_{i,j} \right]  \delta \mathbf{p} \leq 0$ for each string $\{i,j\}$. Additionally the rigid body translation vectors, given earlier, considers the net axis translations $\left[X\right] \delta \mathbf{p}$, $\left[Y\right] \delta \mathbf{p}$ and $\left[Z\right] \delta \mathbf{p}$ of the centroid $\overline{p}$ of the point set. Similarly, the rotation vectors, given earlier, consider (to linear order) the net turns $\left[\Tilde{X}\right] \delta \mathbf{p}$, $\left[\Tilde{Y}\right] \delta \mathbf{p}$ and $\left[\Tilde{Z}\right] \delta \mathbf{p}$ around the axes.

\subsection{Linearization stability criterion}
It is hypothesized here that a stable tensegrity needs to satisfy the linearization stability criterion stated as
\begin{enumerate}[label=(\alph*)]
\item The set of all rod or string constraint covectors $[C_{i,j}]$ along with the translation and rotation vectors span the $6n$-dimensional space $\mathcal{L}_{n}$ of all possible linear constraints on the configurations $\mathbf{p}$. Thus, for any $\delta \mathbf{p} \neq 0$, there is atleast one non-zero value $\Delta_{i,j} = [C_{i,j}] \ \delta \mathbf{p}$. If $\delta \mathbf{p}$ is not the translation or rotation, then this $\Delta_{i,j}$ must be for a string $\{i,j\}$ and in that case, $\Delta_{i,j} < 0$. In other words, if any $\delta \mathbf{p} \neq 0$ leads to $\Delta_{i,j} = [C_{i,j}] \ \delta \mathbf{p} = 0$, then this $\delta \mathbf{p}$ is a null eigenvector. Since null eigenvectors represent the soft modes, these cannot be permitted.
\item The row vector $\mathbf{0}$ can be expressed as a linear combination with all $a_{i,j} > 0$
\begin{equation}
\begin{split}
\mathbf{0} &= \rho_0 [C_{0,1}] + \cdots + \rho_r [C_{2r,2r+1}] + \cdots + \rho_n [C_{2n-2,2n-1}] + \sum_{\{i,j\} \in \mathcal{S}, i<j} {a_{i,j}} [C_{i,j}] 
\end{split}
\label{eq:eq16}
\end{equation}
However, from the equality constraints, it is known that $\left[ {C}_{2r,2r+1} \right]  \delta \mathbf{p} = 0$ and thus
\begin{equation}
\sum_{\{i,j\} \in \mathcal{S}, i<j} {a_{i,j}} [C_{i,j}] \ \delta \mathbf{p} = \mathbf{0}
\end{equation}
However, from the inequality constraint $\Delta_{i,j} = \left[ {C}_{i,j} \right]  \delta \mathbf{p} \leq 0$. For all $a_{i,j} > 0$, the $\Delta_{i,j} \neq 0$ cannot all be negative. Thus, only possible solution would be if $\delta \mathbf{p} = 0$.
\end{enumerate}

In the $6n$-dimensional space, there are six co-vectors from the rigid body translations and rotations, $n$ from the rods, and $\sigma$ from the strings. To satisfy the condition (a), it is necessary to satisfy
\begin{equation}
6 + n + \sigma \geq 6n \ \Rightarrow \ \sigma \geq 5n-6.
\end{equation}
However, the exact equality constraint is satisfied only if the rod and string co-vectors, together with those from the rigid body modes for a basis in $\mathcal{L}_n$. In such a case, \Cref{eq:eq16} implies that $\rho_r$ and $a_{i,j}$ are all zeros. This contradicts the second criterion (b). Thus, if the constraints form a basis, there necessarily exists a $\delta \mathbf{p}$ which satisfies all the equalities and inequalities and thus leading to an unstable structure.

However, alternatively, if
\begin{equation}
6+n+\sigma = 6n+1 \ \Rightarrow \ \sigma = 5(n-1)
\end{equation}
then the stability condition  (a) implies that there is only one row vector that is a linear combination of the other and this is unique, with exactly one set of co-efficients. Thus, for some $\{I,J\}$ between any two points $\mathbf{p}_i$ and $\mathbf{p}_j$, there exists a non-zero value
\begin{equation}
[C_{I,J}] = \sum_{\text{rod or string} \ \{i,j\} \neq \{I,J\}} {b_{i,j} [C_{i,j}] \ \delta \mathbf{p}}
\end{equation}

Such a set of coefficients $b_{i,j}$ can be easily found using Gauss elimination. Thus, if $5(n-1)$ strings are considered, the test simplifies to the fact that ``all $b_{i,j}$ for string $\{i,j\}$ are negative''. For such a structure, the co-efficients $\rho_r$ and $a_{i,j}$ are also directly related to the compressive and tensile forces in the member and given to be
\begin{equation}
\rho_r = \frac{c_r}{\| \mathbf{v}_{2r,2r+1} \|} \ \text{ and } \ a_{i,j} = \frac{t_{ij}}{\| \mathbf{v}_{i,j} \|}
\end{equation}

Thus, since one of the rows is a linear combination of the others, it is possible to obtain a common scalar multiplier, say $\lambda$, for all the co-efficients. If this holds and all the constraints are in the set, whether or not the rank is maximum, then by uniqueness these co-efficients must correspond to the tensions and compressions, up to some $\gamma$. Without loss of generality, this scalar multiplier can be considered to be
\begin{equation}
\gamma = \frac{1}{\sqrt{\rho_0^2 + \rho_1^2 + \cdots + \rho_{n-1}^2 + \sum_{\{i,j\} \in \mathcal{S}, i<j} {\left( a_{ij} \right)^2 }}}
\end{equation}
then the compressions and tensions can be given to be
\begin{equation}
c_r = \gamma \rho_r \| \mathbf{v}_{2r,2r+1} \| \ \text{ and } \ t_{ij} = \gamma a_{i,j} \| \mathbf{v}_{i,j} \|
\end{equation}

If more than $5(n-1)$ strings are used, then setting the tensions are more complicated. Thus, adjusting the stiffness for a over-determined tensegrity can require a lot of careful tension adjustment and maintenance. Thus, one can conclude that for stability of a $n$-rod tensegrity, atleast $5(n-1)$ strings are required. However, as seen in the earlier sections, most of the designed tensegrities satisfy the condition of $6+n+\sigma \leq 6n$ and thus leading to unstable configurations.
% Section 5 (10-segrity)
\section{Form-finding strategy}
\begin{figure}[!htb]
\centering
\includegraphics[scale=0.15]{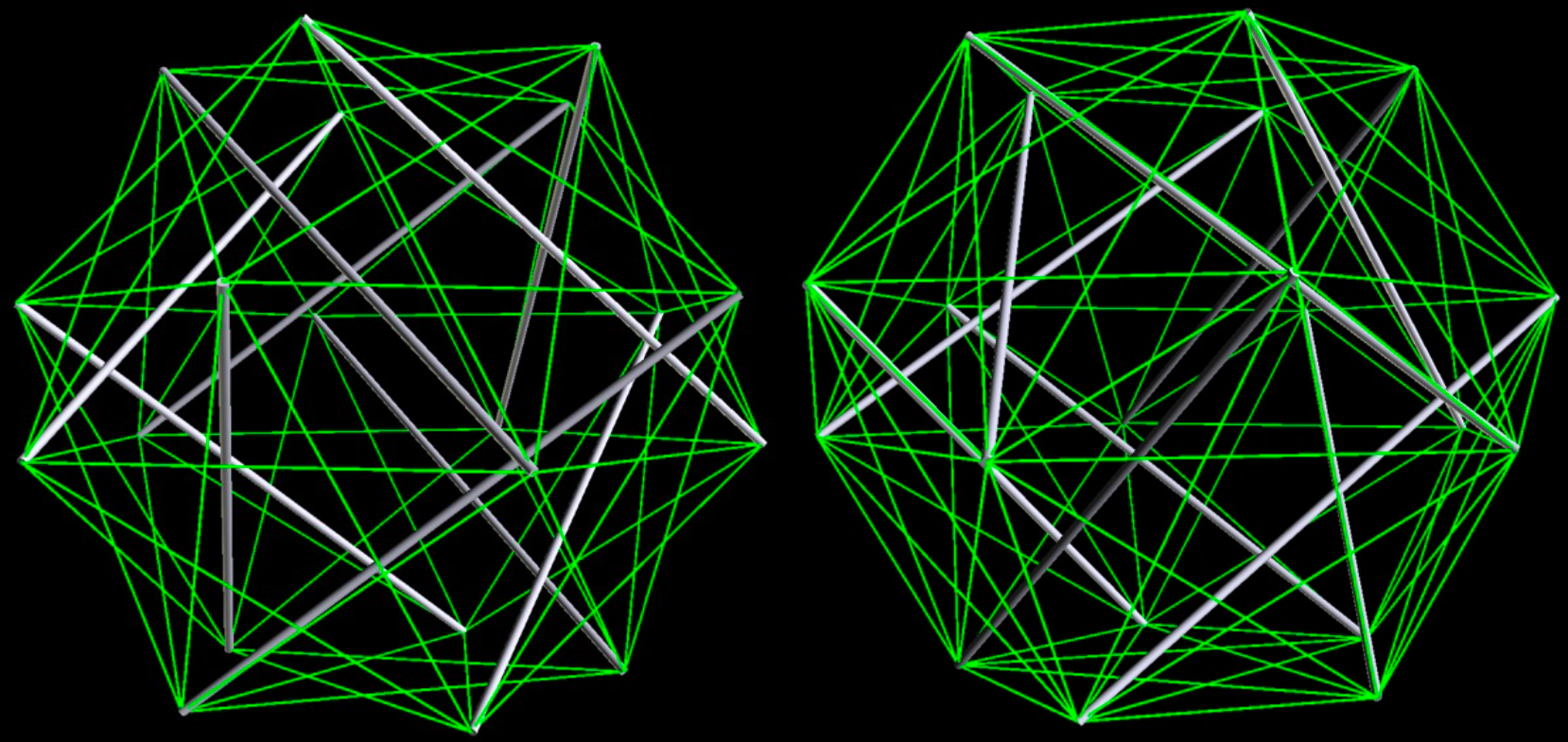}
\caption{A regular dodecahedron tensegrity with 10 rods.\label{dodecahedron}}
\end{figure}

Consider the regular dodecahedron with 10 rods shown in Fig.~\ref{dodecahedron}.
The shape has several options for strings that include 30 face edges, 60 face diagonals, and another 90 connections (not shown in the figure) that cross via the inside of the shape. As discussed in the earlier sections, $5\left(n-1\right) = 45$ strings are required to obtain linearization stability in the model. If one were to allow only the face edges for strings, this is already 53 trillion ways. If face edges and diagonals were allowed, there are $10^{26}$ possible options, and allowing all the inner strings means $6\times 10^{42}$ options. This combinatorial nightmare means that a simple one-by-one test is not feasible, but a better form-finding strategy is necessary.

%For $n$ rods in 3D, the shape and position is a $6n$-dimensional space with $6n$-degrees of freedom/unknowns. There are $n$ rod equalities and six constraints from rigid body motions. Considering this, there are still $5n-4$ unknowns.

If the two ends of each rod $r$ are indexed as $2r$ and $2r+1$, then the unordered list of $2n$ end locations $p_i$ are given to be $p_{2r}=\left(x_{2r},y_{2r},z_{2r}\right)$ and $p_{2r+1}=\left(x_{2r+1},y_{2r+1},z_{2r+1}\right) \ \forall r = 0, 1, \cdots, n-1$. As earlier, the configuration of the tensegrity can be given to be $P = \left(p_0,p_1,\cdots,p_{2n-1}\right)$. Additionally, the tangent vectors $\delta \mathbf{p} = \left( \delta x_{0}, \delta y_{0}, \delta z_{0}, \cdots, \delta x_{2n-1}, \delta y_{2n-1}, \delta z_{2n-1} \right)$ representing the motion at $P$ are all defined in the $6n$-dimensional space $\mathcal{V}$.

Linearizing, each rod $r=0,1,\cdots,n-1$ gives the co-vector
\begin{equation}
\mathbf{e}_{i} = \left[ 0 \ \cdots \ 0 \ x_{2r}-x_{2r+1} \ y_{2r}-y_{2r+1} \ z_{2r}-z_{2r+1} \ x_{2r+1}-x_{2r} \ y_{2r+1}-y_{2r} \ z_{2r+1}-z_{2r} \ 0 \ \cdots \ 0 \right]
\label{eq:eqrowvec}
\end{equation}
to be used in the equality constraint $[C_{2r,2r+1}] \ \delta \mathbf{p} =0$. These vectors $\{\mathbf{e}_i\}$ from \Cref{eq:eqrowvec} are automatically independent and mutually orthogonal. Additionally, they are also orthogonal and independent to the rigid body motion co-vectors unless all the rod ends are co-linear.
%To recall, these are denoted as
%\begin{equation}
%\begin{split}
%\mathbf{e}_{n} &= \left[1 \ 0 \ 0 \ 1 \ 0 \ 0 \ \cdots \ 1 \ 0 \ 0 \right] \\
%\mathbf{e}_{n+1} &= \left[0 \ 1 \ 0 \ 0 \ 1 \ 0 \ \cdots \ 0 \ 1 \ 0 \right] \\
%\mathbf{e}_{n+2} &= \left[0 \ 0 \ 1 \ 0 \ 0 \ 1 \ \cdots \ 0 \ 0 \ 1 \right] \\
%\mathbf{e}_{n+3} &= \left[0 \ z_1 \ -y_1 \ 0 \ z_2 \ -y_2 \ \cdots \ 0 \ z_{2n-1} \ -y_{2n-1} \right] \\
%\mathbf{e}_{n+4} &= \left[-z_1 \ 0 \ x_1 \ -z_2 \ 0 \ x_2 \ \cdots \ -z_{2n-1} \ 0 \ x_{2n-1} \right] \\
%\mathbf{e}_{n+5} &= \left[y_1 \ -x_1 \ 0 \ y_2 \ -x_2 \ 0 \ \cdots \ y_{2n-1} \ -x_{2n-1} \ 0 \right]
%\end{split}
%\end{equation}
The complete set comprising of $n+6$ co-vectors from rods and rigid body motions can be denoted by the set $\mathcal{E}$. 

There are a total of $2n\left(n-1\right)$ unordered pairs possibilities
\begin{equation}
\left\{i,j\right\} \neq \left\{ 2r,2r+1 \right\} \ \text{, any } r
\end{equation}
for string connections. For the above dodecahedron in \Cref{dodecahedron}, this implies 180 possibilities, as discussed earlier. It is possible that one might aim to select string locations from all these pairs or from a subset $\mathcal{T}$ chosen by some \textit{adhoc} criterion that limits the length or the position. Without loss of generality, the size of this subset can be stated to be $L+1$. If each of the pair is listed as $\{ i_{\ell}, j_{\ell} \} \ \forall \ \ell = 0, 1, \cdots, L$, then there exists a set $\left( \ell_0, \ell_1, \cdots, \ell_L \right)$ associated with the vectors $\mathbf{v}_{\ell}$ and given as $\ell_0 \mathbf{v}_{0} + \cdots + \ell_{L} \mathbf{v}_{L}$.

\Cref{eq:eq16} states that
\begin{equation}
\begin{split}
\mathbf{0} &= \rho_0 [C_{0,1}] + \cdots + \rho_r [C_{2r,2r+1}] + \cdots + \rho_n [C_{2n-2,2n-1}] + \sum_{\{i,j\} \in \mathcal{S}, i<j} {a_{i,j}} [C_{i,j}] 
\end{split}
\end{equation}
and since it known that ${a_{i,j}} > 0$, it can be re-written as
\begin{equation}
\begin{split}
&\mathbf{0} = \rho_0 [C_{0,1}] + \cdots + \rho_r [C_{2r,2r+1}] + \cdots + \rho_n [C_{2n-2,2n-1}] + \sum_{k = 0}^{L} {\lambda_k} [C_{k}] \\
\Rightarrow & \mathbf{X} = -\rho_0 [C_{0,1}] - \cdots + \rho_r [C_{2r,2r+1}] - \cdots - \rho_n [C_{2n-2,2n-1}] = \sum_{k = 0}^{L} {\lambda_k} [C_{k}]
\end{split}
\end{equation}
where all $\lambda_k \geq 0$ and $\lambda_0 + \lambda_1 + \cdots + \lambda_L = 1$ without loss of generality. Here, a mapping is used to relate the unordered pairs with the co-efficients of the string co-vectors.

The co-efficients $\left( \lambda_0, \lambda_1, \cdots, \lambda_L\right)$ form a  $(L+1)$-dimensional subspace known as $L$-\textbf{simplex}. Since the interest is to determine $5(n-1)$ non-zero and positive values, this implies the intersection of the above defined simplex $\mathcal{S}$ with a space $\mathcal{E}$ of dimension $5n-4$. The intersection of this space $\mathcal{E}$ and the simplex $\mathcal{S}$ at a face leads to some $\lambda_{\ell}$ to be zero and one can leave out the corresponding $\mathbf{v}_{\ell}$, i.e. the string inequality constraint. 

%Upon elimination, and determination of all non-zero $\lambda$'s, the vector $\mathbf{X}$ can be determined to check the existence of $\rho_i$'s such that

%where $[C_{i,j}]$ are co-vectors from the equality (rod) and rigid-body modes.

\begin{figure}[!htb]
\centering
\includegraphics[scale=0.25]{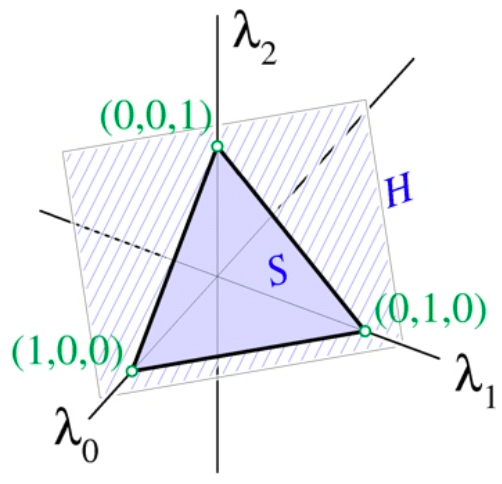}
\caption{The affine plane $H$ contains all the $\left(\lambda_0, \lambda_1, \lambda_2 \right)$ with $\sum_{K} {\lambda_K} = 1$. The triangle $S$ shows where they are all positive. If $L=3$, the set S becomes a tetrahedron in the 3D hyperplane $H$ of $\left(\lambda_0, \lambda_1, \lambda_2, \lambda_3 \right)$-space and so on\label{LSimplex}}
\end{figure}

\begin{figure}[!htb]
\centering
\includegraphics[scale=0.25]{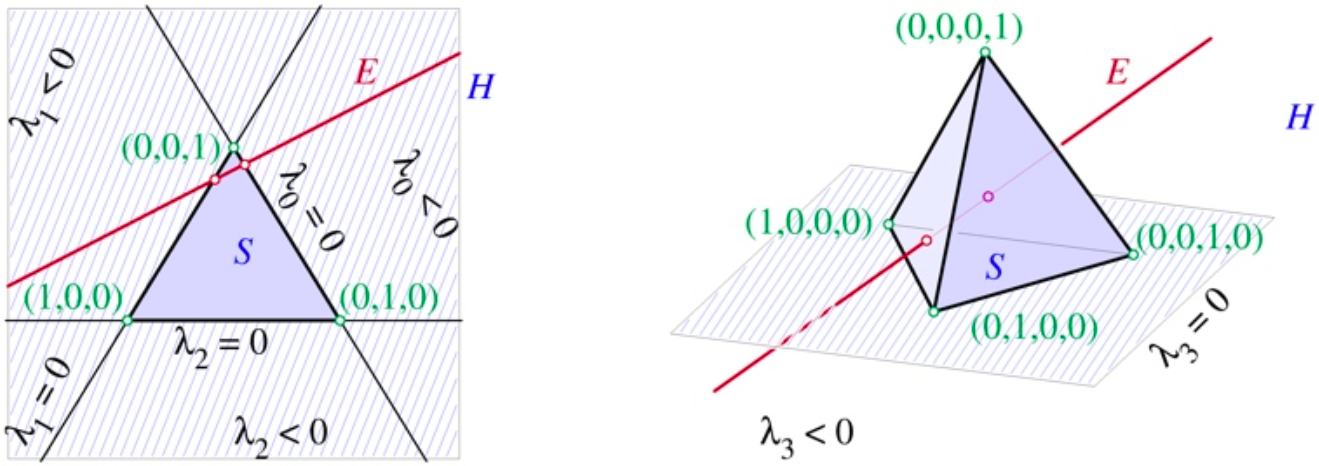}
\caption{Intersection of line with a 2- (left) and 3- (right) simplex.\label{SimplexWithLine}}
\end{figure}

For example, for the dodecahedron from \Cref{dodecahedron} and if all ordered pairs are considered, then $L=180$. It is evident that visualizations of such larger dimensional spaces can be harder and thus it is prudent to develop the concept using a lower dimensional simplex. \Cref{LSimplex} shows the simplex for $L=2$.
The \Cref{SimplexWithLine} shows the intersection of a space with the simplex. Only here, the space $\mathcal{E}$ is a line $E$ intersecting with a simplex $\mathcal{S}$ in 2D and 3D and ensure one or more $\lambda$'s to be zero.

\begin{figure}[!htb]
\centering
\includegraphics[scale=0.25]{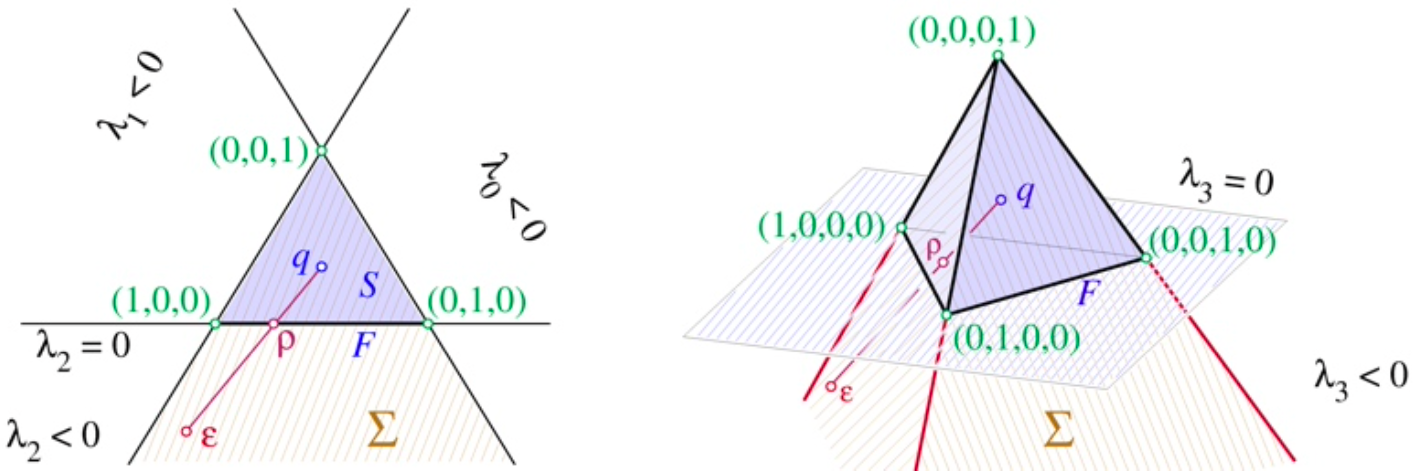}
\caption{Determination of the point $\rho$ which is on the face of the Simplex $\mathcal{S}$ and also on the line segment $\overline{q\varepsilon}$. Here, the point $\varepsilon$ is the nearest point to the centroid of the simplex.\label{search}}
\end{figure}

The search is posed as a nonlinear iterative procedure that is pictorially depicted in the \Cref{search} 
and described below
\begin{itemize}
\item Define the centroid of the simplex. The centroid $q \in \mathcal{S}$ can be defined as 
\begin{equation}
q = \left(m,\cdots,m\right), \text{ where } m = \frac{1}{L+1}
\end{equation}
\item Determine a point, say $\varepsilon = \left(\varepsilon_0, \cdots, \varepsilon_{\ell}, \cdots, \varepsilon_{L+1} \right)$ such that there exists $\rho_r$ that satisfies
\begin{equation}
\mathbf{X} = -\rho_0 [C_{0,1}] - \cdots + \rho_r [C_{2r,2r+1}] - \cdots - \rho_n [C_{2n-2,2n-1}] = \sum_{k = 0}^{L} {\varepsilon_k} [C_{k}]
\end{equation}
and $\varepsilon \in \mathcal{H}$ is the nearest point to $q \in \mathcal{S}$.
\item If the chosen point $\varepsilon$ is a point of the simplex $\mathcal{S}$, then the search is complete
\item If the chosen point $\varepsilon$ is not a point of the simplex $\mathcal{S}$, then it is necessary to determine the face $F$ through which the line $Q$ from $q$ to $\varepsilon$ leaves $\mathcal{S}$. If $\mathbf{r}$ is the vector from $q$ to $\varepsilon$, then $Q = \left\{ q + t\mathbf{r} \ | \ t \in \mathbb{R} \right\}$. Here, each hyperplane $\lambda_{\ell}$ meets $Q$ at
\begin{equation}
\rho = q + q + t_{\ell} \mathbf{r} = \left( \frac{m}{m-\varepsilon_{\ell}} \right) \mathbf{r}
\end{equation}
and the segment $\overline{q\varepsilon}$ leaves $\mathcal{S}$ at the smallest positive $t_{\ell}$ and this corresponds to the most negative $\varepsilon_{\ell}$. Referring to this index with \underline{$\ell$}, one would consider the space $\Sigma$ as the convex set of all points where $\lambda_{\ell \neq \underline{\ell}} \geq 0$. Here, either $\varepsilon$ is in $\Sigma$ with $\lambda_{\underline{\ell}}$ as the only negative $\lambda_{\ell}$ or not.
\item If $\varepsilon$ is in $\Sigma$, then the string $\mathbf{v}_{\ell}$ can be dropped and the remaining set of strings leave the search as viable as it was.
\item In the case that $\varepsilon$ is not in $\Sigma$, then this implies that $\mathcal{E}$ does not meet $\mathcal{E}$ at a face. Thus, one of the $\varepsilon_{\ell} < 0$ is safe to omit. 
\item The procedure terminates when atleast one branch reaches a solution with $\varepsilon_{\ell} > 0$. If every branch ends with a solution that is $\varepsilon_{\ell} < 0$, then there is no solution to the original problem.
\end{itemize}
The above provides a recursive strategy that facilitates a solution with every step with branching, if required to omit $\ell$ with $\varepsilon < 0$. This procedure is much faster than a brute force approach that tries all combinations of $5(n-1)$ strings among the candidate set $\mathcal{T}$.

% Section 6 (Results)
\section{Results and discussions}
The famous three-rod tensegrity that has been often used in studying the mechanics of tensegrity structures has been considered here to demonstrate the applicability of the developed methods to find stable structures. In order to differentiate between the stable and unstable structure, the below nomenclature is employed for convenience.
\begin{enumerate}
\item Three-rod with nine-strings (which has a swinging mode): 9-segrity
\item Three-rod with ten-strings (fully stable): 10-segrity
\end{enumerate}
The rod length for both the models are considered to be 4 units to provide a one-to-one comparison. The geometry of the 9-segrity considered is as given below. Consider a scale that puts the bottom points on the unit circle at height $z=0$. The bottom three points are
\begin{equation}
\begin{split}
p_0 &= \left(1,0,0\right) \\
p_1 &= \left(\cos{\left(\frac{2\pi}{3}\right)},\sin{\left(\frac{2\pi}{3}\right)},0\right) = \left(-1/2,\sqrt{3}/2,0\right)\\
p_2 &= \left(\cos{\left(\frac{-2\pi}{3}\right)},\sin{\left(-\frac{2\pi}{3}\right)},0\right) = \left(-1/2,-\sqrt{3}/2,0\right)
\end{split}
\end{equation}
The top three points are given to be
\begin{equation}
\begin{split}
q_0 &= \left(\cos{\left(\theta\right)},\sin{\left(\theta\right)},h\right) \\
q_1 &= \left(\cos{\left(\theta + \frac{2\pi}{3}\right)},\sin{\left(\theta + \frac{2\pi}{3}\right)},h\right)\\
q_2 &= \left(\cos{\left(\theta - \frac{2\pi}{3}\right)},\sin{\left(\theta - \frac{2\pi}{3}\right)},h\right)
\end{split}
\end{equation}
Here, two configurations are considered with $\theta = -100^{\circ}$ and $210^{\circ}$. The height of the tensegrity is adjusted to ensure a rod length of 4 units using the relation
\begin{equation}
2\left(1-\cos{\theta} \right) + h^2 = r^2
\end{equation}

\begin{figure}[!htb]
\centering
\includegraphics[scale=0.25]{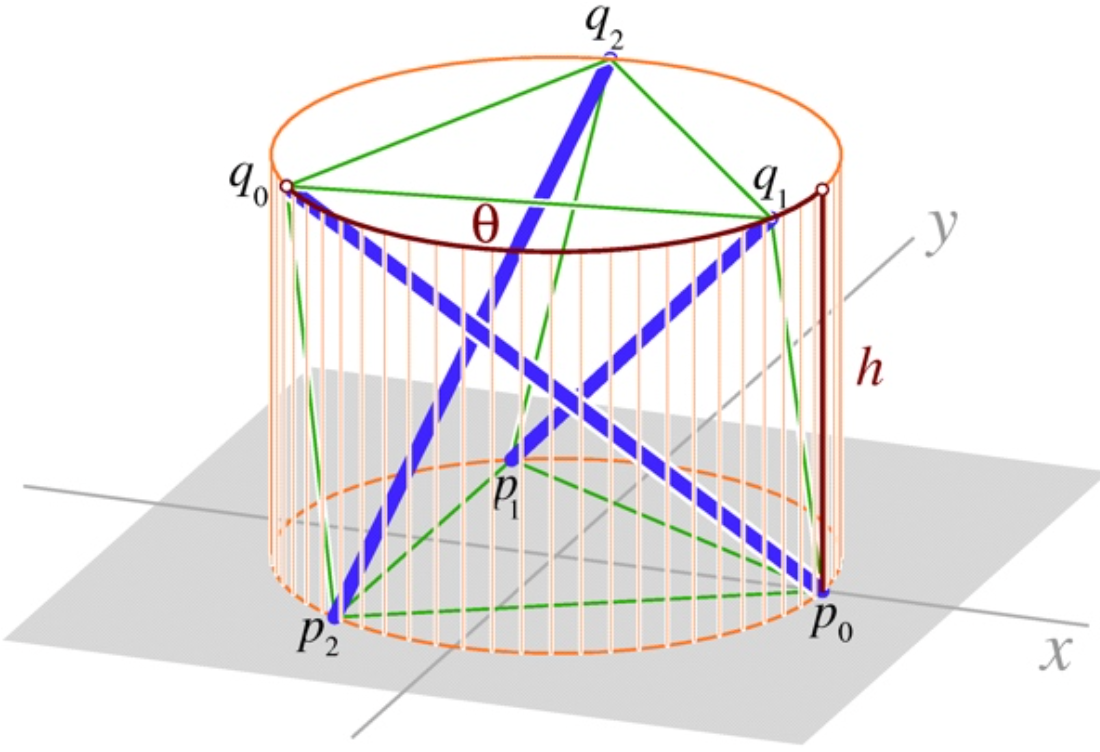}
\caption{Depiction of the 9-segrity. The co-ordinates $h$ and $\theta$ on a subset of the shapes possible for a 9-segrity, with ends on the radius-1 cylinder. The lower ends are all $\sqrt{3}$ apart; so are the top ones, at height $h$. The rod and cable lengths follow from how far $\theta$ at the top twists, relative to the base.\label{3rod9segrity}}
\end{figure}
As shown in the \Cref{3rod9segrity}, the rods are positioned such that the pairs of rod ends are given by $\left(p_0,q_0\right)$, $\left(p_1,q_1\right)$ and $\left(p_2,q_2\right)$
and the nine string pairs, without repetition, are given by $\left(p_0, p_1 \right)$, $\left(p_0, p_2 \right)$, $\left(p_0, q_1\right)$, $\left(p_1, p_2 \right)$, $\left( p_1, q_2\right)$, $\left(p_2, q_0\right)$, $\left( q_0, q_1\right)$, $\left(q_0,q_2 \right)$, $\left(q_1,q_2  \right)$.

\subsection{9-segrity: Theoretical observations}
The 9-segrity has an inherent symmetry. Considering the squared distance of the point $q_0$ w.r.t. $p_0$, considering a fixed rod length of $r$ units for all rods,
\begin{equation}
\left( \cos{\theta} - 1 \right)^2 + \left( \sin{\theta} \right)^2 + h^2 = 2-2\cos{\theta} + h^2 \Rightarrow 2\left(1-\cos{\theta} \right) + h^2 = r^2
\end{equation}
Similarly, the squared distance of point $q_0$ from $p_1$ is 
\begin{equation}
\left( \cos{\theta} + \frac{1}{2} \right)^{2} + \left( \sin{\theta} - \frac{\sqrt{3}}{2} \right)^{2} + h^2 = 2 - 2 \cos{\left(\theta + \frac{2\pi}{3}\right)} + h^2 
\end{equation}
and thus joining them by a cable of length $c$ gives a constraint
\begin{equation}
2\left[1-\cos{\left(\theta+\frac{2\pi}{3}\right)} \right] + h^2 \leq c^2
\end{equation}

\begin{figure}[!htb]
\centering
\includegraphics[scale=0.6]{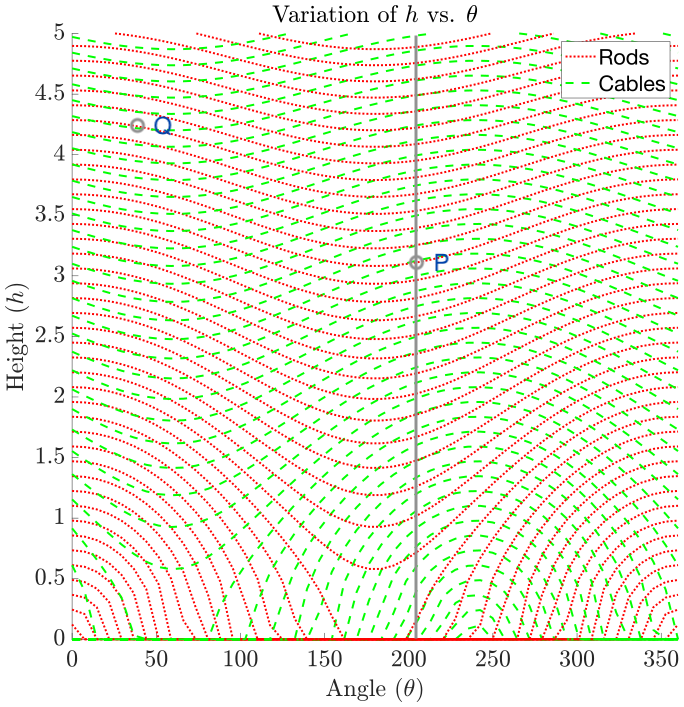}
\caption{The variation of $h$ vs. $\theta$ for various fixed r (red) and c (green). Along each particular red curve, the square of the distance from $p_0$ to $q_0$ is constant. Along each green curve the square of the distance between $p_1$ to $q_0$ is s constant.\label{constraintunroll}}
\end{figure}

It is easier to see how these two constraints interact if one unrolls the cylinder onto a flat diagram as shown in \Cref{constraintunroll}.
For any particular rod length $r$, the rod forces the upper end to be on the particular corresponding red curve, but able to move along it as $\theta$ changes. Similarly, upon choosing a particular cable length $c$, the cable forces the upper end to be on or below a particular corresponding green curve. Wherever a red and green curve cross at an angle, the $\theta$ can be varied to move along the red curve, and go downward from the green one. Thus, the point is not fixed and can move without violating the constraints. At a point like $P$, the two curves are tangent and the shape is free to move along the common tangent line. The linear approximations given by
\begin{equation}
\begin{split}
\left[ \frac{\partial\left(\text{red}\right)}{\partial \theta} \ \frac{\partial\left(\text{red}\right)}{\partial h} \right] &= \left[ 2\sin \theta \ \ 2h \right] \\
\left[ \frac{\partial\left(\text{green}\right)}{\partial \theta} \ \frac{\partial\left(\text{green}\right)}{\partial h} \right] &= \left[ 2\sin \left(\theta + \frac{2\pi}{3} \right) \ \ 2h \right] 
\end{split}
\end{equation}
The above functions are identical at points like $P$ where $\theta = 210^{\circ}$ and the above reduces to
\begin{equation}
\begin{split}
-\delta \theta + 2h \delta h &= 0 \\
-\delta \theta + 2h \delta h &\leq 0 \\
\end{split}
\end{equation}
where the second automatically holds true if the first does. Looking beyond the linear approximation at $P$, the green curve is tangent to the red curve from below. Thus, to second order, therefore, following the red curve strictly increases the cable length, moving to higher green curves. This is forbidden by the constraint, so the shape is fixed.

However, there are also points like $Q$, where a green curve touches the red from above. This gives another situation where a linear analysis, like above, has a constant tangent and thus requires higher order terms. In this case, moving along the red curve quadratically lowers the green value from a maximum at $Q$ and the tensegrity can fall apart. However, if build exactly, while the forces will balance each other but lead to an unstable equilibrium.

\begin{figure}[!htb]
\centering
\includegraphics[scale=0.75]{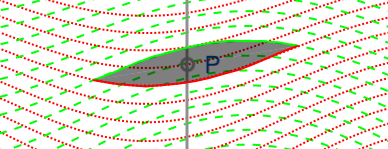}
\caption{Perturbation from the red and green curves tangent at $P$. A similar shape occurs if they both move up, or both down, by different amounts.\label{zoom2p}}
\end{figure}

Further, as iterated throughout literature, firm tensegrities are hard to make, except with elastic bands. If the $r$ and $c$ do not lead to curves that meet exactly and tangentially at $P$ , they will likely either cross each other (twice) or fail to meet at all. If they do not meet, then the measured cable is too short, and the tensegrity cannot be made with it. \Cref{zoom2p} shows the area near the point $P$.
Here, the crossing pair of curves is shown distinctly. A small perturbation from passing through $P$ produces a much larger grey area between the curves. This amplification comes precisely from the ideal tangency. As second order contact is perturbed the widest separation of the two curves, near the contact point, grows like the errors in $r$ and $c$. If the perturbation produces crossings, then the distance between the crossings grows like the square root of this separation. Note here, that if the cross is a small number $\delta$, then its square root is bigger than $\delta$. Thus, the distance grows drastically with fabrication error. All of this allow a large amount of ``swinging'', compared to the fabrication error. Thus, very careful fabrication is necessary to avoid this.

As proposed in several literature, one intends to use such a tensegrity as a deployable structure, particularly of interest to space applications. If the number of length constraints are too small to meet the full-rank convexity condition, their typical situation automatically allows movement away from the target shape. The lengths must be very precisely and non-linearly co-ordinated to keep them in the relationship needed for the degenerate condition of tangency and a swinging-free shape. At every stage of adjusting them, they must not deviate into the typical kind of crossing. Further on, any adjustment mechanism has moving parts, and is thus subject to wear. Preserving tangency through the process of mutually non-linear length adjustments, and maintaining this precision through a lifetime of service cycles, is a narrow and expensive target. Thus, such structures are not optimal nor suitable as deployable structures.

\subsection{9-segrity with an additional cable is still a 9-segrity}
\begin{figure}[!htb]
\centering
\includegraphics[scale=0.25]{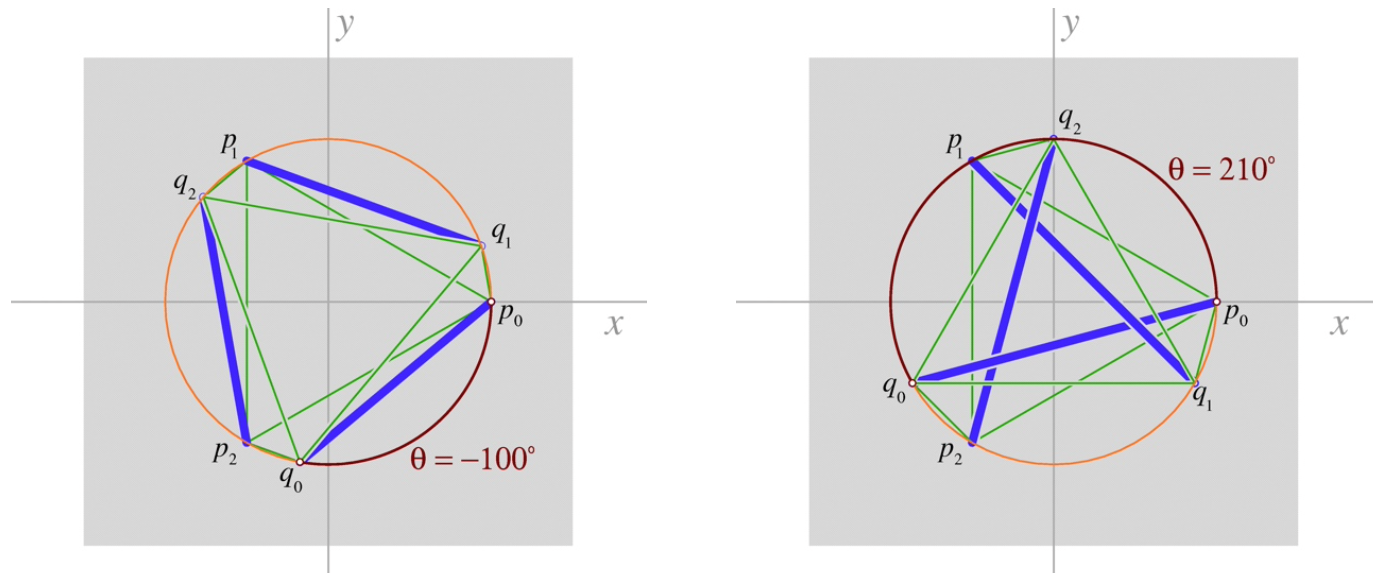}
\caption{A view of the 9-segrity from directly above. There are clearly no cables that can hold rods in these positions against an inward swing. Only equilibrium angle allowed is $\theta = 210^{\circ} (right)$.\label{stopturning}}
\end{figure}

A stable and firm tensegrity is possible with three rods and ten cables. However, as shown later in \Cref{10segrity}, this cannot be achieved by just adding one cable to the linearly degenerate system of the three-rod and nine-string tensegrity. This section outlines the discussion and reasoning behind this conclusion. From a simple observation of the top view of the 9-segrity shown in \Cref{stopturning}, it is evident that there are no cables that can prevent it from swinging inwards.

\begin{figure}[!htb]
\centering
\includegraphics[scale=0.3]{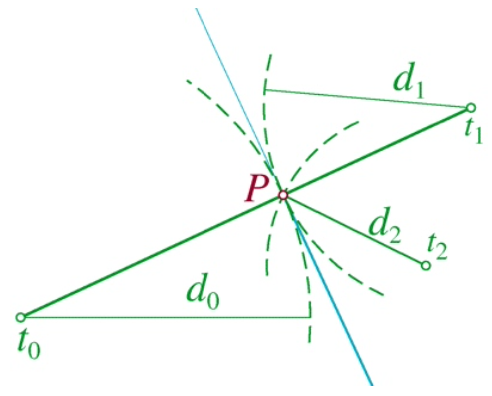}
\caption{Point $P$ upon addition of one additional cable $d_2$ from $t_2$.\label{P3wires}}
\end{figure}

\begin{figure}[!htb]
\centering
\includegraphics[scale=0.15]{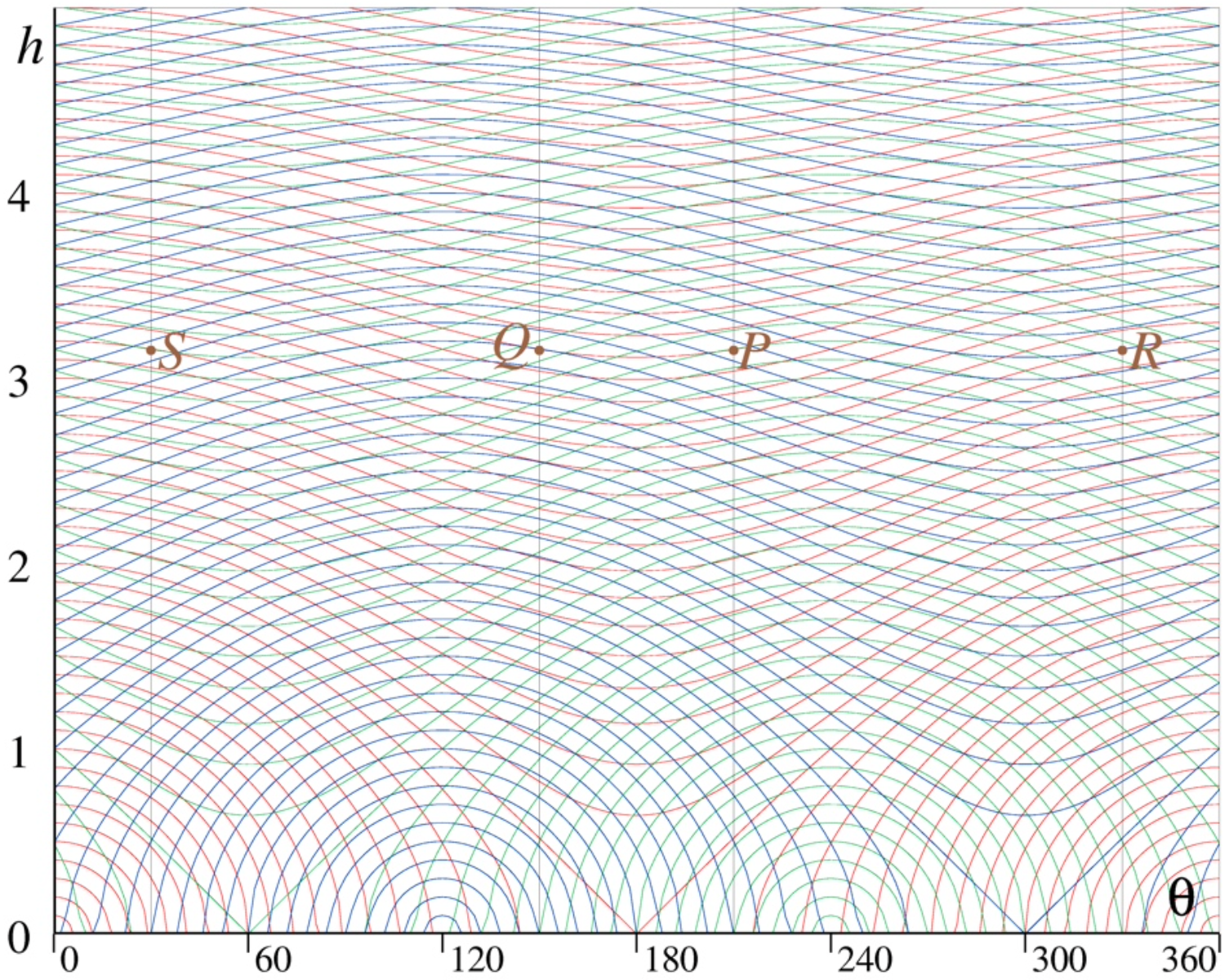}
\caption{Curves of equal distance (through the unit cylinder) from $q_0$ to $p_0$ (red), $p_1$ (green) and $p_2$ (blue). It is easy to check that the green and red curves are tangent at $\theta = 30^{\circ}, 210^{\circ}$, the blue and red curves at $\theta = 150^{\circ}, 330^{\circ}$.\label{constraintunroll2}}
\end{figure}

Hypothetically, if one were to add one additional string to the 9-segrity, the stability at point $P$ is as shown in \Cref{P3wires}.
Here, one additional cable $d_2$ is added from $t_2$. This stops the point $P$ from going up the line of degenerate vector. While, it cannot cross the tangent to the $d_2$ circle centered at $t_2$, it can still move along the grey area outlined earlier. The system fails the full-rank convexity test. A tenth cable, unless can be tightened enough to change the length of the existing nine cables, would have zero tension, and be almost slack. This is further illustrated by rolling out the constraint equations, as earlier, in \Cref{constraintunroll2}.
Using three equal cables to join more ends ($p_2$ as well as $p_1$ to $q_0$ , plus $p_0$ to $q_1$ and $p_1$ to $q_2$) will over-constraint the system but still does not get the structure out of the certainly-soft, hard-to-build zone. From \Cref{constraintunroll2}, if neither the green function nor the blue can increase, and one can move only along a particular red curve, still any motions that these constraints do permit, definitely leads to an instability. If none of these motions are permitted, that does not prove stability either.

In \Cref{constraintunroll2}, one can see that not increasing the blue values also prevents $P$ from moving to the right along its red equality curve, but not to the left! Holding $P$ against this motion depends on higher order terms. Again, anywhere between $P$ and $R$, moving to the left along a red curve actually decreases both the green and blue values and thus these constraints cannot hold a point still either. In the wide band on the left, between $S$ and $Q$, moving right along a red curve decreases both green and blue, and is allowed by the cables. Between the lines through $P$ and $Q$, however, not increasing the green values prevents leftward motion along any red curve, while not increasing blue values blocks it on the right. 

A target range is in the strip $150^{\circ} < \theta < 210^{\circ}$ between $P$ and $Q$. Suppose there is a small error in choosing $(\theta,h)$ and instead one chooses $(\theta',h')$, then this small error is sufficient to cause a soft motion or collapse. IN other words, if $P'$ is obtained instead of $P$ such that $\| P' - P \| < \varepsilon$, then the convexity condition holds true for $P'$ too and thus allows for the swinging between $P$ and $P'$.

\subsection{10-segrity: A stabilized model}
\begin{figure}[!htb]
\centering
\includegraphics[scale=0.15]{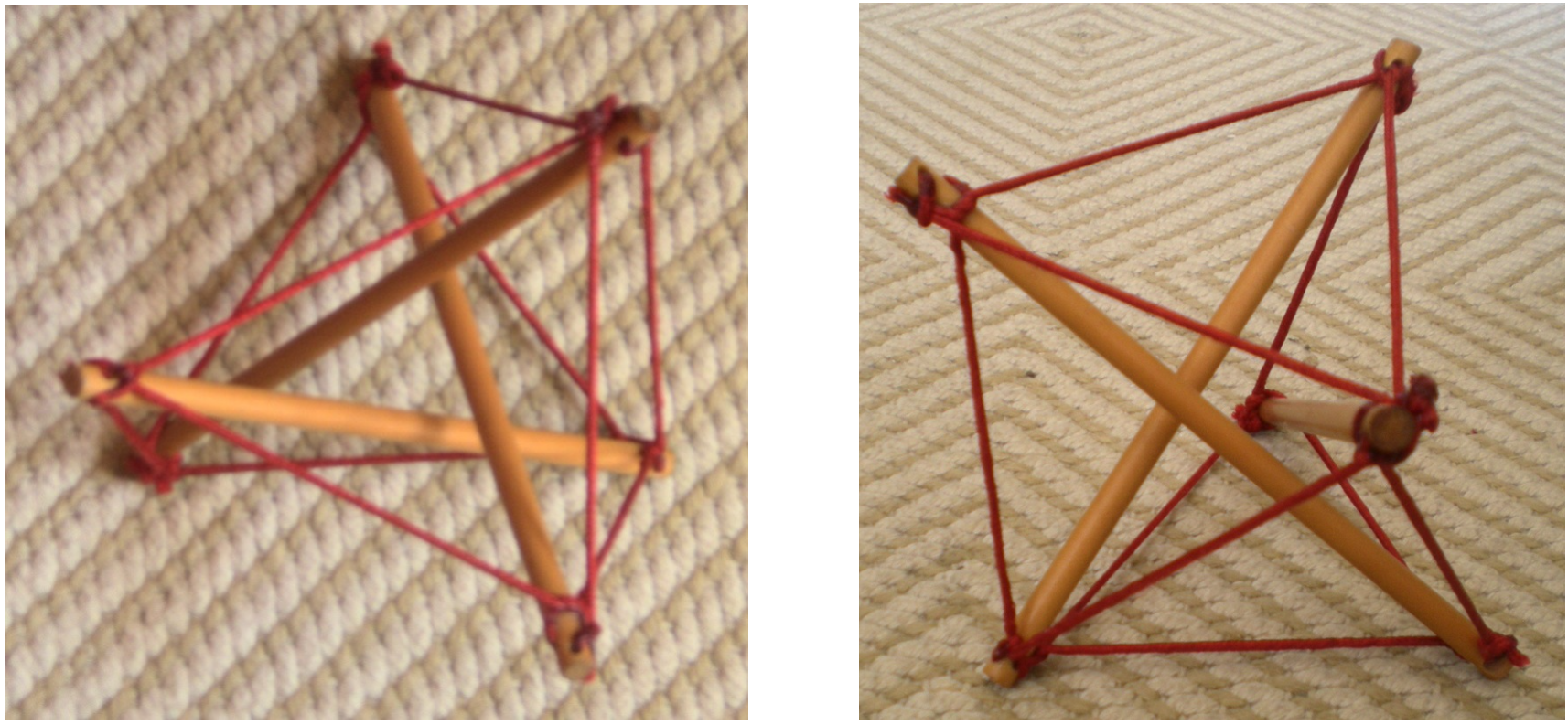}
\caption{Firm three-rod tensegrity with ten cables. Looking down (left) and from front (right).\label{10segrity}}
\end{figure}

Using the proposed form-finding technique, an alternative three-rod model consisting of ten strings is proposed here, as shown in \Cref{10segrity}, that is stable and free of swinging modes. 
The configuration of the 10-segrity is given by the unordered sets for one set of rod ends
\begin{equation}
\begin{split}
p_0 &= (0.5,0,-2) \\ %N1
p_1 &= (0,-2,0.5) \\ %N2
p_2 &= (2,0.5,0) %N5
\end{split}
\end{equation}
and the other set of rod ends given by
\begin{equation}
\begin{split}
q_0 &= (0.5,0,2)\\ %N3
q_1 &= (0,2,0.5) \\ %N4
q_2 &= (-2,0.5.0) %N6
\end{split}
\end{equation}
The rods are position such that the pairs of rod ends are given by $(p_0,q_0)$, $(p_1,q_1)$ and $(p_2,q_2)$ and the set of string pair ends are given by $(p_0,p_1)$, $(p_0,p_2)$, $(p_0,q_1)$, $(p_1,q_0)$, $(p_1,p_2)$, $(p_1,q_2)$, $(q_0,p_2)$, $(q_0,q_2)$, $(q_1,p_2)$ and $(q_1,q_2)$.

In order to demonstrate the stability of the proposed 10-segrity, a modal analysis is considered and discussed here.
\begin{figure}[!htb]
\centering
\includegraphics[scale=0.3]{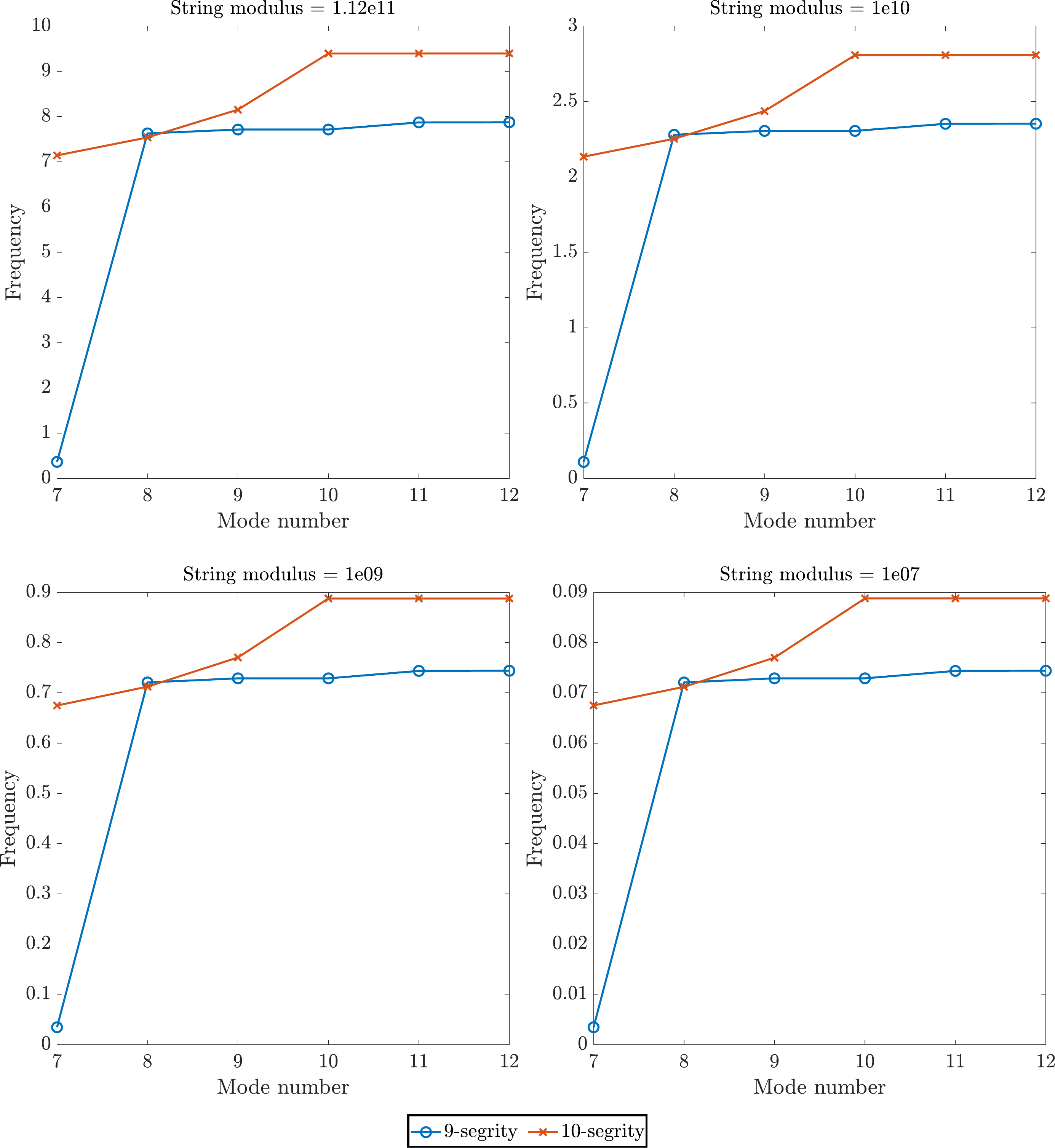}
\caption{Comparison of natural frequencies for the 9- and 10-segrity for various modes (7 - 12). The first six modes represent rigid body motions. For the 9-segrity, here $\theta = 210^{\circ}$ \label{EigenValues}}
\end{figure}

\subsubsection{Modal analysis (No constraints)}
In order to compare the stability of the proposed model with the original 9-segrity with the proposed 10-segrity, a modal analysis is performed to extract the first 12 eigenmodes.

The 9-segrity considered here has a $\theta = 210^{\circ}$ since this is the most stable configuration possible. In order to keep the comparison acceptable, both the models are considered with the same rod lengths, i.e. 4 units. The rods are considered to be made of steel and strings of Kevlar. Both as assumed to have a circular cross-section. The materials are considered to be linearly elastic in nature. The properties of the steel rod is considered to have a Youngs' modulus of $210$ GPa, Poisson ratio of 0.3, density of 8000 kg/m$^3$ and a cross-sectional radius of 0.01m. The Kevlar strings are considered to have a Poisson ratio of 0.36, density of 1440 kg/m$^3$ and cross-sectional radius of 0.001m. Four values of  Youngs' modulus of the Kevlar string are considered to check the influence of material properties: 112 GPa, 10 GPa, 1 GPa and 10 MPa. However, the strings are considered to have no stiffness under compression. The rods and strings are modeled using 3D, second-order Timoshenko beam elements. In this work, we consider \texttt{B32} element from Abaqus standard \cite{AbaqusManual}. The rods are discretized with 20 elements and the strings with 14-15 elements. 

In this first case, there are no boundary conditions applied on both the models. As expected, the natural frequencies for the first six modes are zeros indicating a rigid body mode. Thus, only the modes 7 - 12 are considered for analysis. The variation of the natural frequency of vibration for modes 7 - 12 as a function of the Youngs' modulus of Kevlar are shown in \Cref{EigenValues}. As discussed in earlier literature, it is clearly evident here that the lowest mode (other than the rigid body modes) is a soft / swinging mode of deformation with a natural frequency at 0.3655 rad/s, i.e. nearly zero.

\begin{figure*}[!htb]
\centering
\begin{subfigure}{.12\textwidth}
  \centering
\end{subfigure}
\begin{subfigure}{.24\textwidth}
  \centering
  \includegraphics[width=0.9\textwidth]{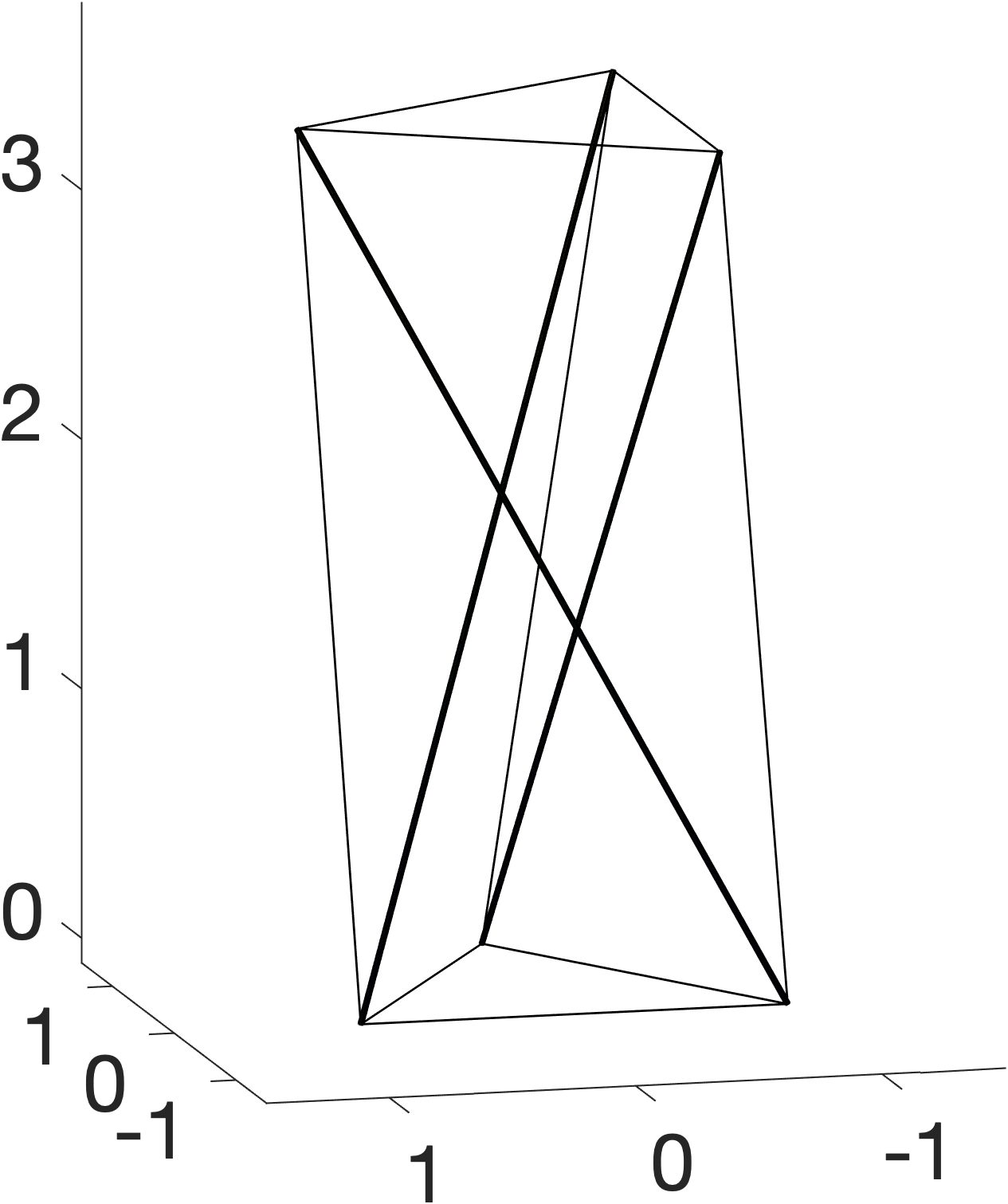}
  \caption{Original shape}
\end{subfigure}
\begin{subfigure}{.24\textwidth}
  \centering
  \includegraphics[width=0.9\textwidth]{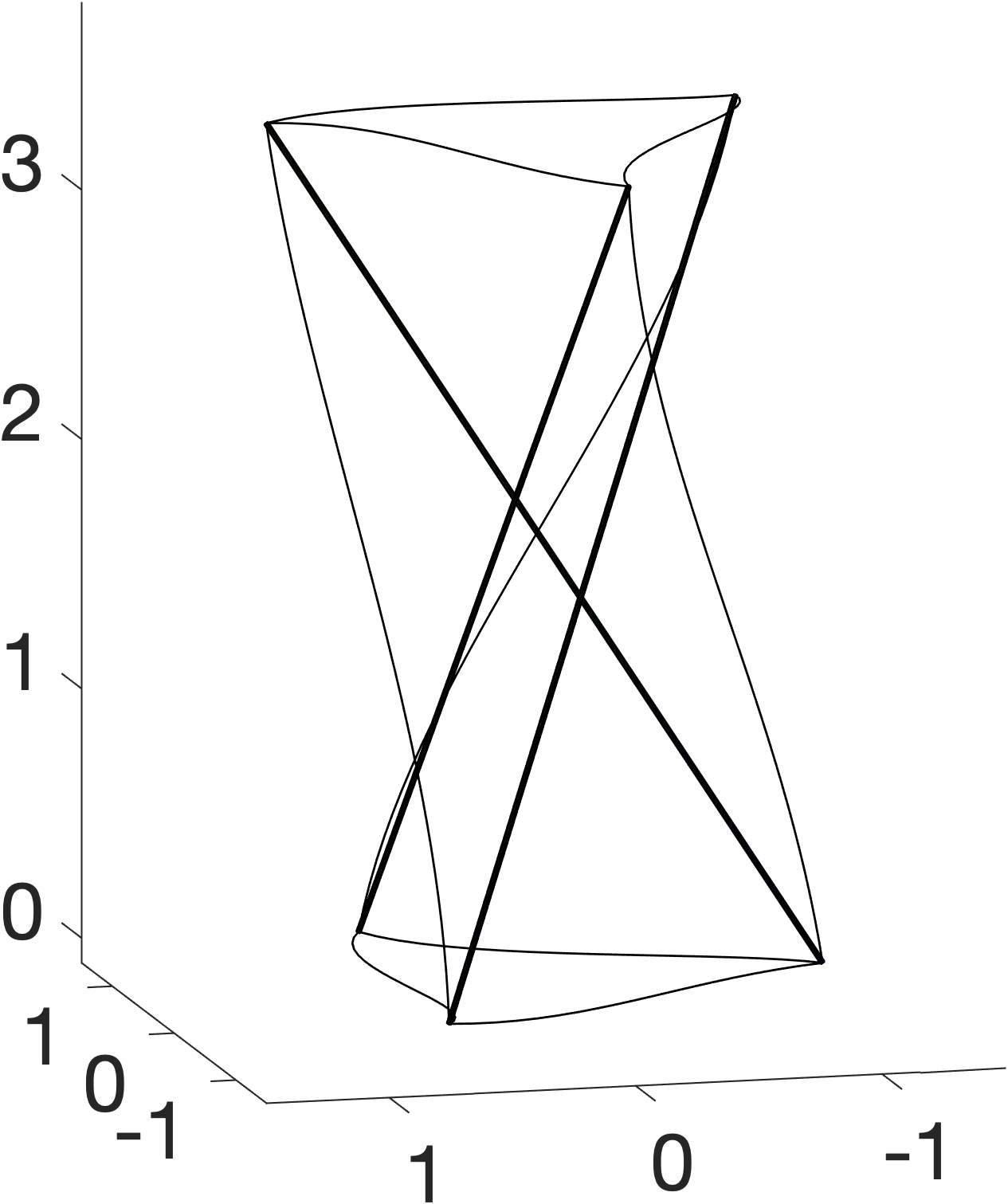}
  \caption{Mode 7 (0.3655 rad/s)}
\end{subfigure}
\begin{subfigure}{.24\textwidth}
  \centering
  \includegraphics[width=0.9\textwidth]{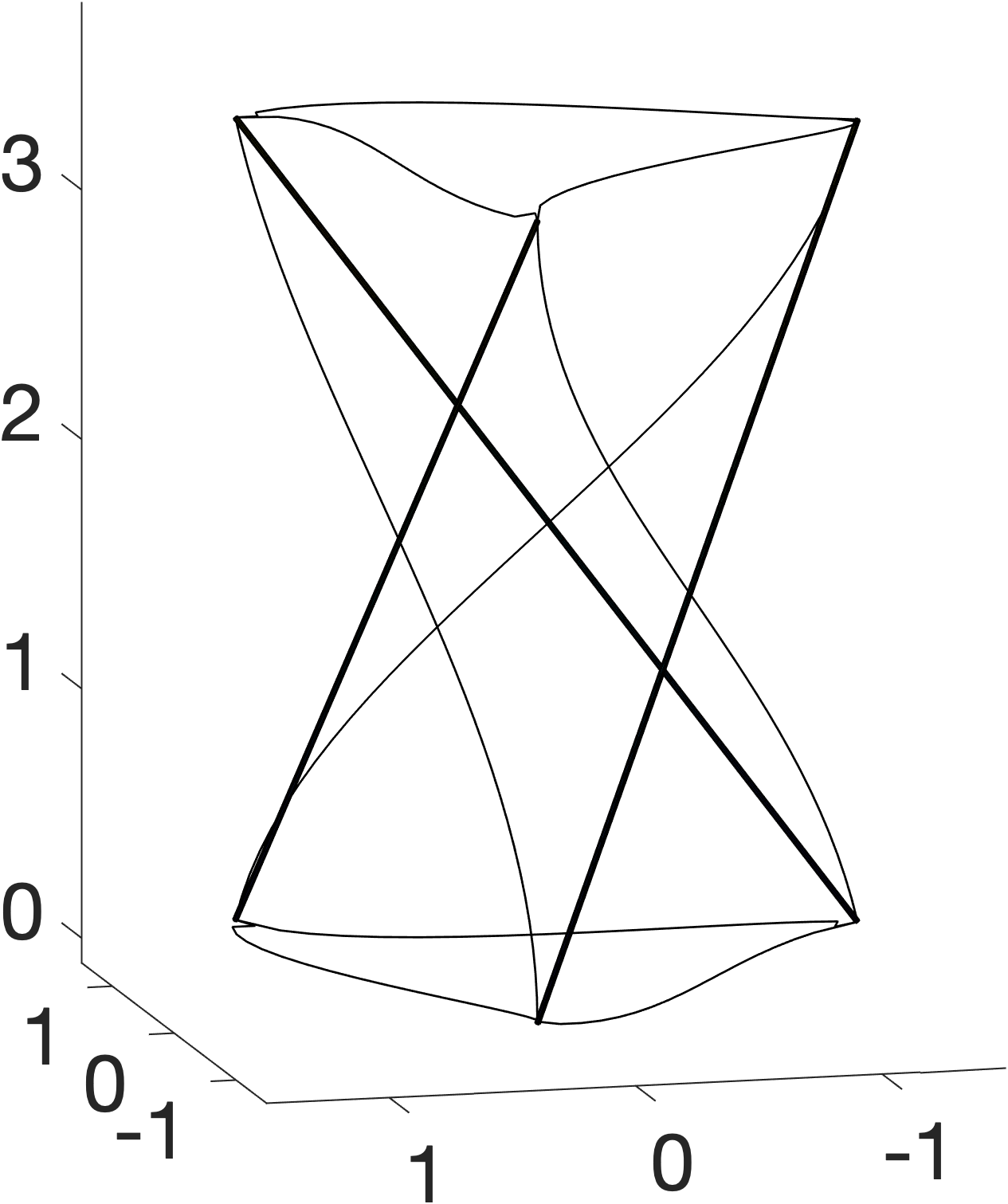}
  \caption{Mode 7 (0.3655 rad/s)}
\end{subfigure}
\begin{subfigure}{.12\textwidth}
  \centering
\end{subfigure}\\
\begin{subfigure}{.24\textwidth}
  \centering
  \includegraphics[width=0.9\textwidth]{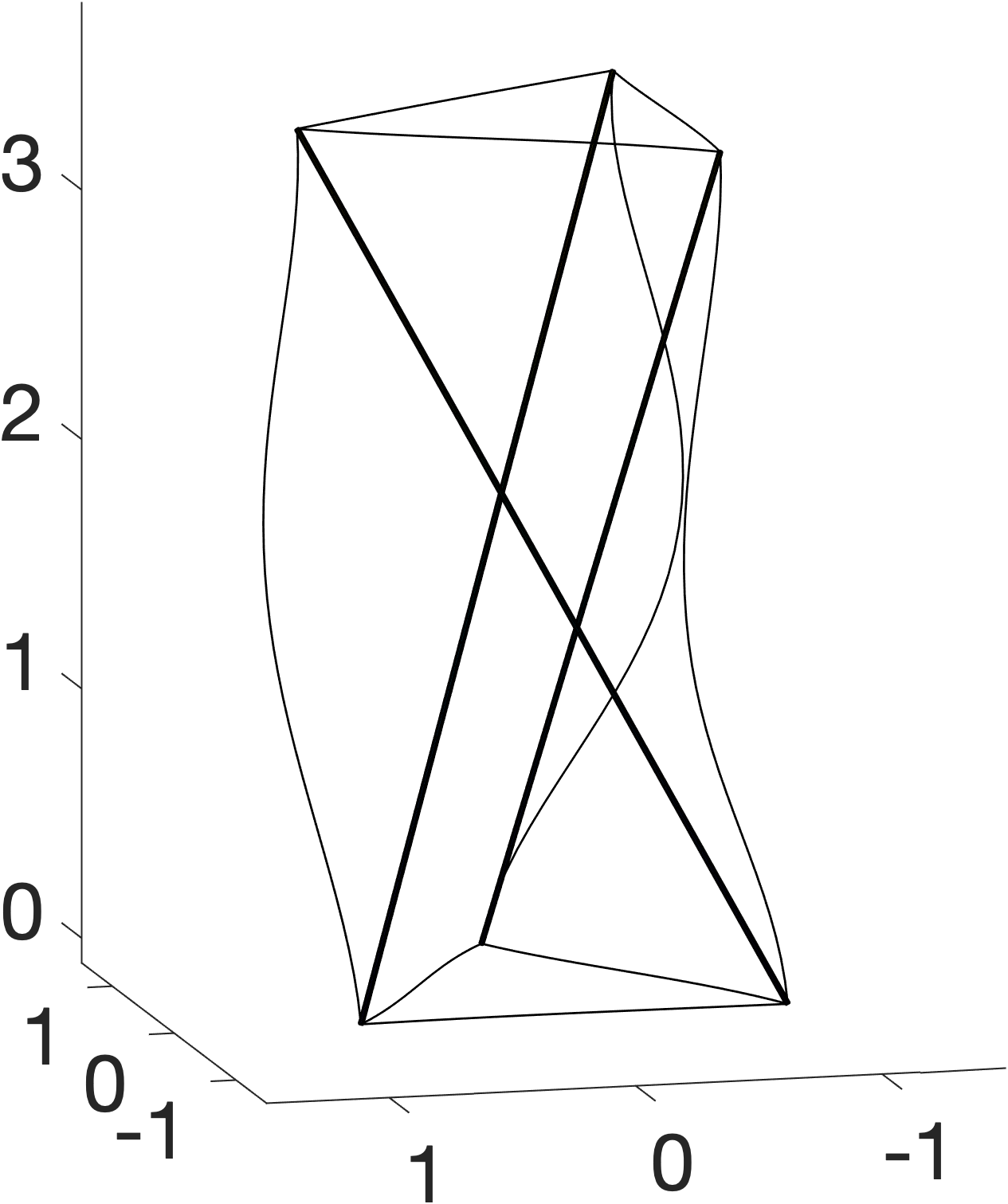}
  \caption{Mode 8 (7.6274 rad/s)}
\end{subfigure}
\begin{subfigure}{.24\textwidth}
  \centering
  \includegraphics[width=0.9\textwidth]{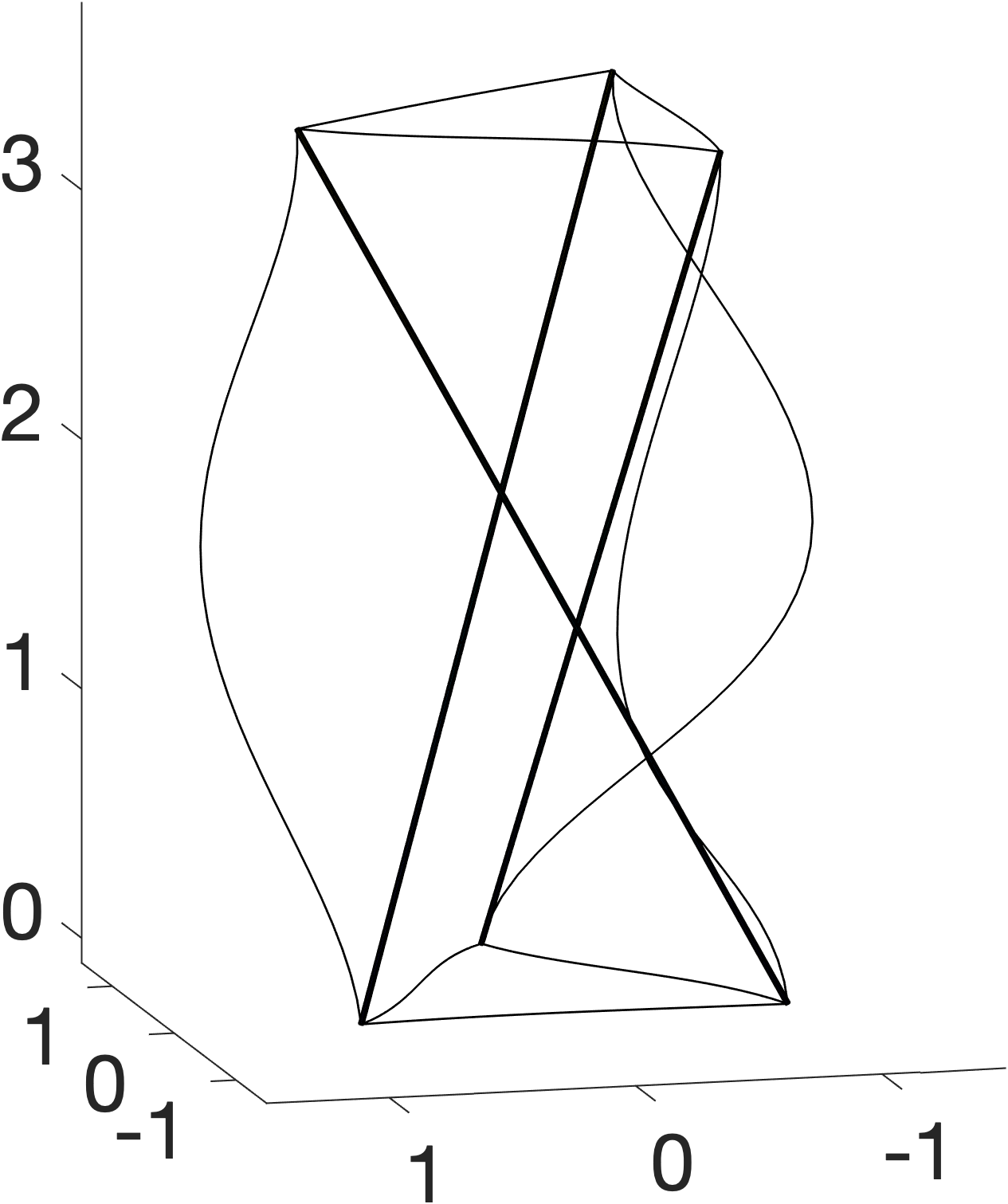}
  \caption{Mode 8 (7.6274 rad/s)}
\end{subfigure}
\begin{subfigure}{.24\textwidth}
  \centering
  \includegraphics[width=0.9\textwidth]{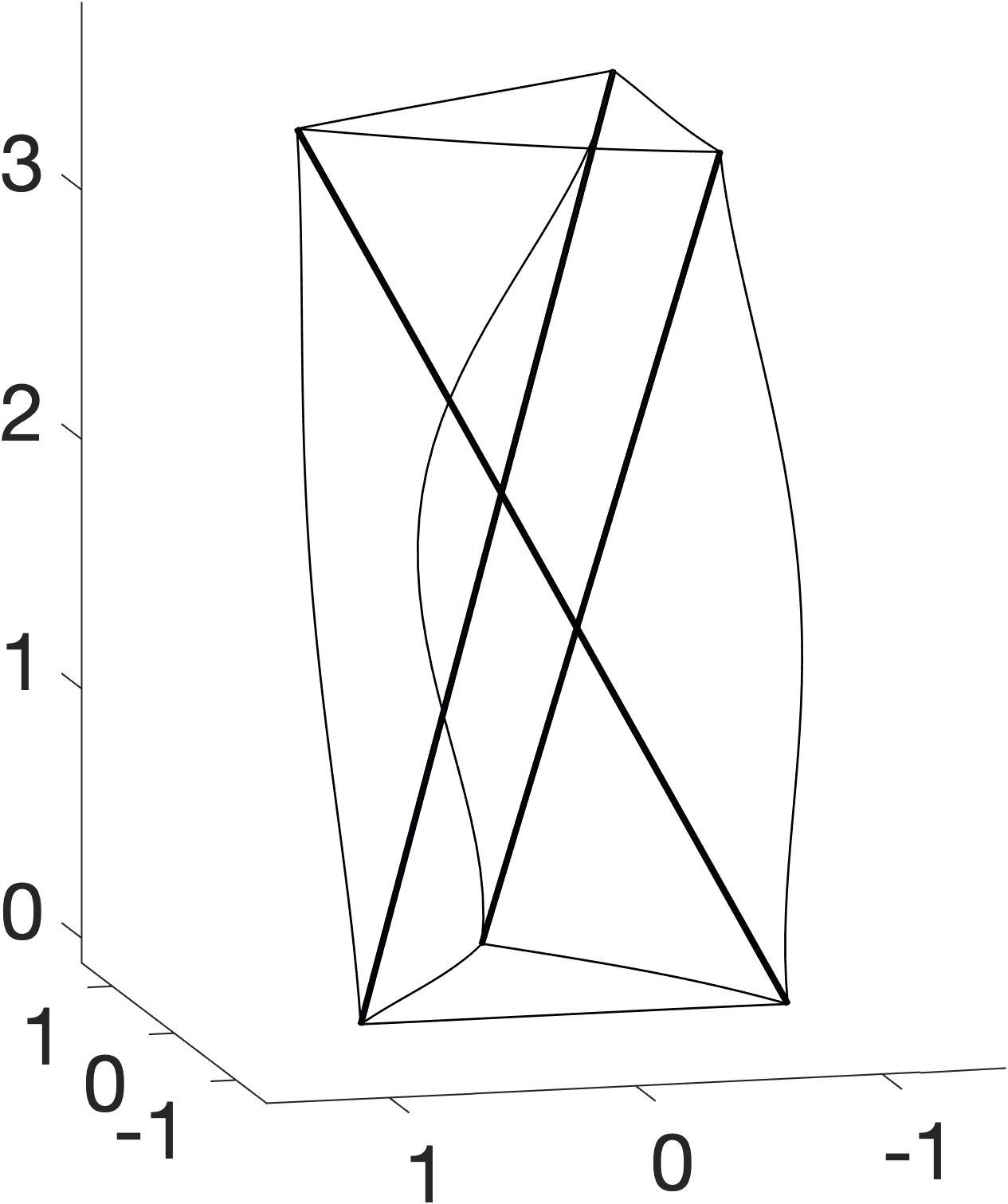}
  \caption{Mode 9 (7.7131 rad/s)}
\end{subfigure}
\begin{subfigure}{.24\textwidth}
  \centering
  \includegraphics[width=0.9\textwidth]{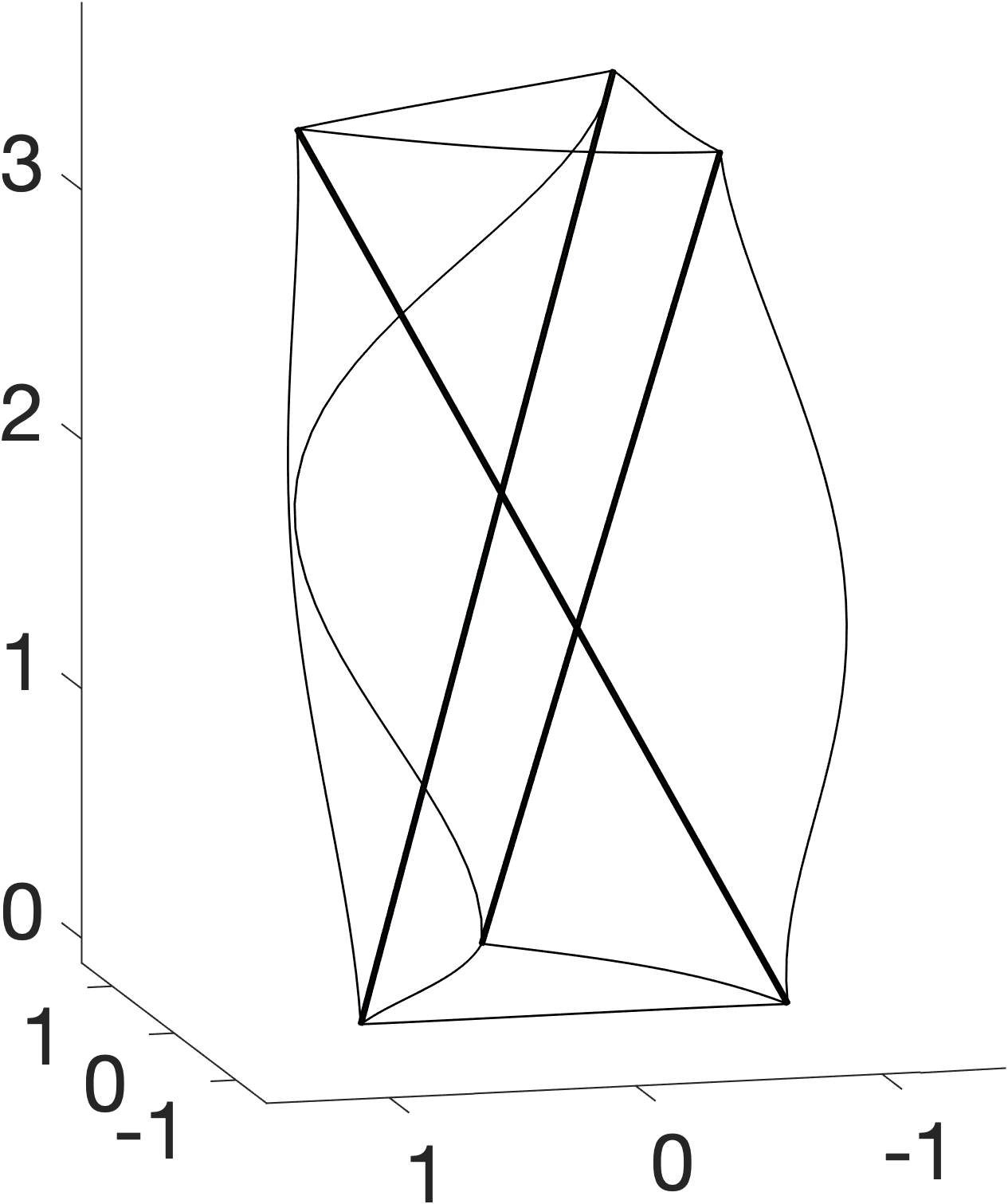}
  \caption{Mode 9 (7.7131 rad/s)}
\end{subfigure}
\caption{Visualization of different modes of deformation of the 9-segrity. No boundary condition is applied.\label{Fig:9segrityModeshapes}}
\end{figure*}

\begin{figure*}[!htb]
\centering
\begin{subfigure}{.12\textwidth}
  \centering
\end{subfigure}
\begin{subfigure}{.24\textwidth}
  \centering
  \includegraphics[width=0.9\textwidth]{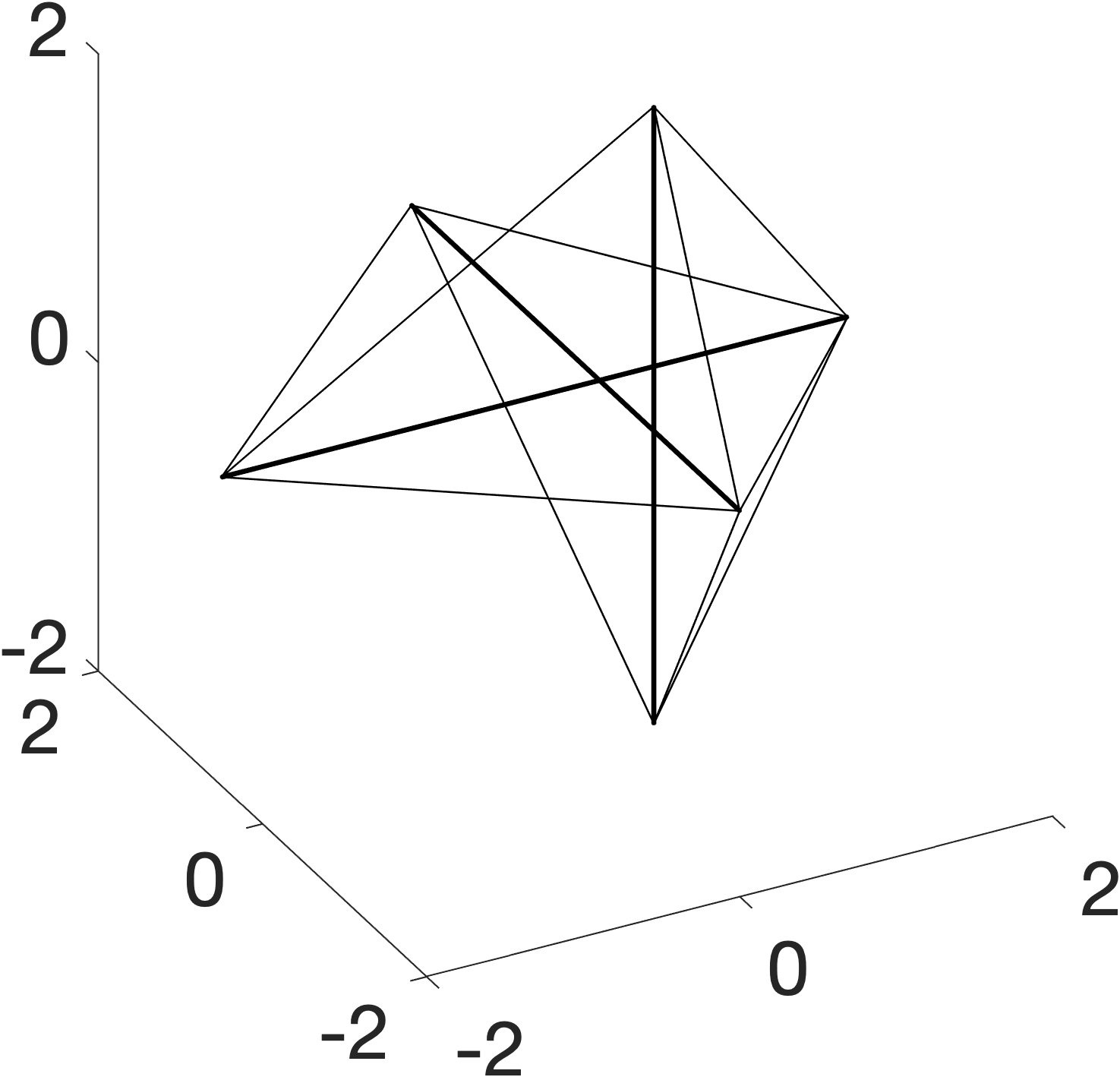}
  \caption{Original shape}
\end{subfigure}
\begin{subfigure}{.24\textwidth}
  \centering
  \includegraphics[width=0.9\textwidth]{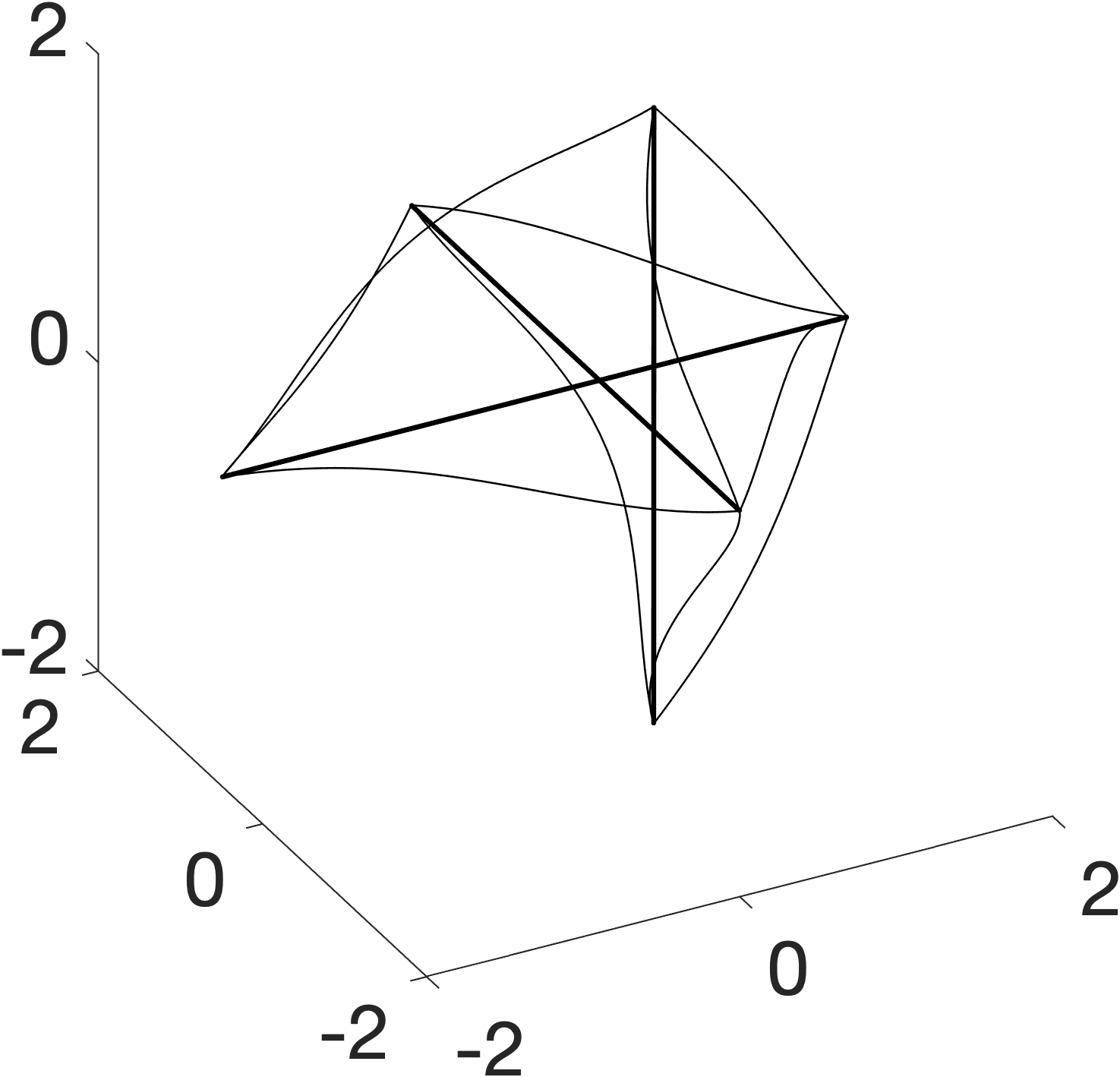}
  \caption{Mode 7 (7.1414 rad/s)}
\end{subfigure}
\begin{subfigure}{.24\textwidth}
  \centering
  \includegraphics[width=0.9\textwidth]{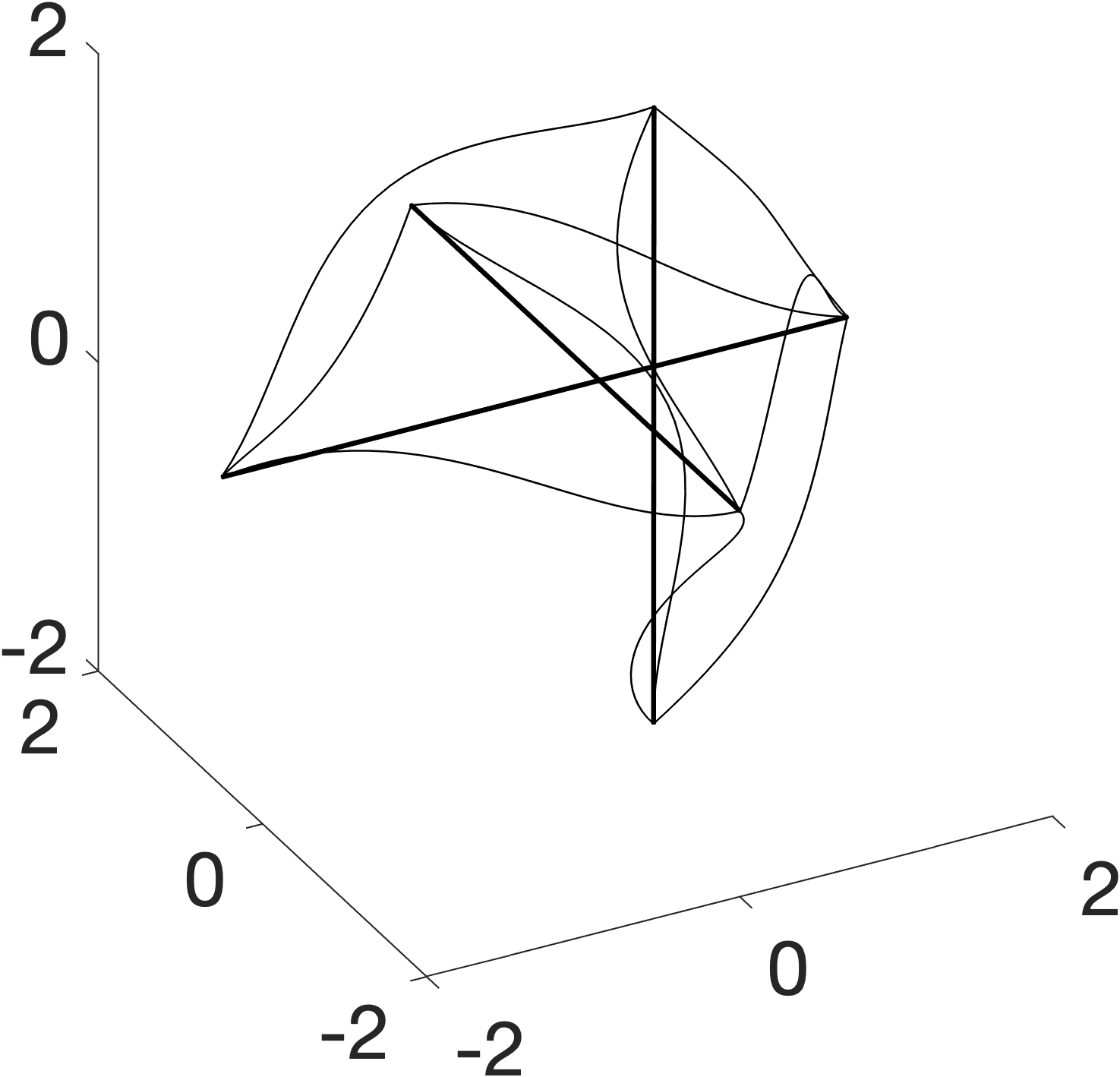}
  \caption{Mode 7 (7.1414 rad/s)}
\end{subfigure}
\begin{subfigure}{.12\textwidth}
  \centering
\end{subfigure}\\
\begin{subfigure}{.24\textwidth}
  \centering
  \includegraphics[width=0.9\textwidth]{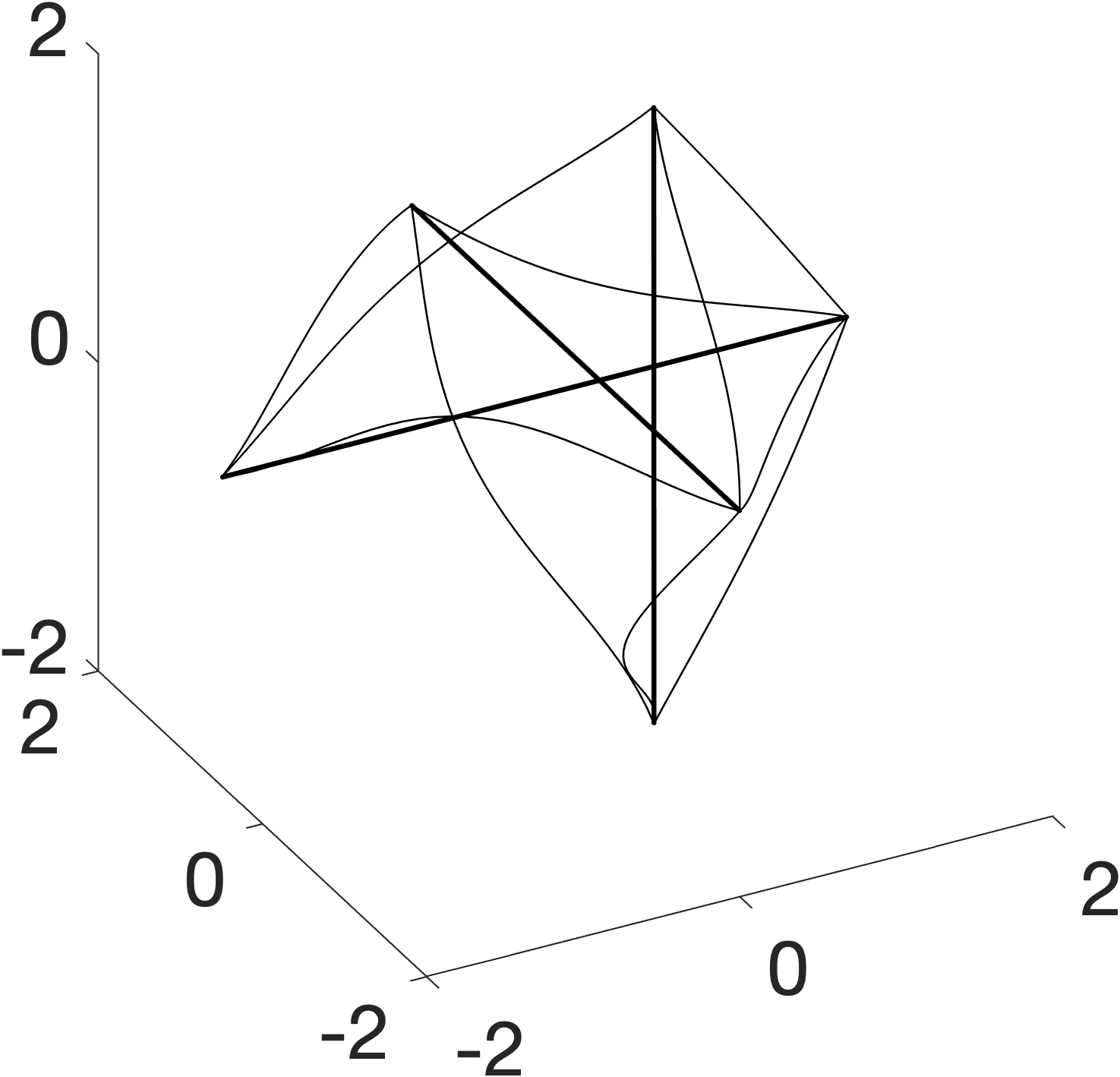}
  \caption{Mode 8 (7.5383 rad/s)}
\end{subfigure}
\begin{subfigure}{.24\textwidth}
  \centering
  \includegraphics[width=0.9\textwidth]{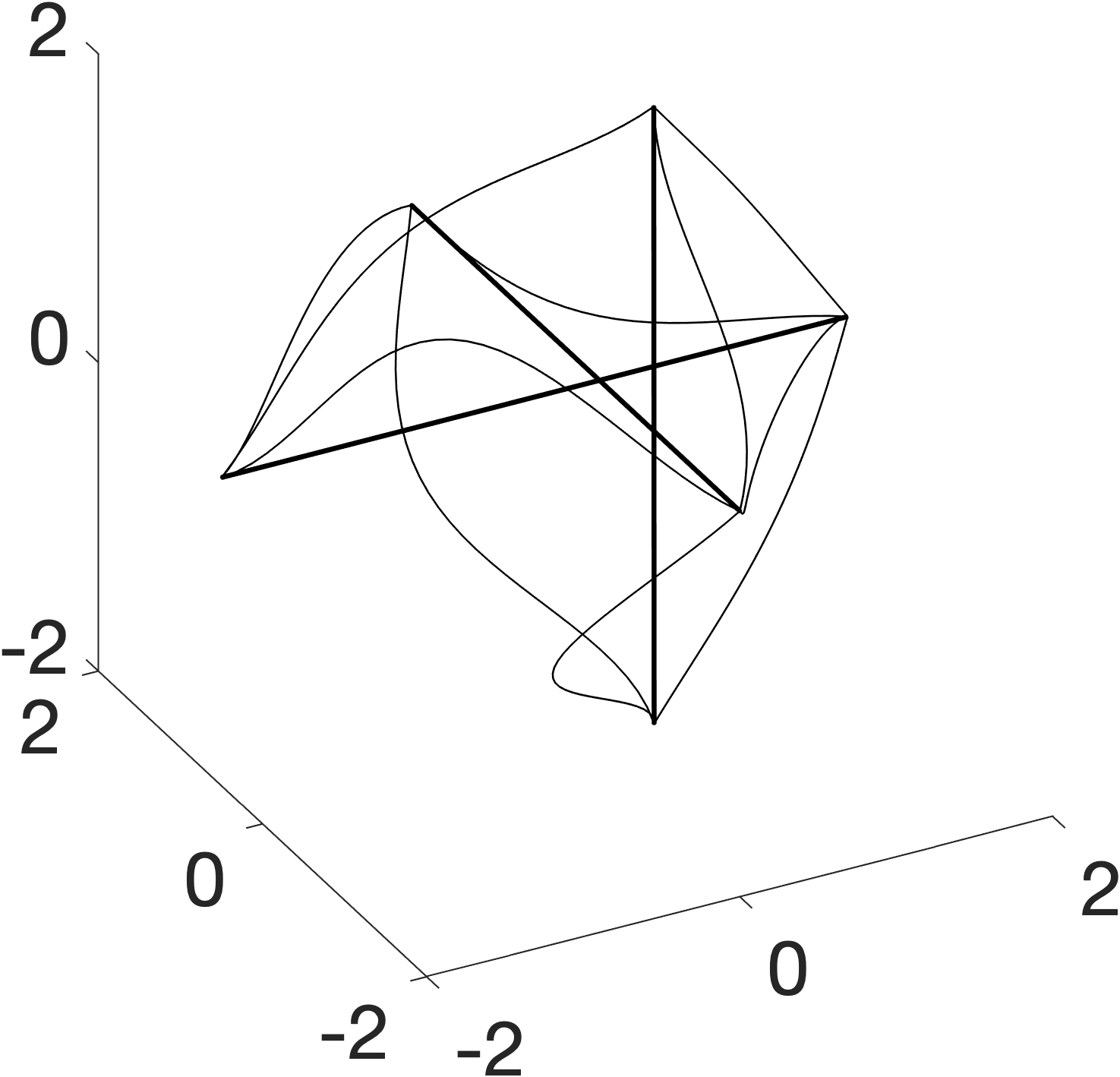}
  \caption{Mode 8 (7.5383 rad/s)}
\end{subfigure}
\begin{subfigure}{.24\textwidth}
  \centering
  \includegraphics[width=0.9\textwidth]{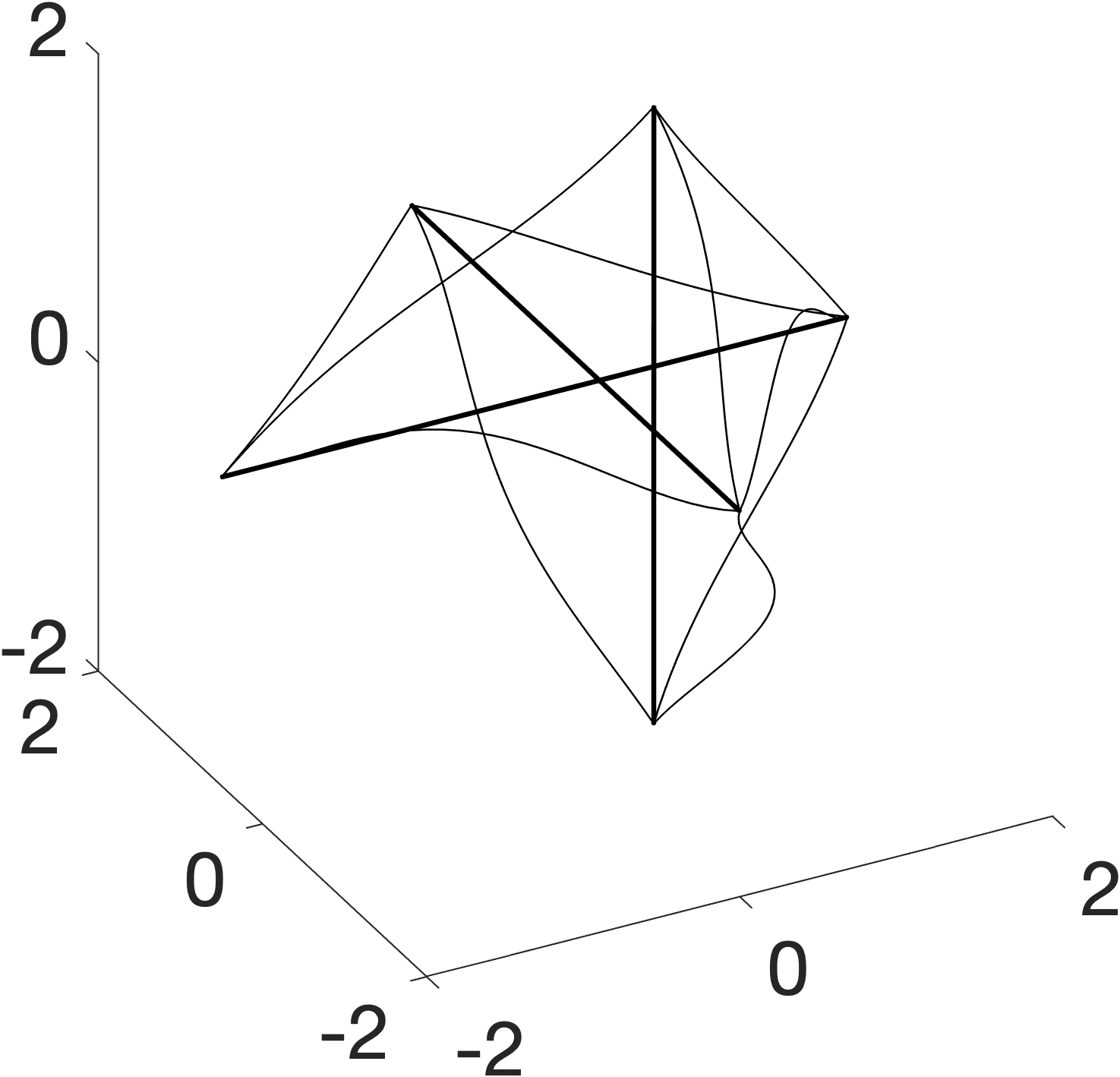}
  \caption{Mode 9 (8.1529 rad/s)}
\end{subfigure}
\begin{subfigure}{.24\textwidth}
  \centering
  \includegraphics[width=0.9\textwidth]{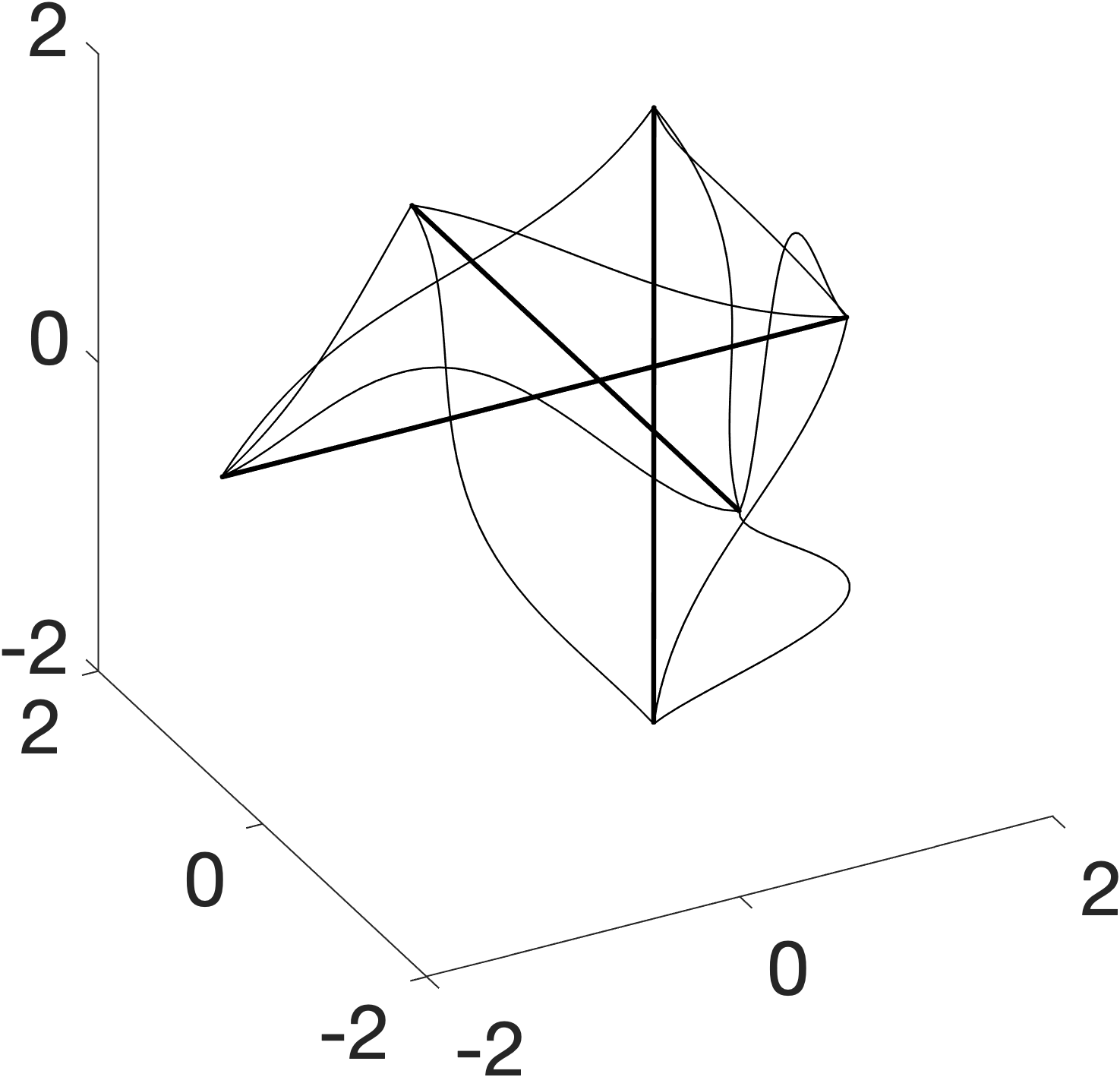}
  \caption{Mode 9 (8.1529 rad/s)}
\end{subfigure}
\caption{Visualization of different modes of deformation of the 10-segrity. No boundary condition is applied.\label{Fig:10segrityModeshapes}}
\end{figure*}

The original shape and the 7 / 8 / 9-th mode shapes for both the 9-segrity and the 10-segrity are visualized in the \Cref{Fig:9segrityModeshapes} and \Cref{Fig:10segrityModeshapes}. Here, the visualization is considered only for the case with the Youngs' modulus of Kevlar being 112 GPa. The Mode 7 for the 9-segrity shows the famous swinging mode, originally discussed by Connelly and co-workers. However, the other higher modes relate direction to individual motion of the string elements rather than the overall structure itself. In contrast, the 10-segrity proposed in this work shows no sign of a soft mode similar to that in the 9-segrity. 

It is further important to note here that once the strings are slack, they are under a compressive load and thus do not have any stiffness. The whole concept of stability of a tensegrity is based on the idea that the strings are always in tension and rods in compression. Once the strings are in compression, they have zero stiffness and they can deform and maintain a deformed configuration without any external work. This also means that they are in continuous equilibrium over a finite range of motion. This is also what can be physically seen that once a string is slack, they can occupy several positions. The positions for the strings shown in higher modes would be only one representation of the several finite shapes that the string can attain.

However, the primary noteworthy point here is that at the first eigenmode of deformation, the slackness of the string can result in the structural instability in the 9-string rather than the 10-string configuration. There are no strings available to prevent the swinging mode of deformation. In contrast, the first eigenmode of the 10-string configuration requires much higher energy and can be prevented by usage of stiff strings and no-prestressing.

\subsubsection{Modal analysis (with constraints)}
\begin{figure*}[!htb]
\centering
\begin{subfigure}{.24\textwidth}
  \centering
  \includegraphics[width=0.9\textwidth]{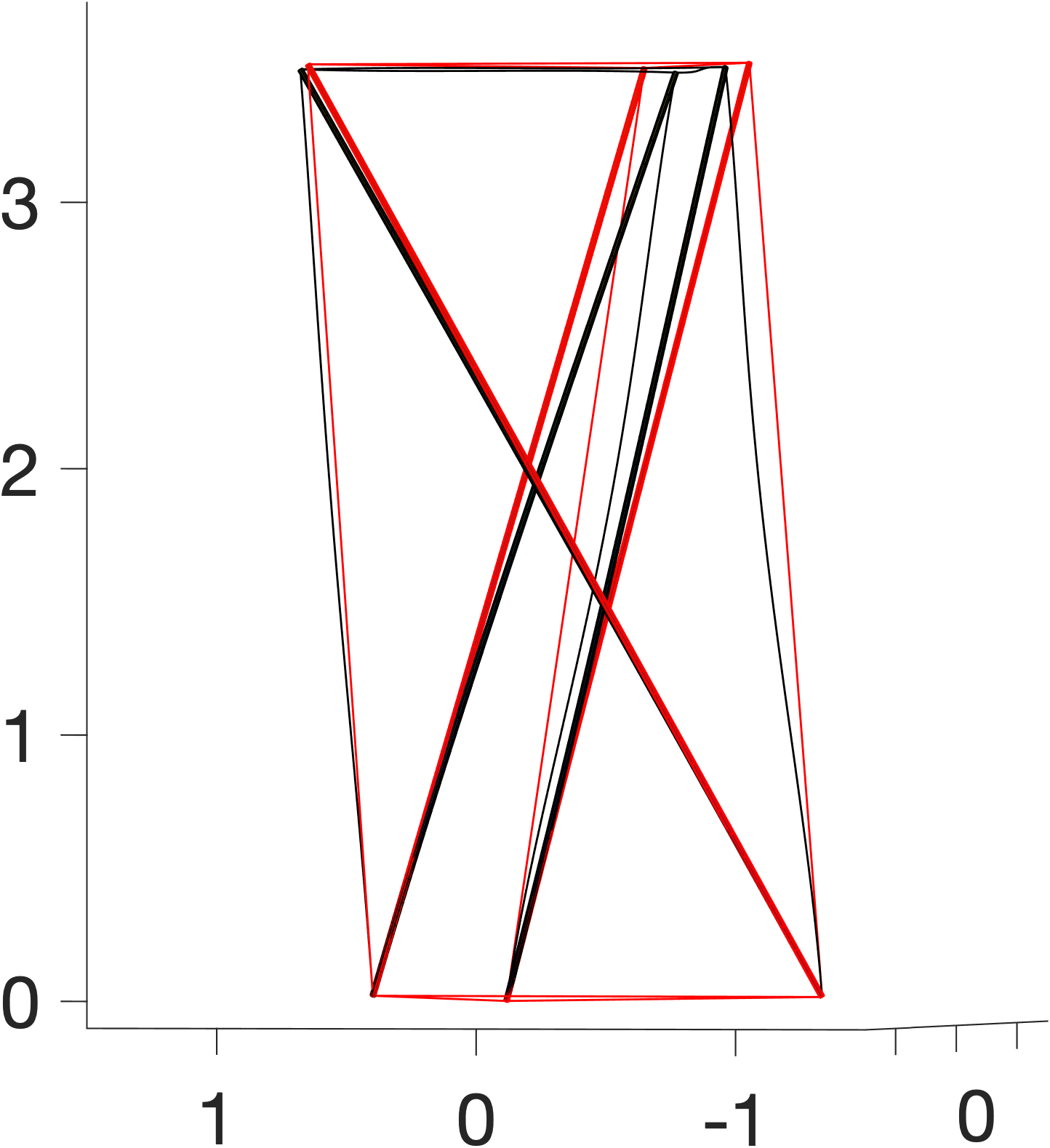}
\end{subfigure}
\begin{subfigure}{.24\textwidth}
  \centering
  \includegraphics[width=0.9\textwidth]{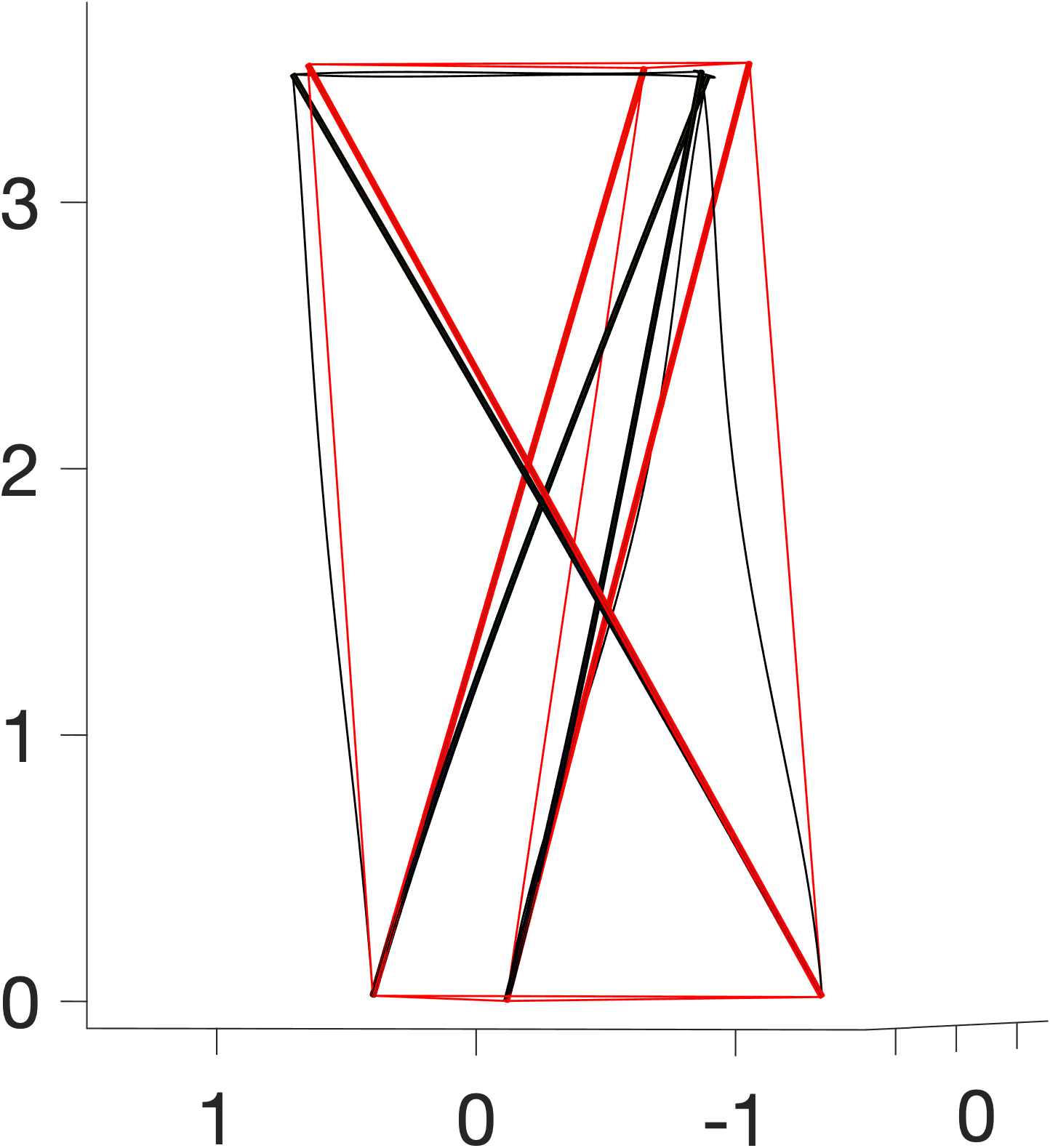}
\end{subfigure}
\begin{subfigure}{.24\textwidth}
  \centering
  \includegraphics[width=0.9\textwidth]{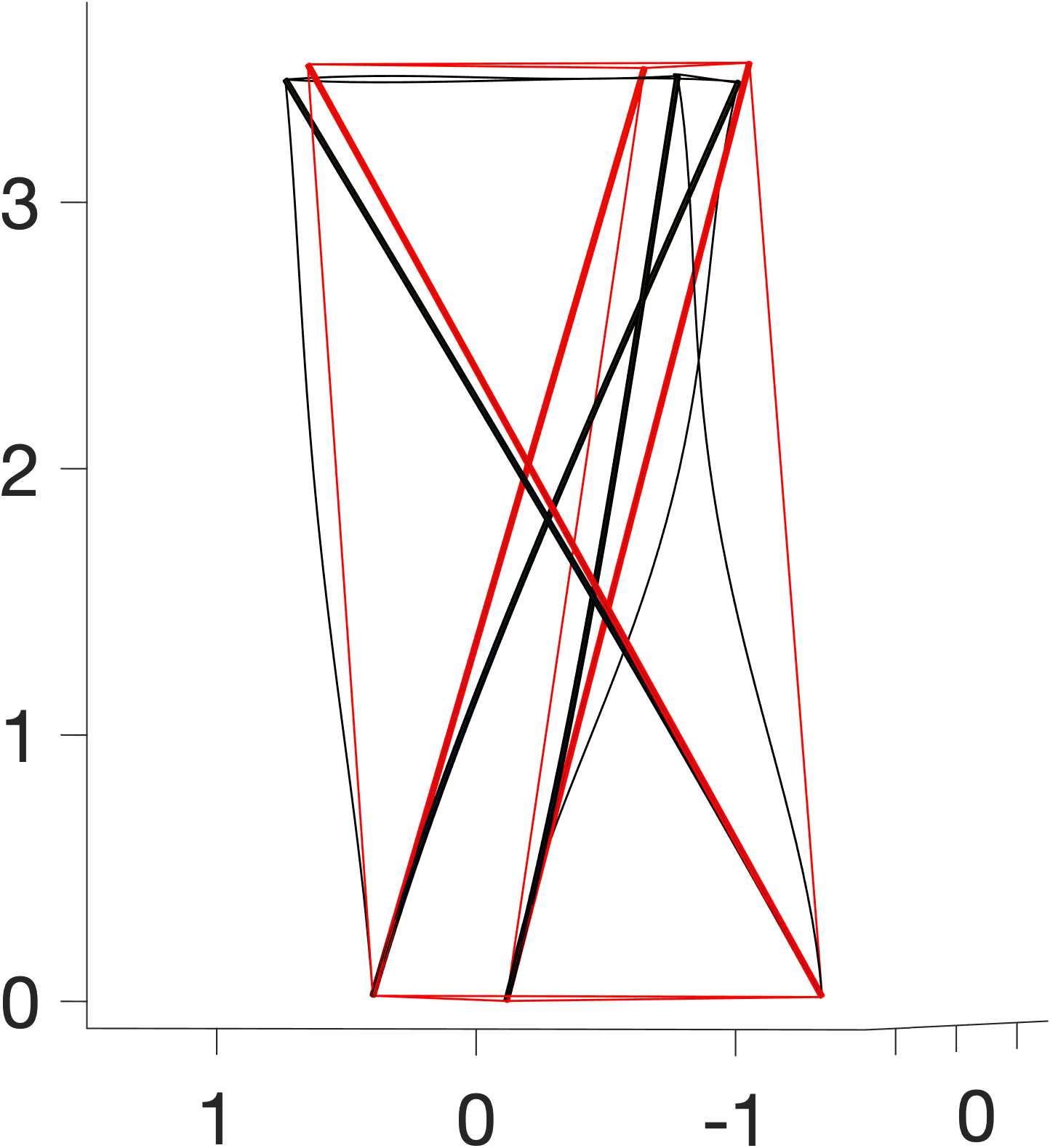}
\end{subfigure}
\begin{subfigure}{.24\textwidth}
  \centering
  \includegraphics[width=0.9\textwidth]{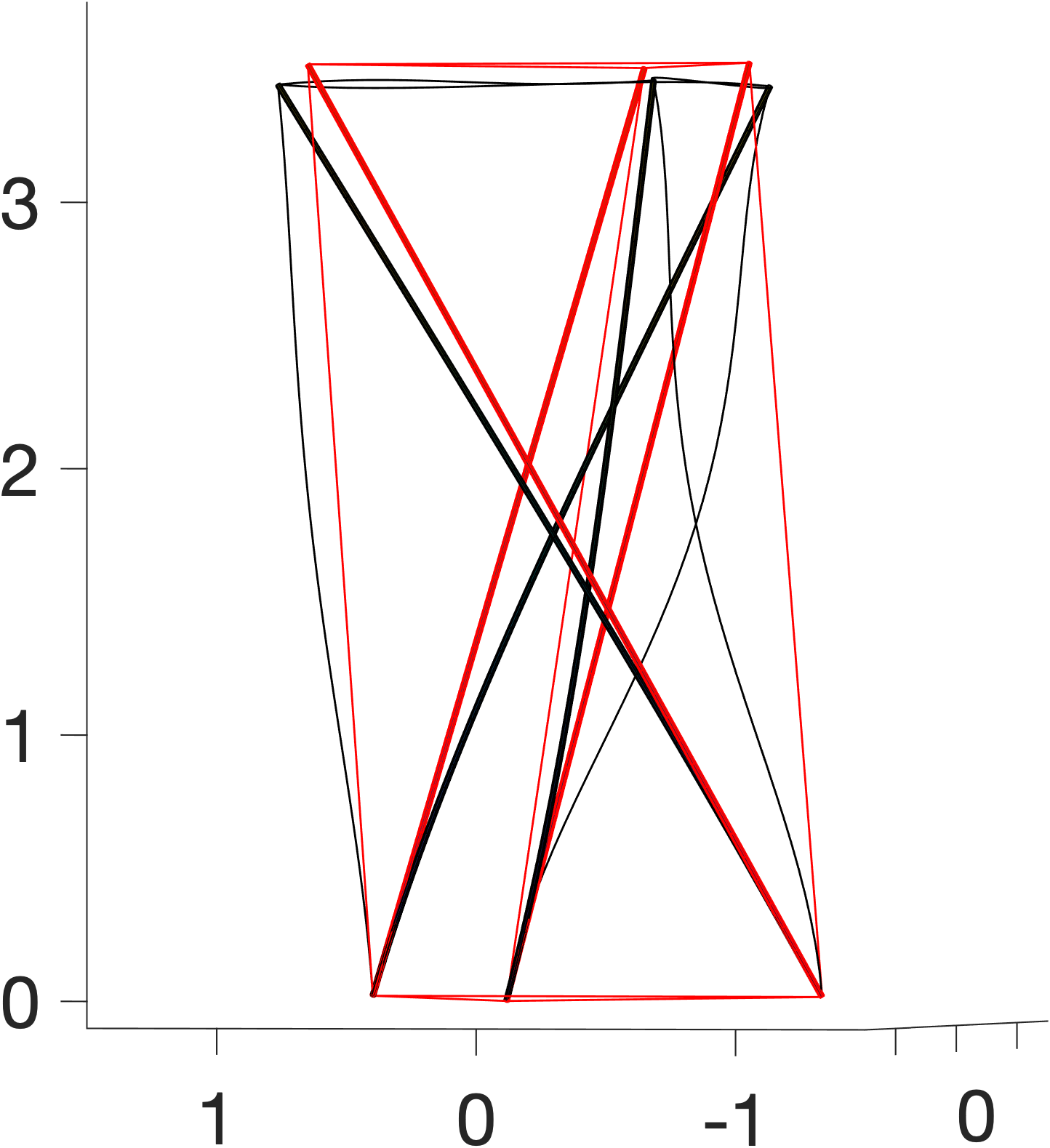}
\end{subfigure}
\begin{subfigure}{.24\textwidth}
  \centering
  \includegraphics[width=0.9\textwidth]{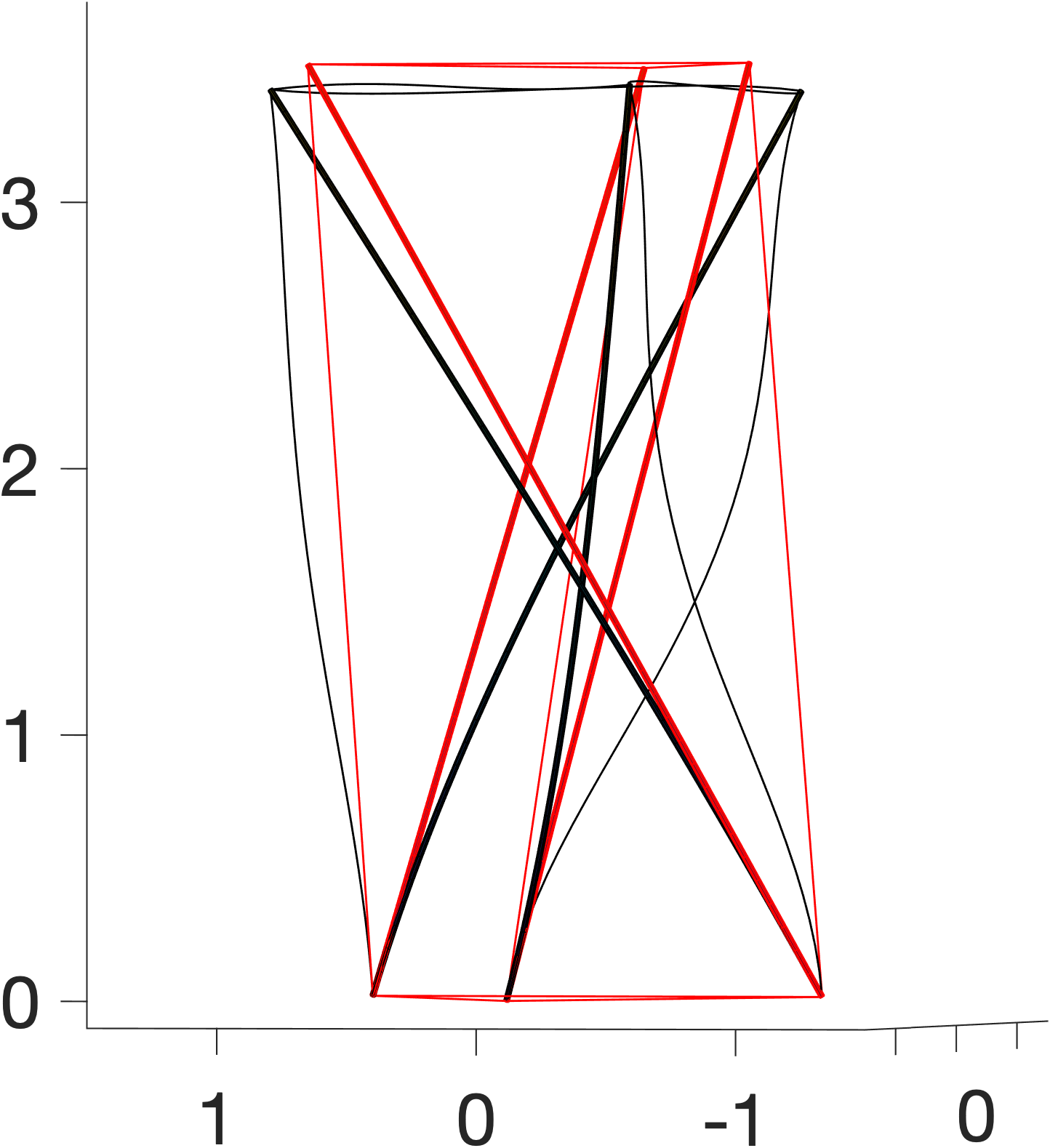}
\end{subfigure}
\begin{subfigure}{.24\textwidth}
  \centering
  \includegraphics[width=0.9\textwidth]{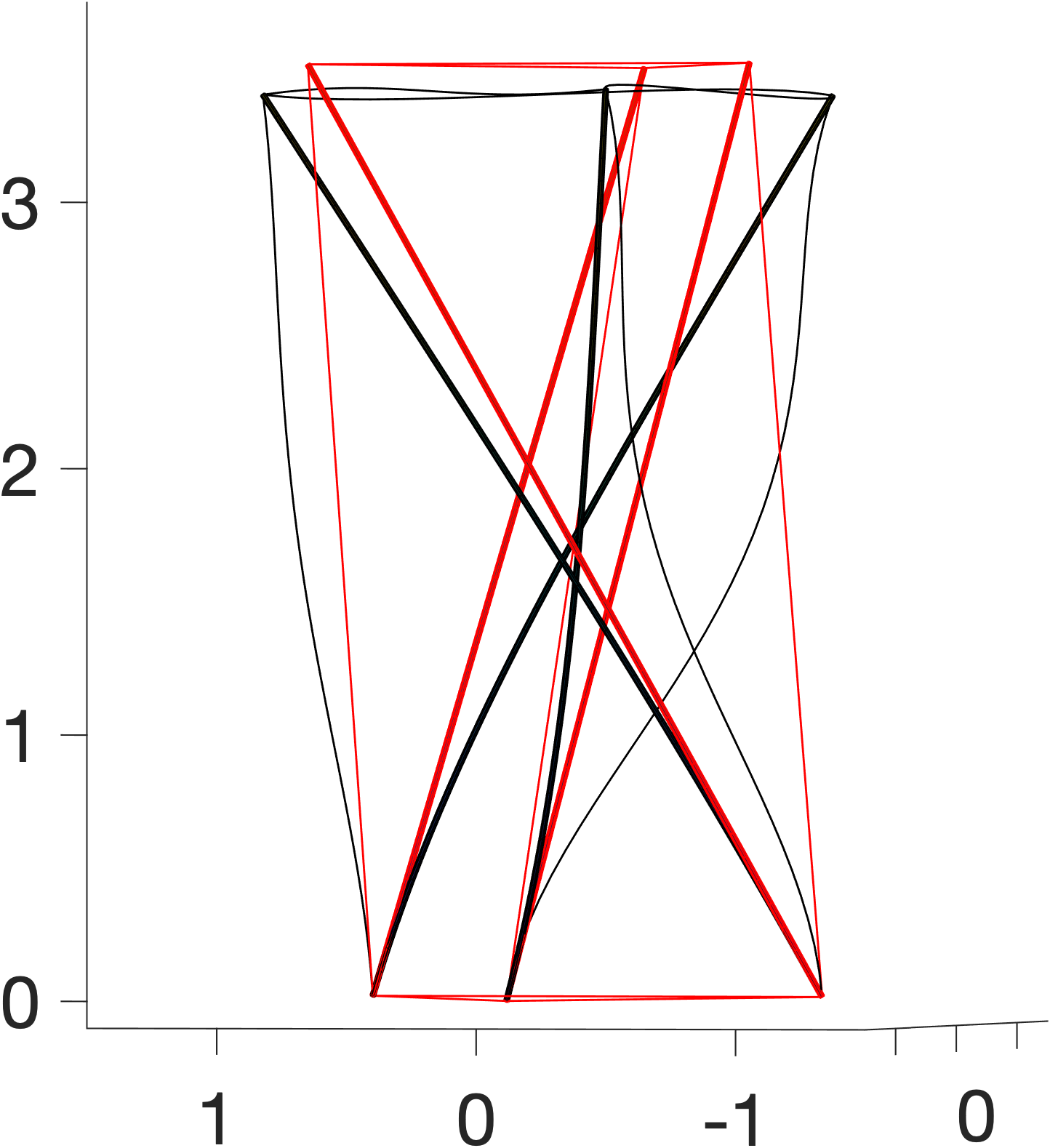}
\end{subfigure}
\begin{subfigure}{.24\textwidth}
  \centering
  \includegraphics[width=0.9\textwidth]{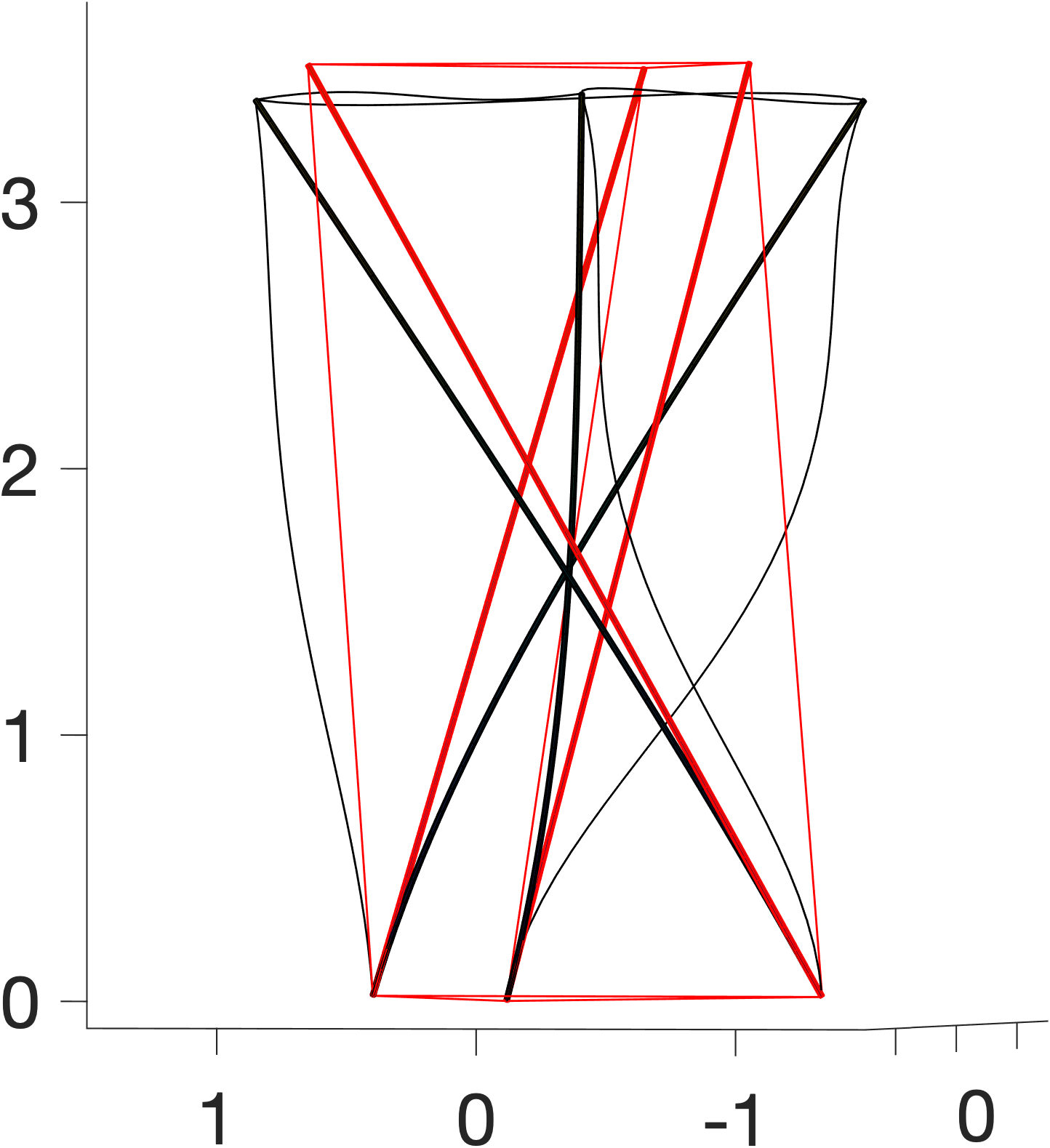}
\end{subfigure}
\begin{subfigure}{.24\textwidth}
  \centering
  \includegraphics[width=0.9\textwidth]{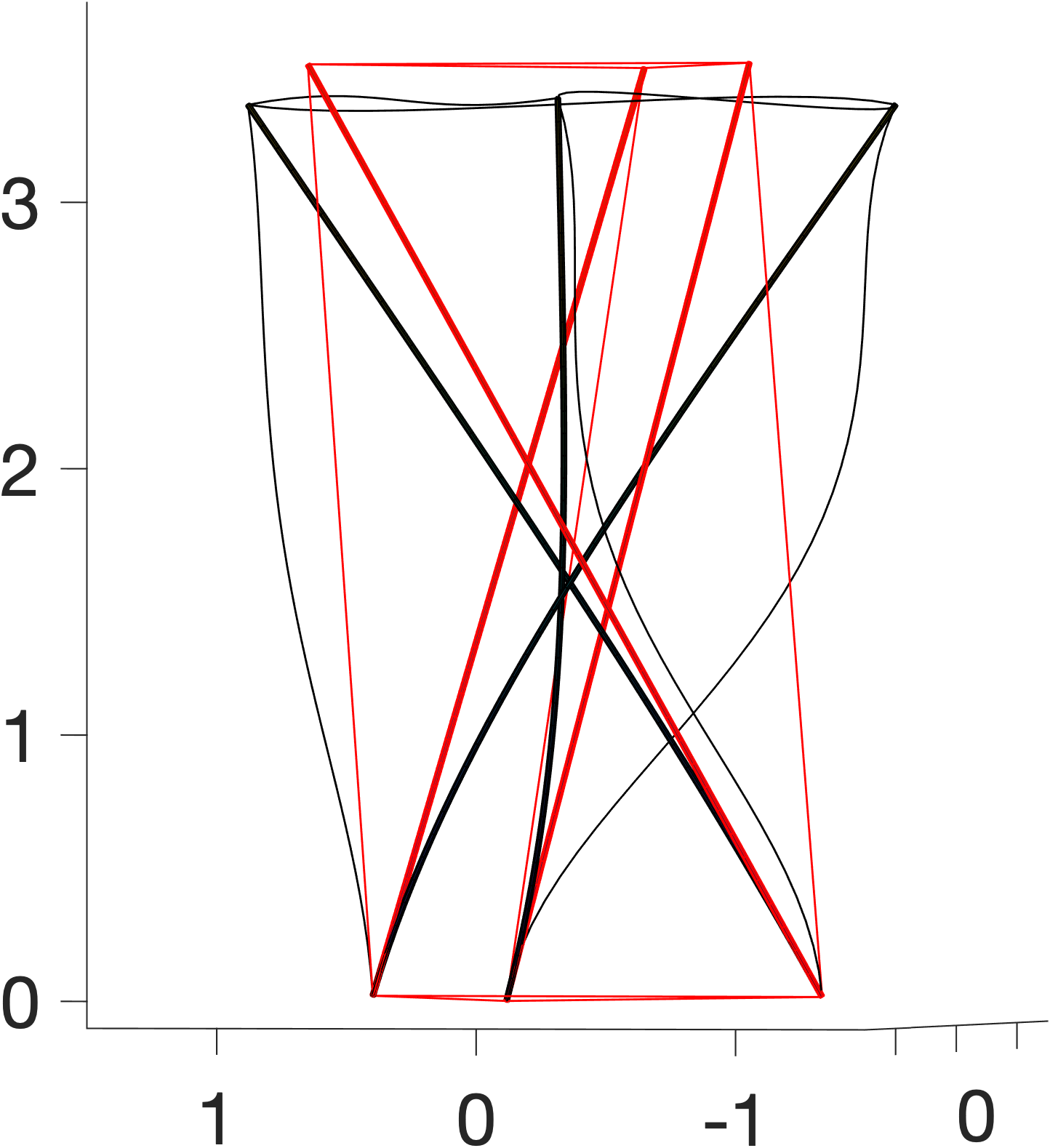}
\end{subfigure}
\caption{Visualization of the lowest mode (also the swinging mode at 0.3655 rad/s) of deformation of the 9-segrity. One end of each rod is fixed.\label{Fig:9segritySwing}}
\end{figure*}

\begin{figure*}[!htb]
\centering
\begin{subfigure}{.24\textwidth}
  \centering
  \includegraphics[width=0.9\textwidth]{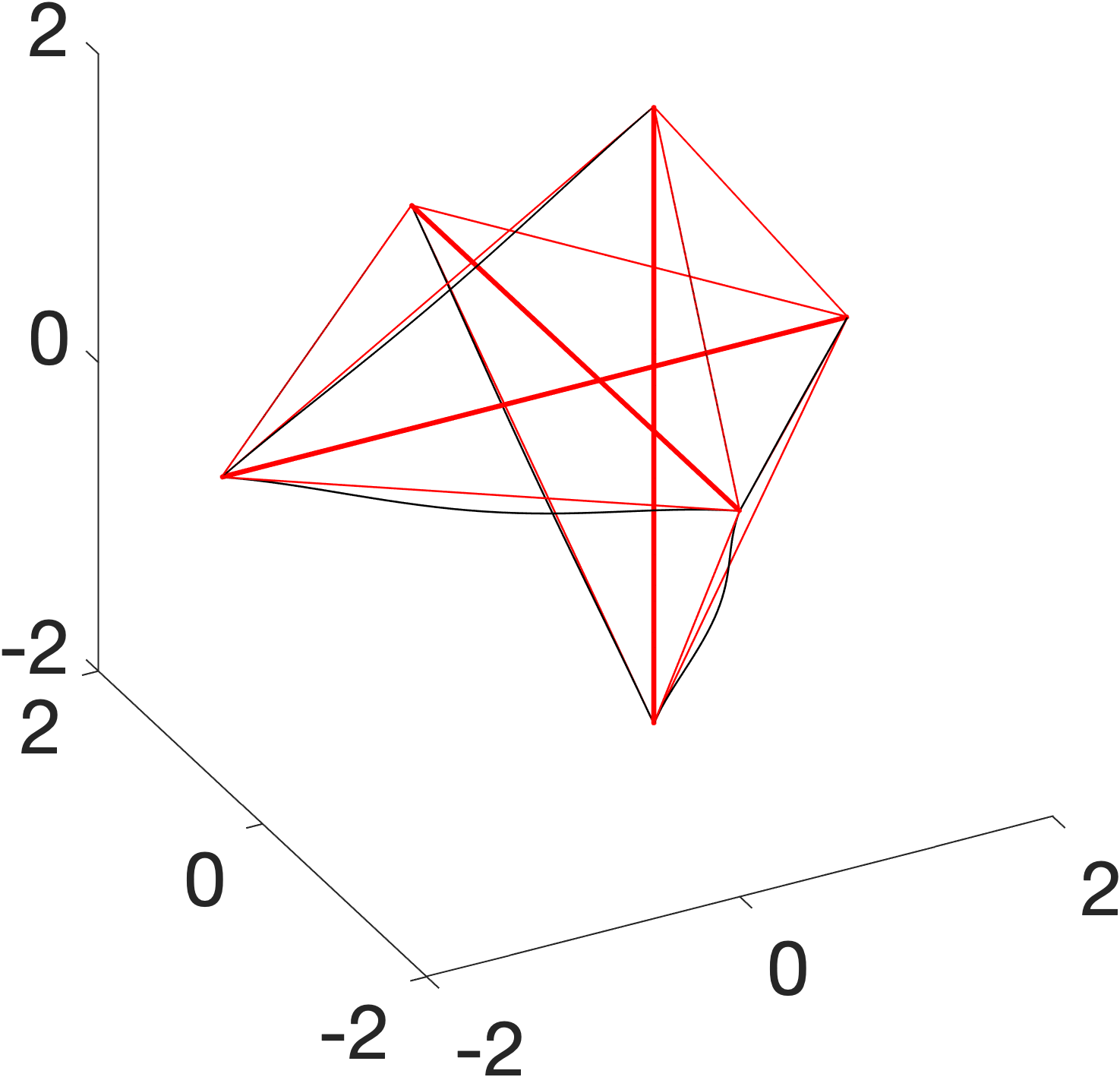}
\end{subfigure}
\begin{subfigure}{.24\textwidth}
  \centering
  \includegraphics[width=0.9\textwidth]{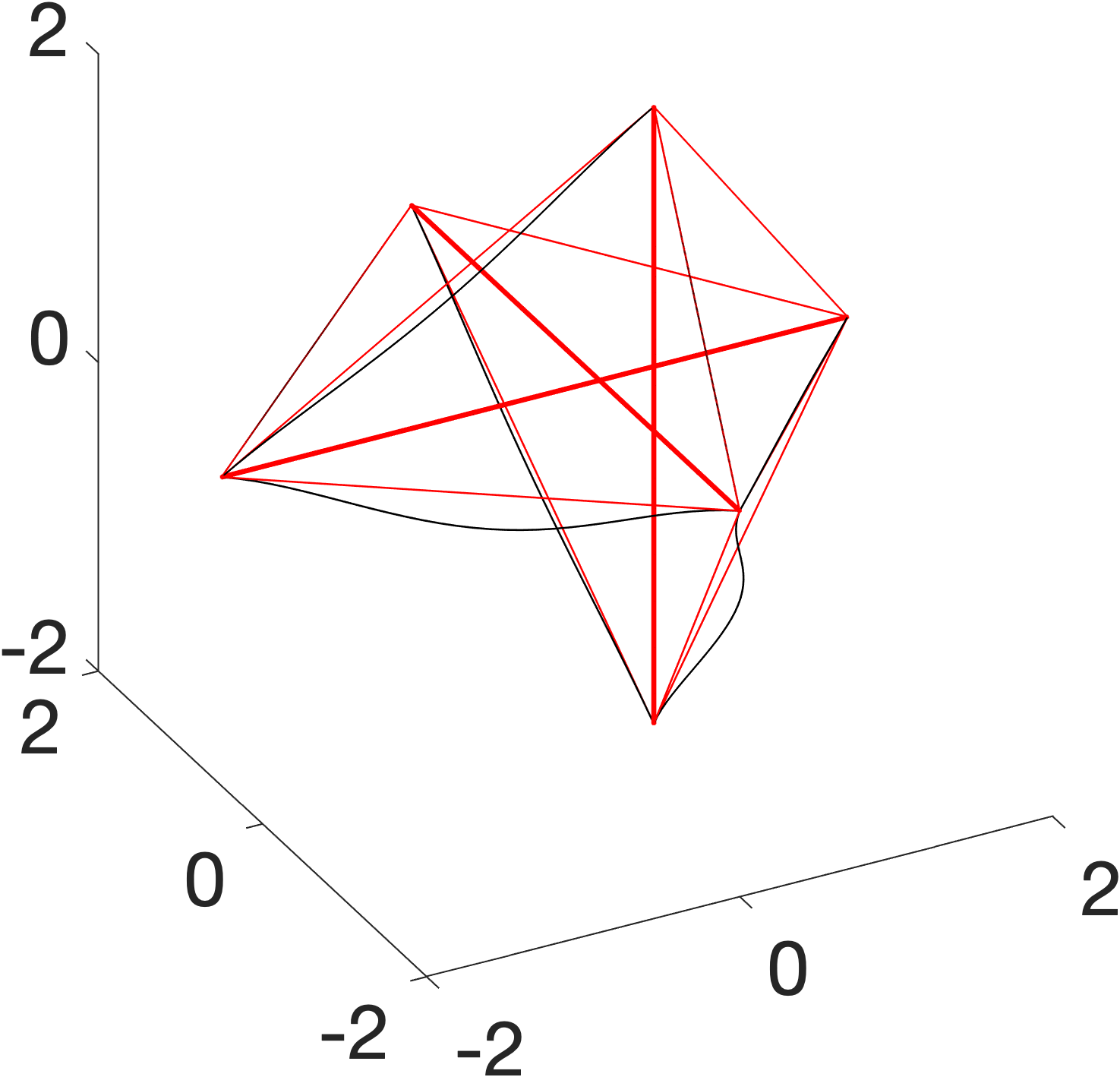}
\end{subfigure}
\begin{subfigure}{.24\textwidth}
  \centering
  \includegraphics[width=0.9\textwidth]{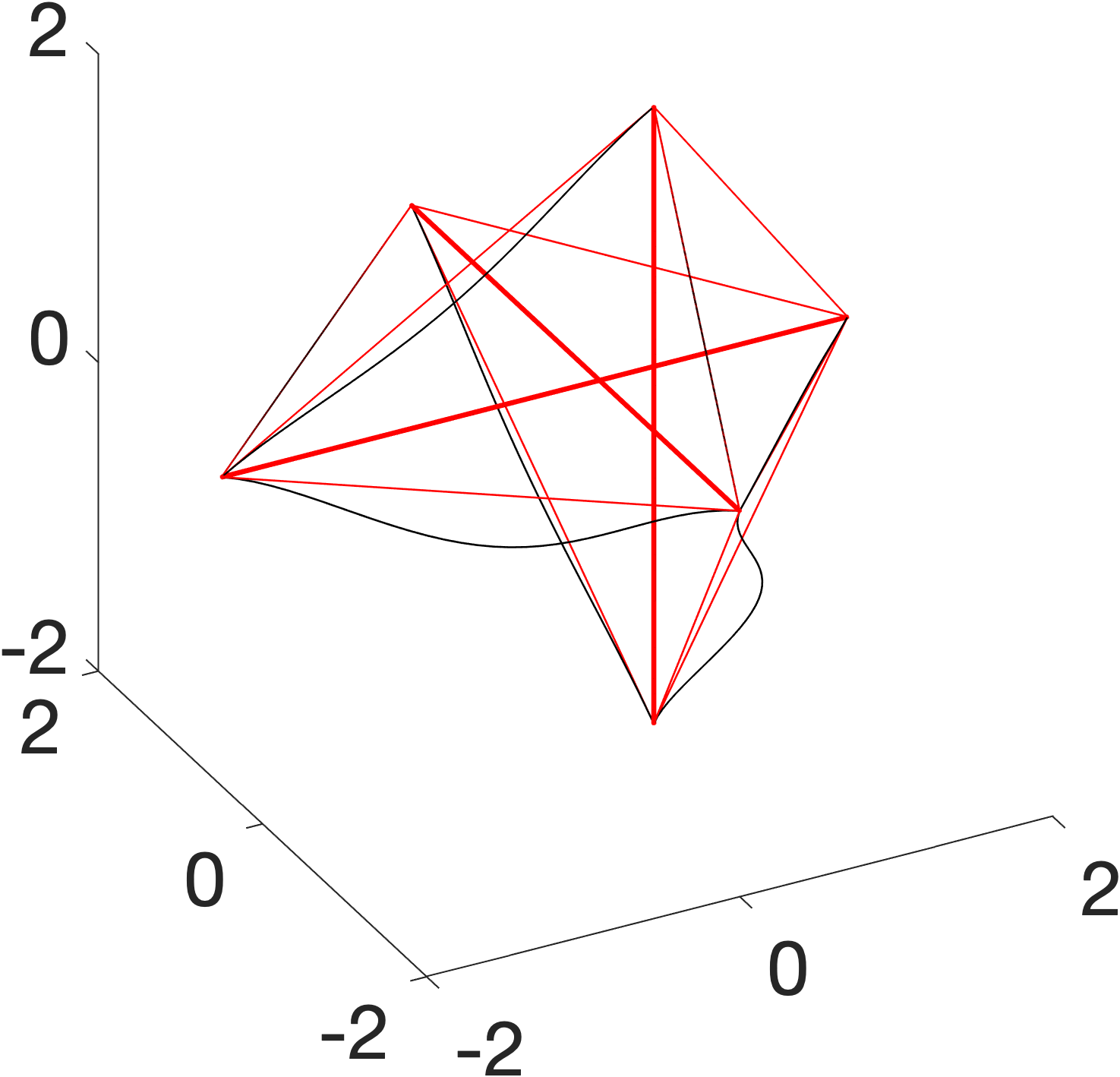}
\end{subfigure}
\begin{subfigure}{.24\textwidth}
  \centering
  \includegraphics[width=0.9\textwidth]{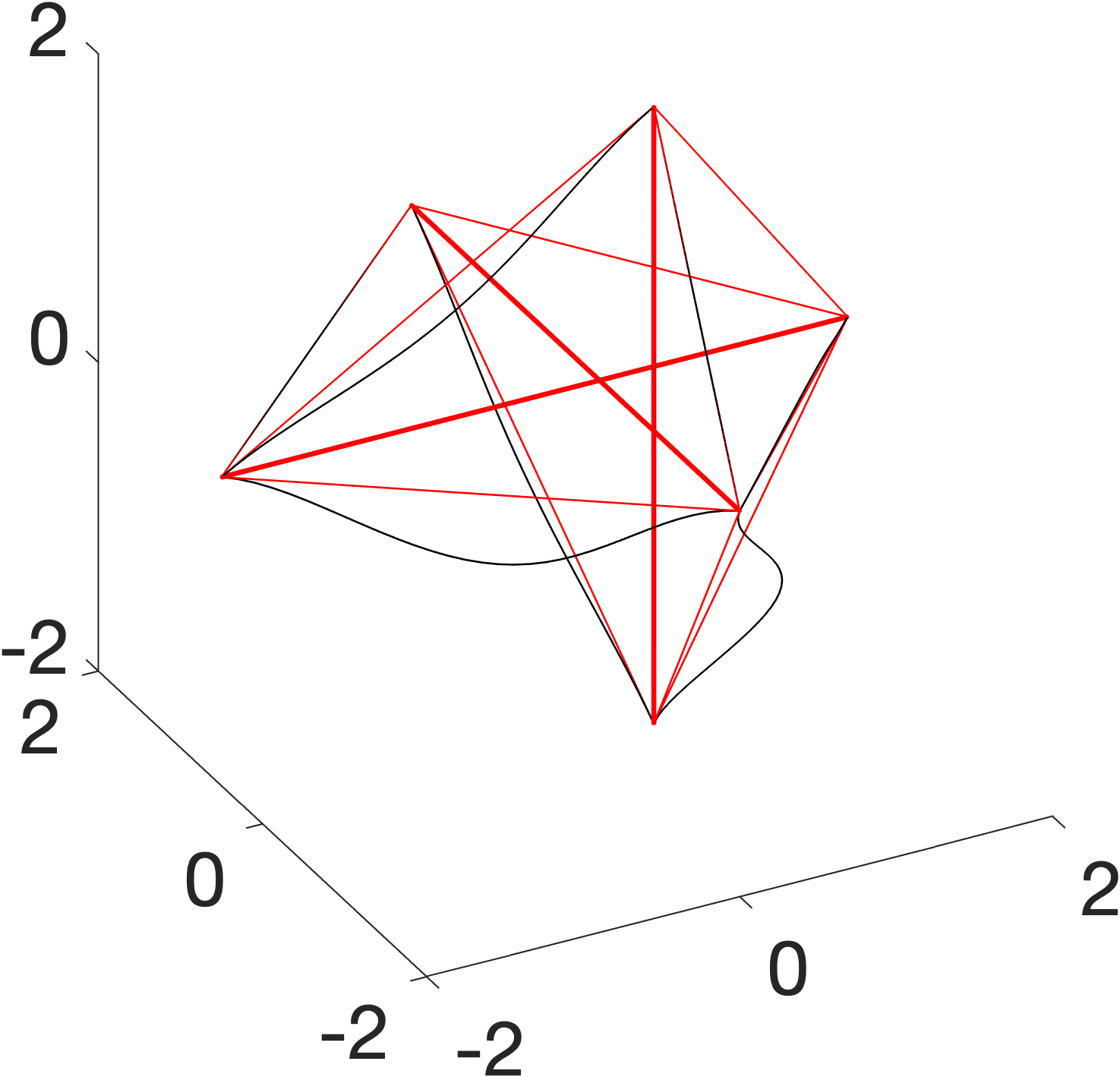}
\end{subfigure}
\begin{subfigure}{.24\textwidth}
  \centering
  \includegraphics[width=0.9\textwidth]{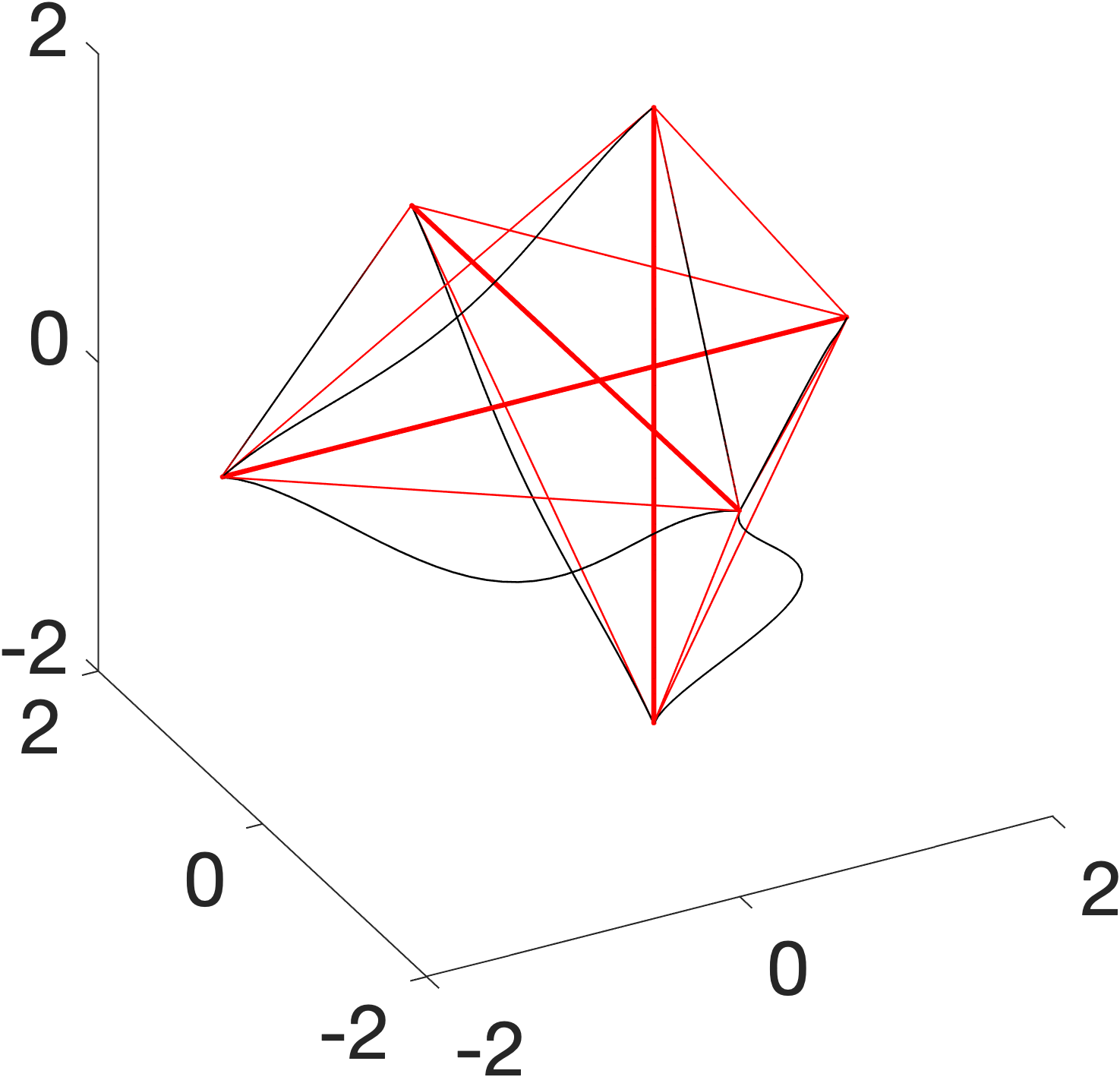}
\end{subfigure}
\begin{subfigure}{.24\textwidth}
  \centering
  \includegraphics[width=0.9\textwidth]{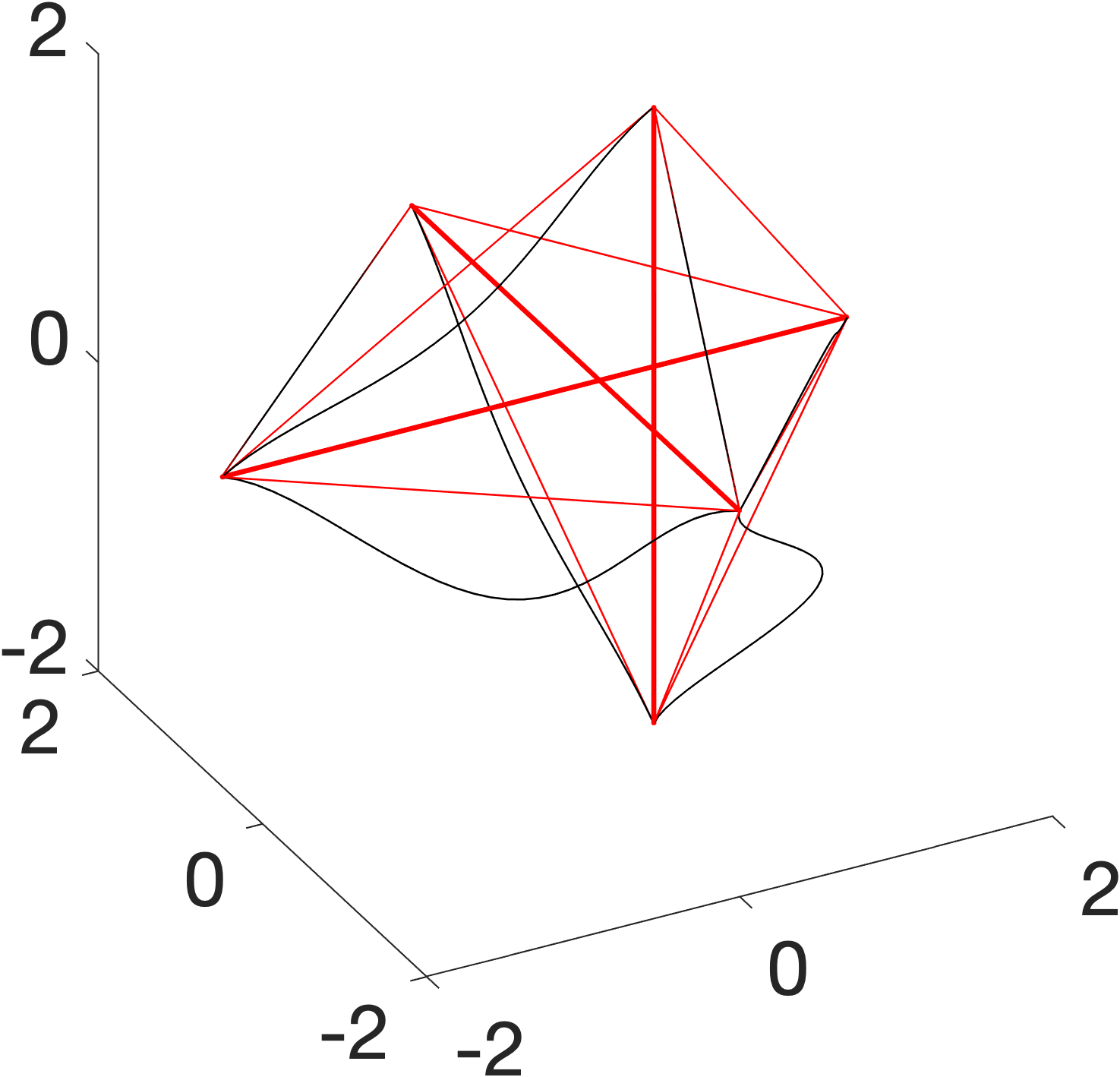}
\end{subfigure}
\begin{subfigure}{.24\textwidth}
  \centering
  \includegraphics[width=0.9\textwidth]{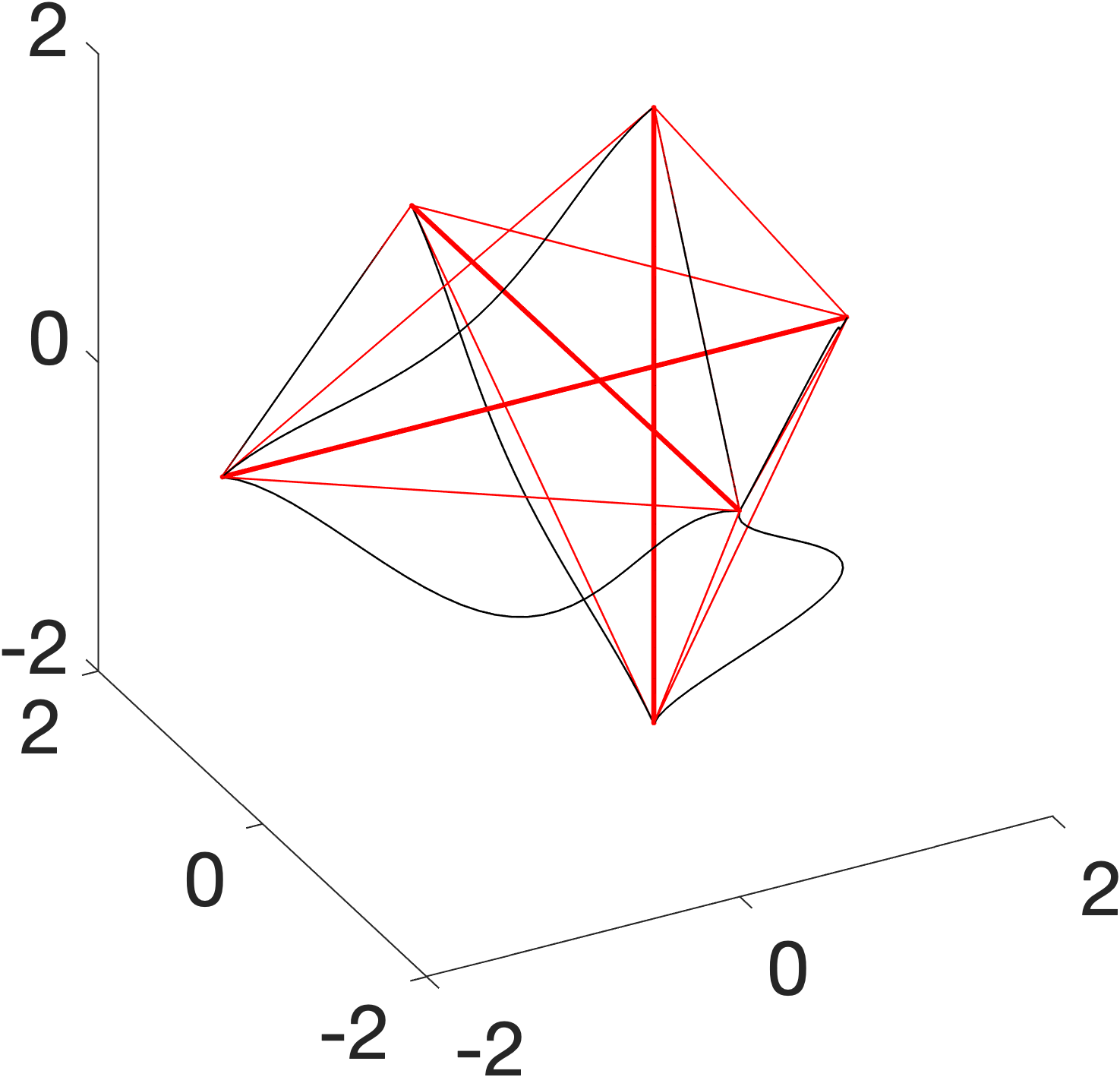}
\end{subfigure}
\begin{subfigure}{.24\textwidth}
  \centering
  \includegraphics[width=0.9\textwidth]{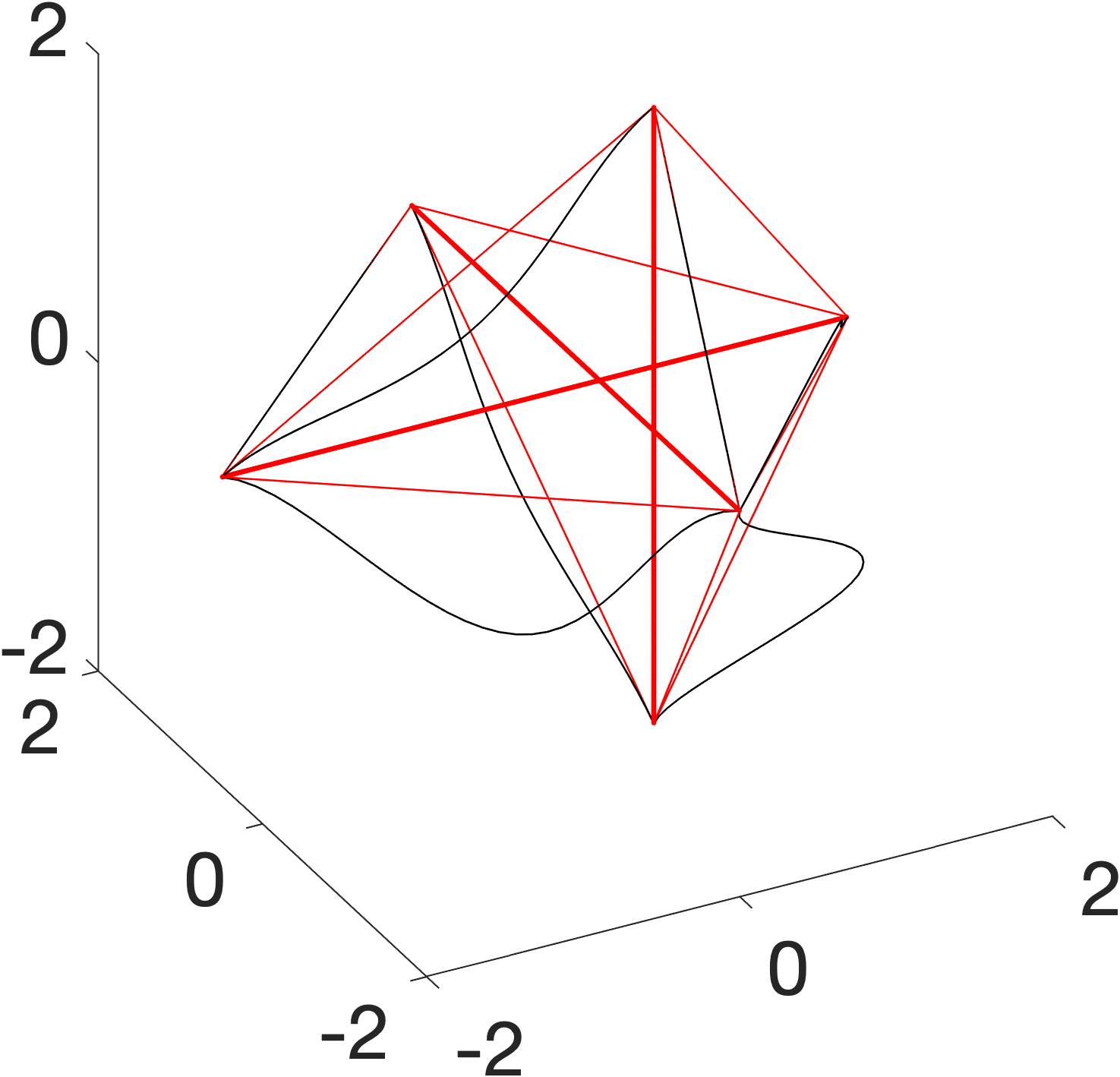}
\end{subfigure}
\caption{Visualization of the lowest mode of deformation of the 10-segrity at 7.1414 rad/s. One end of each rod is fixed.\label{Fig:10segritySwing}}
\end{figure*}

In order to visualize the swinging mode better, one end of each of the rods are fixed. This is equivalent to fixing one end of the tensegrity on the ground or to another structure. The lowest mode for the 9-segrity and the 10-segrity are visualized in \Cref{Fig:9segritySwing} and \Cref{Fig:10segritySwing} respectively. The videos of these visualizations are also enclosed along with this paper. The 9-segrity shows a twist in its top surface with the rods moving apart. As discussed in the earlier section, no 10-th string can be added to this structure that can prevent this torsional motion. In contrast, the proposed 10-segrity demonstrates are more stable behavior.

As discussed in the previous sub-section, it is again important to note here that once a string is slack, they can occupy several positions. The positions for the strings shown in higher modes would be only one representation of the several finite shapes that the string can attain. However, the primary noteworthy point again is that such weaker configurations are possible easily in the lower-energy modes in the existing 9-string rather than the 10-string structures.

% Section 7 (Conclusions and future outlook)
\section{Conclusions and future work}
This work presents a mathematical framework for the design of tensegrities. The famous three-rod tensegrity, used here, shows a vividly clear torsional mode of vibration, often also referred to a swinging soft mode. The proposed form-finding strategy facilitates to find a three-rod and ten-string structure that are stable to the second-order and suitable for engineering.

The proposed form-finding methodology uses mechanics as a foundation to ensure that the resulting designs are stable and free of swinging modes. The work provides a simple to use relation between the number of rods, beads and strings to check for the stability of the resulting structure. This has been demonstrated through the famous three-rod tensegrity. Further on, this work also provides a comprehensive review to show that most of the designs produced till date do not satisfy this relation and thus, significant room for improvements exists. Such a mathematical framework has the potential to be further extended to introduce swinging modes to leverage the instabilities for engineering purposes as well.

One of the areas for the future includes consideration of contact behavior at the joints of rods and strings. This work uses beam models and considers a common node at the ends. Joining the strings through a common node, instead of joint constraint, is expected to allow the bending mode of deformations to be transferred to the strings as well. This has been addressed in this work through the usage of no-compression constitutive model for the strings. Any bending deformation will necessarily introduce compression but however such compressive loads cannot be sustained by the strings and thus non-physical and will not be permitted. However, a string element as introduced in the work of \citet{GiFiTa2017} would be more appropriate alongside usage of appropriate joints and will be considered in the future work. Further on, contacts at joints can introduce additional non-linearities that warrants additional future investigation. Some of the recent works \cite{bach2021, bachwi2022} in this area demonstrate novel modeling architectures that can facilitate modeling of such lattice structures by considering the role of joints. 

One of the areas for application of such stable tensegrity towers is as deployable space structures. Deployable space structures can help save space and weight, that form the primary constraints for space engineering. While mechanisms have been explored for deployment of these structures, extremely fine tuning is required to actually maintain them if unstable configurations were used. The proposed mathematical framework can be used to design and develop stable tensegrity designs that can be used in deployable space structures.

In the recent years, lattice-based metamaterials have been proposed as viable alternatives for impact absorbers or to selectively shield certain frequencies. With the advent of 3D printing, tensegrity structures will find many potential applications in these areas. Thus, the above proposed mathematical framework can provide an easy way to design them from ground-up. 

% Acknowledgements
% Abstract to be restricted to 100 words
\section*{Acknowledgements}
I (Ajay) would like to thank Dr. Tim Poston for introducing me to the world of tensegrities, his expertise, advice and discussions in the early course of this study that has been extremely valuable towards completion of this work.

%% References
% % \bibliographystyle{Classes/elsarticle-num}
%  %\bibliographystyle{Classes/elsarticle-num-names}
\bibliographystyle{Classes/elsarticle-harv}
\bibliography{Extras/Reference}
\pagebreak

%\input{Extras/Vitae}

%% appendix sections are then done as normal section
%\appendix
%\input{Extras/Appendix}

\end{document}